\documentclass[11pt]{article}
 \usepackage{benstyle}
 
\usepackage[font={small}]{caption}

\title{Infinite-dimensional $\ell^1$ minimization and function approximation from pointwise data}
\author{Ben Adcock \\ Department of Mathematics \\ Simon Fraser University \\ Canada}

\begin{document}
\maketitle

\begin{abstract}
We consider the problem of approximating a smooth function from finitely-many pointwise samples using $\ell^1$ minimization techniques.  In the first part of this paper, we introduce an infinite-dimensional approach to this problem. Three advantages of this approach are as follows.  First, it provides interpolatory approximations in the absence of noise.  Second, it does not require \textit{a priori} bounds on the expansion tail in order to be implemented.  In particular, the truncation strategy we introduce as part of this framework is independent of the function being approximated, provided the function has sufficient regularity.  Third, it allows one to explain the key role weights play in the minimization; namely, that of regularizing the problem and removing aliasing phenomena.  In the second part of this paper we present a worst-case error analysis for this approach.  We provide a general recipe for analyzing this technique for arbitrary deterministic sets of points.  Finally, we use this tool to show that weighted $\ell^1$ minimization with Jacobi polynomials leads to an optimal method for approximating smooth, one-dimensional functions from scattered data.
\end{abstract}

\section{Introduction}\label{s:introduction}

Many problems in science and engineering require the approximation of a smooth function from a finite set of pointwise samples.  Although a classical problem in approximation theory, in the last several years there has been an increasing focus on the use of convex optimization techniques for this task \cite{DoostanOwhadiSparse,HamptonDoostanCSPCE,MathelinGallivanCSPDErandom,PengHamptonDoostantweighted,RauhutWardSpher,Rauhut,RauhutWardWeighted,YanGuoXui_l1UQ,KarniadakisUQCS}.  This is driven in part by applications such as uncertainty quantification, wherein the dimension is typically high and the amount of data severely limited.  As dimension increases, smooth multivariate functions are increasingly well-represented by their best $k$-term approximation in certain orthogonal expansions (e.g.\ multivariate Legendre polynomials).  Hence the expectation is that these techniques will yield improvements over more standard approaches -- such as discrete least squares and interpolation -- at least when the dimension is sufficiently high and the data points arise from appropriate sampling distributions.  A number of recent studies, such as those listed above, have shown this to be the case.

\subsection{Current approaches}\label{ss:currentapproaches}
Let $f$ be a function in $L^2(D)$, where $D \subseteq \bbR^d$, $\{ \phi_i \}_{i \in \bbN}$ be an orthonormal basis and write $f = \sum_{i \in \bbN} x_i \phi_i$, where $x = \{ x_i \}_{i \in \bbN} \in \ell^2(\bbN)$ is the infinite vector of coefficients of $f$.  If $\{ t_n \}^{N}_{n=1}$ is a finite set of points, the problem is to approximate $x$, and therefore $f$, from the data $\{ f(t_n) \}^{N}_{n=1}$.  

Since $x$ is an infinite vector, in order to compute an approximation to $f$ it is necessary to truncate in some way.  In the usual formulation (see \cite{DoostanOwhadiSparse,HamptonDoostanCSPCE,MathelinGallivanCSPDErandom,PengHamptonDoostantweighted,RauhutWardSpher,Rauhut,RauhutWardWeighted,YanGuoXui_l1UQ,KarniadakisUQCS}), one introduces a fixed $M \geq N$ and seeks to approximate the first $M$ coefficients $x_1,\ldots,x_M$ of $x$.   If $A = \{ \phi_i(t_n) \}^{N,M}_{n=1,i=1} \in \bbC^{N \times M}$ then the standard (weighted) $\ell^1$ minimization problem is as follows:
\be{
\label{intro_bad_min}
\min_{z \in \bbC^M} \| z \|_{1,w}\ \mbox{subject to $\| A z - y \| \leq \delta$},\qquad y = \{ f(t_n) \}^{N}_{n=1}.
}
Here $\| z \|_{1,w} = \sum^{M}_{i=1} w_i | z_i|$ is the $\ell^1_w$-norm on $\bbC^M$ with weights $w_i > 0$.  The parameter $\delta$ handles the truncation, and is normally chosen so that $\{ x_i \}^{M}_{i=1}$ is feasible for \R{intro_bad_min}.  That is,
\be{
\label{finite_bad}
\max_{n=1,\ldots,N} \left | f(t_n) - \sum^{M}_{i=1} x_i \phi_i(t_n) \right | \leq \delta.
}
In other words, the error introduced by truncating the infinite expansion to a vector of length $M$ is viewed as noise in the data.  

Unfortunately, this formulation raises a number of issues, which we describe next.  Overcoming these issues is the goal of this paper.

\pbk
(i) In order to choose $\delta$, one must have an \textit{a priori} estimate for the truncation error $|f(t) - \sum^{M}_{i=1} x_i \phi_i(t) |$.  Note that the approximation error resulting from \R{intro_bad_min} is sensitive to the choice of $\delta$ \cite{KarniadakisUQCS}.  In practice, cross-validation techniques have been proposed to empirically determine the truncation error \cite{DoostanOwhadiSparse,HamptonDoostanCSPCE,KarniadakisUQCS}.  Yet such techniques are time-consuming, largely lack theoretical support and may not always result in an accurate estimation.

\pbk (ii) The approximation $\tilde{f}$ of $f$ obtained from \R{finite_bad} does not interpolate the data.  In the absence of noise, interpolatory solutions are often desirable in applications since they ensure that the approximation exactly fits the underlying function $f$ at the points at which $f$ is known.

\pbk (iii) The approximation $\tilde{f}$ can be dependent on the choice of weights, and is prone to aliasing (also known as overfitting) if the weights are chosen inappropriately \cite{RauhutWardWeighted}.

\pbk (iv) Besides some specific cases where such techniques are known to perform extremely well -- such as when the coefficients $\{ x_i \}^{M}_{i=1}$  are sparse and the data points $\{ t_n \}^{N}_{n=1}$ are chosen randomly according to the orthogonality measure of the basis $\{\phi_i\}_{i \in \bbN}$ \cite{AdcockCSFunInterp,HamptonDoostanCSPCE,RauhutWardSpher,Rauhut,RauhutWardWeighted} -- very little is known about the approximation error $| f(t) - \tilde{f}(t) |$.  In particular, for general $f$ (not necessarily having sparse coefficients) and arbitrary \textit{deterministic} scattered points $\{ t_n \}^{N}_{n=1}$ the quality of the approximation $\tilde{f}$ to $f$ is largely unknown.

\pbk
Issue (iv) has ramifications for a variety of applications where the primary limitation is the availability of data -- that is, where it is time-consuming or expensive to acquire more samples --  as opposed to data-rich scenarios where processing speed is the key concern (in which case classical techniques such as least-squares fitting are likely superior).  If weighted $\ell^1$ minimization techniques are to find wide use in practice, then it is beneficial to have error bounds for both ideal (i.e.\ random sample points) and non-ideal conditions (i.e.\ fixed, deterministic sample points).

\subsection{Our contributions}
The purpose of this paper is to address these issues.  In \S \ref{s:minproblems} we first propose an infinite-dimensional weighted $\ell^1$ minimization problem which removes the need for \textit{a priori} knowledge of magnitude of the expansion tail.  In the absence of noise, its solutions are exactly interpolatory, unlike solutions of \R{intro_bad_min}.    As one might expect, however, such an infinite-dimensional minimization problem cannot be solved numerically.  Hence we next introduce a truncation strategy based on a user-controlled parameter $K \in \bbN$.  This leads to finite-dimensional minimization problem over $\bbC^K$, reminiscent of \R{intro_bad_min} but with a number of key differences.  First, unlike \R{intro_bad_min}, it requires no knowledge of the expansion tail, and second, it retains the interpolatory property of the infinite-dimensional problem.  In \S \ref{s:truncation} we show how to select the parameter $K$ in a manner independent of $f$, whenever $f$ has sufficient regularity, and dependent only on the basis $\{ \phi_i \}_{i \in \bbN}$ and data points $\{ t_n \}^{N}_{n=1}$.

Formulating the minimization problem in an infinite-dimensional setting also allows us to address issue (iii).  In \S \ref{s:need} we first show that unweighted $\ell^1$ minimization is largely unsuitable for the function interpolation from scattered, deterministic data, since it leads to an \textit{aliasing} phenomenon.  Specifically, without weights it is possible for the optimization problem to have infinitely-many solutions which interpolate $f$ at the data points, but which do not approximate $f$ to any accuracy away from these points.  Fortunately, this problem can be completely resolved by the introduction of slowly-growing weights.  In effect, these weights regularize the optimization problem and ensure that such bad solutions of the unweighted problem, while still feasible, are no longer minimizers of the weighted problem.  Through subsequent analysis we quantify how fast the weights need to grow to resolve this phenomenon, and demonstrate this result with numerical examples.

Issues (i)--(iii) are the focus of the first half of this paper (\S \ref{s:prelim}--\ref{s:truncation}).  In the second half, we consider (iv).  More precisely, we pose and answer the following two questions:

\pbk (a) How well can one approximate a function $f$ using weighted $\ell^1$ minimization from its samples taken on an arbitrary deterministic grid of $N$ points?

\pbk (b) How does this approximation perform in comparison to other techniques, such as least squares?

\pbk
Note that we do not assume any sparsity of the coefficients $x = \{ x_i \}_{i \in \bbN}$ of $f$, although we do assume some mild decay of $x_i$ as $i \rightarrow \infty$ (otherwise the weighted $\ell^1$ problem does not make sense).  We also do not assume any structure to the data: the points $\{ t_n \}^{N}_{n=1}$ are deterministic and can be arbitrarily distributed in the domain.  As is standard in scattered data approximation,  we classify the error in terms of its density (or fill distance) \cite{WendlandScatteredBook}.

Our motivation for examining (a) and (b) is the following.  Least-squares fitting is a classical and widely-used technique, but is well known to be intensive in the number of samples required to achieve stability and accuracy (see \cite{DavenportEtAlLeastSquares,MiglioratiNobileLowDisc,MiglioratiEtAlFoCM} and references therein).    Conversely, under certain conditions (sparsity and random sampling) $\ell^1$ techniques are known to give very good approximations from relatively few samples.  However, in certain practical scenarios -- such as when using legacy data -- one may not have the luxury to choose the data points in a way to deliver the best accuracy of the approximation.  Moreover, while functions in high dimensions tend to have sparse coefficients in polynomial bases \cite{CohenDeVoreSchwabFoCM,CohenDeVoreSchwabRegularity,DoostanOwhadiSparse}, in low (in particular, one) dimensions polynomial coefficients usually exhibit rapid decay, but typically little sparsity.  Since $\ell^1$-based techniques are computationally more intensive than classical methods such as least-squares fitting, this raises the following question: is it still worth using $\ell^1$ techniques even when the data points are scattered and sparsity is not assured?    

In \S \ref{s:lin_approx_err} we present a general mathematical framework for answering these questions.  We introduce a linear approximation error analysis for the infinite-dimensional weighted $\ell^1$-minimization, which allows it to be compared  directly to existing techniques.  In particular, we reduce (b) to a question about the behaviour of three particular quantities that depend on the data points $\{ t_n \}^{N}_{n=1}$, the expansion basis $\{ \phi_i \}_{i \in \bbN}$ and the weights $\{ w_i \}_{i \in \bbN}$.  Analyzing these quantities for each specific problem setup provides an answer to (b).

To illustrate the various aspects of this framework, in the final part of this paper (\S \ref{s:AlgPoly_Examp} and \S \ref{s:TrigPoly_Examp}) we consider several examples, including one-dimensional Jacobi polynomial approximations from scattered data points.  In this case, we prove the following:

\thm{
\label{t:intro_thm}
For $\alpha,\beta > -1$ let $\{ \phi_i \}_{i \in \bbN}$ be the orthonormal Jacobi polynomial basis \R{Jacobi_ON} on $[-1,1]$ and let $T = \{ t_n \}^{N}_{n=1}$ be a set of $N$ scattered points in $[-1,1]$.  Let $h$ be the density of the points, defined by \R{h_def}, and suppose that the truncation parameter 
\be{
\label{trunc_choice}
K \gtrsim h^{\frac{1}{2r}} \xi^{-1-\frac{1}{r}},
}
for some $r \in \bbN$, where $\xi$ is the minimal separation between the points $T$.  Fix weights $w = \{ w_i \}_{i \in \bbN}$ with  $w_{i} = \| \phi_i \|_{L^\infty} i^{\gamma}$ for some $\gamma > \max\{1/2 - q,0\}$, where $q$ is as in \R{Jacobi_ON_growth}, and let $f = \sum_{i \in \bbN} x_i \phi_i$ with $x \in \ell^1_{\tilde{w}}(\bbN)$ for $\tilde{w}_i = \sqrt{i} (w_i)^2$.  Then given measurements $y = \{ f(t_n) \}^{N}_{n=1}$ one can compute, via weighted $\ell^1$ minimization with weights $w$, an approximation $\hat{x}$ to the coefficients $x$ satisfying
\be{
\label{intro_err}
\| x - \hat{x} \| \lesssim \| x - P_M x \|_{1,w} + \| x - P_K x \|_{1,\tilde{w}},
}
where $P_K x = \{ x_1,\ldots,x_K,0,0,\ldots \}$ and $P_M x = \{x_1,\ldots,x_M,0,0,\ldots\}$, provided
\be{
\label{intro_h_scale}
h^{-1} \gtrsim M^2 \log M.
}
Moreover, the approximation $\tilde{f} = \sum^{K}_{i=1} \hat{x}_i \phi_i$ exactly interpolates the data: $\tilde{f}(t_n) = f(t_n)$, $\forall n$.
}

This theorem demonstrates the key aspects of this paper.  (i): the truncation parameter $K$ is determined independently of $f$, and its contribution to the overall error is clarified by \R{intro_err}.  In particular, if the data is roughly equally-spaced, then $h, \xi = \ord{1/N}$ and it suffices to take $K = \ordu{N^{1+\frac{1}{2r}}}$ for any $r>0$.  (ii): the approximation $\tilde{f}$ exactly interpolates $f$ in the absence of noise.  Note that noise can also be dealt with within our framework; we exclude it here for ease of presentation.  (iii): one gets an explicit criterion for how to choose the weights.  (iv): the estimate \R{intro_err} for the approximation error depends only on the density $h$ of the deterministic points $T$ which can be arbitrarily distributed in the domain.

As we discuss in \S \ref{s:AlgPoly_Examp}, the estimates \R{intro_err} and \R{intro_h_scale} demonstrate not just good performance of this approach for scattered data, but in fact near-optimal performance.  As we explain, no stable method which is convergent as $h \rightarrow 0$ can exhibit an error bound depending on $x - P_M x$ measured in some norm with $M$ growing faster than $h^{-1/2}$ as $h \rightarrow \infty$.  Our numerical results support this conclusion, and in fact show that weighted $\ell^1$ minimization performs rather better in practice and similarly to an oracle least-squares fit.

\subsection{Relation to previous work}\label{ss:relation}
A theory for reconstruction of sparse polynomials from random pointwise samples was developed in a series of papers by Rauhut \& Ward \cite{Rauhut,RauhutWardWeighted}.  Extension and application of this work in uncertainty quantification has been considered in \cite{AdcockCSFunInterp,ChkifaDownwardsCS,DoostanOwhadiSparse,HamptonDoostanCSPCE,MathelinGallivanCSPDErandom,PengHamptonDoostantweighted,YanGuoXui_l1UQ,KarniadakisUQCS}.  The use of weighted $\ell^1$ minimization was introduced in \cite{PengHamptonDoostantweighted,RauhutWardWeighted,KarniadakisUQCS}.  We also use weighted minimization in this paper for approximation from determinstic samples, yet for rather different purposes.  Namely, weights are chosen to regularize the minimization problem and remove the aliasing phenomenon.  Typically, this requires only very slow growth of the weights, which we quantify in the paper.  Unlike other works, we do not select weights based on \textit{a priori} information about decay of the polynomial coefficients.  In fact, in \S \ref{s:need} we will show that choosing weights in this way leads to inconsistent and often negligible improvements when the samples are scattered and deterministic.  

The infinite-dimensional framework we introduce in this paper is inspired in part by the framework of infinite-dimensional compressed sensing in Hilbert spaces, due to A.\ C.\ Hansen and the present author \cite{BAACHShannon,BAACHGSCS,BAACHOptimality} (see also \cite{BAGSAIEP} for an overview).  A key difference is the need for weighted minimization in the present setup, due to the lack of continuity of the sampling operator.  We note also that our worst case analysis and comparison to least-squares fitting is similar to that presented in \cite{GSl1} for generalized sampling in the Hilbert space setting.

The examples we use in this paper consist of algebraic and trigonometric polynomials respectively.  Polynomial approximations (so-called polynomial chaos expansions) are popular in areas such as uncertainty quantification \cite{SMUQ,XiuKarniadakisPC}.  However, we stress that the framework and analysis of \S \ref{s:prelim}--\ref{s:lin_approx_err} of this paper is completely general, and can be applied to other bases.  We mention several other examples in \S \ref{s:conclusions}.  Our examples are also one-dimensional.  We do this so as to better elucidate the key ideas, without the notational complexities of the multivariate setting.

On this topic, we wish to clarify that the aim of this paper is not to propose weighted $\ell^1$ minimization as a panacea for function approximation.  In the one-dimensional setting especially there is a wealth of other techniques which are likely superior (see \cite{FEStability,AdcockPlatteMapped,BoydRunge,TrefPlatteIllCond} and references therein).  The advantages of weighted $\ell^1$ minimization come to the fore as the dimension increases, as has been verified empirically in a number of works such as those mentioned previously.  Instead, the purpose of this paper is to first propose a framework for weighted $\ell^1$ minimization that overcomes some existing issues, and second provide a more comprehensive analysis of its approximation capabilities for fixed samples.  We use the one-dimensional case to this end primarily for illustrative purposes.

\section{Preliminaries}\label{s:prelim}
Let $D \subseteq \bbR^d$ be a domain and $\nu(t)$ an integrable nonnegative weight function satisfying $\int_D \nu(t) \D t = 1$.  Let $L^2_{\nu}(D)$ be the space of complex-valued weighted square-integrable functions on $D$ with norm $\nm{\cdot}_{L^2_{\nu}}$ and inner product $\ip{\cdot}{\cdot}_{L^2_{\nu}}$, and suppose that $\{ \phi_i \}_{i \in \bbN} \subseteq L^2_{\nu}(D) \cap L^\infty(D)$ is a set of functions that are orthonormal with respect to $\nu$.  Note that
\be{
\label{phii_unif}
1 = \| \phi_i \|^2_{L^2_{\nu}} \leq \| \phi_i \|^2_{L^{\infty}},\qquad \forall i \in \bbN,
}
where $\nm{\cdot}_{L^{\infty}}$ is the uniform norm on $D$.

\subsection{Scattered data}
For $N \in \bbN$, let $T  = \{ t_n \}^{N}_{n=1} \subseteq \overline{D}$ be a set of $N$ scattered data points.  Our aim is to approximate a function $f : D \rightarrow \bbC$ from the values $\{ f(t_{n}) \}^{N}_{n=1}$.  To ensure an accurate approximation, we require a notion of closeness of the points $T$.  We quantify this by defining the \textit{density}
\be{
\label{h_def}
h = \sup_{t \in D} \min_{n=1,\ldots,N} | t - t_{n} |,
}
(also known as the fill distance \cite{WendlandScatteredBook}) where $\abs{\cdot}$ is the Euclidean distance.  In our analysis later, we will present convergence rates of the various approximations in terms of $h \rightarrow 0$.

Associated to the points $T$ will also be a set of values $\tau_{n} \geq 0$, $n = 1,\ldots,N$, which we refer to as \textit{quadrature weights}.  This is not to be confused with the optimization weights $w_i$ introduced later.  For simplicity, we define these as follows
\be{
\label{quad_weights_1}
\tau_{n} = \int_{V_{n}} \nu(t) \D t,\quad n=1,\ldots,N,\qquad V_{n} = \left \{ t \in D : | t - t_n | \leq | t - t_m | , \ \forall m \neq n \right \},
}
where $V_{n}$ are the Voronoi cells of the points $T$ in $D$.  Given such quadrature weights, we define the following sesquilinear form on $L^2_{\nu}(D) \cap L^\infty(D)$:
\bes{
\ip{f}{g}_h = \sum^{N}_{n=1} \tau_{n} f(t_{n}) \overline{g(t_{n}) },
}
and write $\nm{\cdot}_h = \sqrt{\ip{\cdot}{\cdot}_h}$ for the corresponding seminorm.  Note that the quadrature weights $\tau_n$ are not strictly necessary at this stage, but will play a pivotal role later in the paper.

\subsection{Weighted spaces}\label{ss:weighted_spaces}
For the remainder of this paper, $w = \{ w_i \}_{i \in \bbN}$ will be a set of positive weights satisfying
\be{
\label{weights_growth}
w_i \geq \| \phi_i \|_{\infty} \geq 1,\quad \forall i \in \bbN,
} 
where the latter inequality is due to \R{phii_unif}.  Define the weighted $\ell^p$ spaces by
\bes{
\ell^p_w(\bbN) = \left \{ x = \{ x_i \}_{i \in \bbN} : \ \| x \|_{p,w} : = \left ( \sum_{i \in \bbN} (w_i)^p | x_i |^p \right )^{1/p} < \infty \right \},\quad p > 0.
}
Note that $\| x \|_{p,w} = \| W x \|_{p}$, where 
\be{
\label{W_def}
W = \mathrm{diag}(w_1,w_2,\ldots),
}
is the infinite diagonal matrix of weights.  For the remainder of this paper, we will assume that the function we wish recover $f = \sum_{i \in \bbN} x_i \phi_i \in L^2_{\nu}(D)$ has coefficients $x = \{x_i \}_{i \in \bbN} \in \ell^1_w(\bbN)$.

\subsection{Other notation}
For $\Delta \subseteq \bbN$ we let $P_{\Delta} : \ell^2(\bbN) \rightarrow \ell^2(\bbN)$ be the projection defined by
\bes{
(P_{\Delta} x)_j = x_j,\quad j \in \Delta,\qquad (P_{\Delta} x)_j =0,\quad j \notin \Delta.
}
If $\Delta = \{1,\ldots,K \}$ for some $K \in \bbN$, then we merely write $P_K$.  We also let $\{ e_j \}_{j \in \bbN}$ denote the canonical basis of $\ell^2(\bbN)$, so that
\bes{
P_{\Delta}(\cdot) = \sum_{j \in \Delta} \ip{\cdot}{e_j} e_j.
}
We will allow the slight abuse of notation throughout the paper in thinking of $P_{\Delta} x$ as both an element of $\ell^2(\bbN)$ and $\bbC^{|\Delta|}$.  The intended meaning will be clear from the context.

If $x \in \bbC$, we let $\mathrm{sign}(x) = x / | x |$ be its complex sign with the convention that $\mathrm{sign}(0) =0$.  For $x \in \ell^\infty(\bbN)$ we let $\mathrm{sign}(x) = \{ \sgn(x_i) \}_{i \in \bbN} \in \ell^\infty(\bbN)$ be the corresponding sequence of complex signs of the entries of $x$.  Finally, we use the notation $a \lesssim b$ to mean that there exists a constant $C$ independent of all relevant parameters such that $a \leq C b$.

\subsection{Examples}\label{ss:examples}
As mentioned in \S \ref{ss:relation}, the examples we consider in this paper consist of one-dimensional functions on bounded intervals, which we take to be $D=(-1,1)$ without loss of generality.  In \S \ref{s:conclusions} we briefly discuss extensions to higher dimensions, unbounded intervals and other approximation systems.

\examp{
\label{ex:Jacobi}
If $f$ is smooth, then it is natural to approximate it using a basis of orthogonal polynomials.  Let
\bes{
\nu(t) = \nu^{(\alpha,\beta)}(t) = c^{(\alpha,\beta)}(1-t)^\alpha(1+t)^{\beta},\qquad \alpha,\beta > -1,
}
be the Jacobi weight function, where $c^{(\alpha,\beta)} = \left(  \int^{1}_{-1} (1-t)^\alpha(1+t)^{\beta} \D t \right )^{-1}$ is a normalizing constant, $P^{(\alpha,\beta)}_j$ be the $j^{\rth}$ orthogonal polynomial with respect to this weight function, and
\be{
\label{Jacobi_ON}
\phi_j = \left ( \kappa^{(\alpha,\beta)}_{j-1} \right )^{-1/2} P^{\alpha,\beta}_{j-1},
}
be the corresponding orthonormal polynomial, where $\kappa^{(\alpha,\beta)}_j$ is as in \R{Jacnorm_def}.
One can show that
\be{
\label{Jacobi_ON_growth}
\| \phi_j \|_{L^\infty} = \ord{j^{q+1/2}},\quad j \rightarrow \infty,\qquad \mbox{where $q = \max \{ \alpha,\beta , -1/2 \}$}.
}
See Appendix \ref{a:Jacobi} (several other properties of Jacobi polynomials that will be needed later are also listed therein).  Since the weights $\{ w_i \}_{i \in \bbN}$ introduced in \S \ref{ss:weighted_spaces} are required to satisfy \R{weights_growth}, this means that for this example they must grow at least as fast as $j^{q+1/2}$ as $j \rightarrow \infty$.
}

\examp{
\label{ex:Trigonometric}
Functions that are smooth and periodic can be efficiently approximated using trigonometric polynomials.  In this case, we have $\nu(t) = 1/2$ and define $\{ \phi_i \}_{i \in \bbN}$ to be the Fourier basis
\be{
\label{fourier_basis}
\phi_j(t) =\E^{\I j \pi t},\quad j \in \bbZ.
}
For convenience we index over $\bbZ$ rather than $\bbN$ in this example.  Note that $\| \phi_j \|_{\infty}=1$ and therefore the weights  $w_j$ in this example are required to satisfy $w_j \geq 1$, $\forall j \in \bbZ$.
}

\section{Minimization problems}\label{s:minproblems}
Define the operator $U : \ell^1_{w}(\bbN) \rightarrow \bbC^N$ by $U x = \{ \sqrt{\tau_{n}} g(t_{n}) \}^{N}_{n=1}$ where $g = \sum_{i \in \bbN} x_i \phi_i$.  Note that this operator is bounded.  
We shall also view $U \in \bbC^{N \times \infty}$ as the infinite matrix with entries
\bes{
U_{n,i} = \sqrt{\tau_{n}} \phi_i(t_{n}),\quad n=1,\ldots,N,\ i \in \bbN.
}
From now on, we make no distinction between the operator $U$ and the infinite matrix.

\subsection{Infinite-dimensional weighted $\ell^1$ minimization}
Let $f = \sum_{i \in \bbN} x_i \phi_i$ be a function we wish to recover, where $x = \{ x_i \}_{i \in \bbN} \in \ell^1_w(\bbN)$.  Suppose first that we are given noiseless measurements of $f$, that is, $f(t_{n})$, $n=1,\ldots,N$, and let 
\bes{
y = \left \{  \sqrt{\tau_{n}} f(t_{n}) \right \}^{N}_{n=1} \in \bbC^N,
}
be the vector of measurements normalized by the quadrature weights.  To recover the infinite vector $x$ of coefficients, and therefore $f$, we shall use weighted $\ell^1$ minimization.   In order to avoid issues of truncation (recall \S \ref{ss:currentapproaches}), we first formulate the following infinite-dimensional optimization problem:
\be{
\label{inf_min_noiseless}
\inf_{ z \in \ell^1_w(\bbN)} \| z \|_{1,w}\  \mbox{subject to $U  z = y$}.
}
If $\hat{x} \in \ell^1_w(\bbN)$ is a minimizer of \R{inf_min_noiseless}, then the corresponding approximation $\tilde{f}$ to $f$ is given by
\be{
\label{f_N_inf}
\tilde{f} = \sum_{i \in \bbN} \hat{x}_i \phi_i.
}
In general, the measurements may be noisy.  Suppose we are given 
\bes{
f(t_n) + e_{n},\quad n=1,\ldots,N,
}
where $|e_n| \leq \eta$, $n=1,\ldots,N$, for some known $\eta \geq 0 $.  Write
\be{
\label{y_N_def}
y = \{ \sqrt{\tau_n} \left ( f(t_n) + e_n \right ) \}^{N}_{n=1} \in \bbC^N.
}
In this case, we solve the inequality-constrained optimization problem
\be{
\label{inf_min}
\inf_{ z \in \ell^1_w(\bbN)} \| z \|_{1,w}\  \mbox{subject to $\| U  z - y \| \leq \eta$}.
}
Note that \R{inf_min_noiseless} is just a special case of \R{inf_min} corresponding to $\eta = 0$, and that both \R{inf_min_noiseless} and \R{inf_min} always have a solution, since the feasible set nonempty (specifically, $x$ is always feasible).  Note also that solutions of \R{inf_min_noiseless} are interpolatory in the sense that $\tilde{f}(t_n) = f(t_n)$, $n=1,\ldots,N$, whenever $\tilde{f}$ is given by \R{f_N_inf} with $\hat{x}$ being a minimizer of \R{inf_min_noiseless}.  Conversely, solutions of \R{inf_min} yield approximations $\tilde{f}$ that are interpolatory up to the noise magnitude $\eta$.

Throughout, we shall assume that the noise bound $\eta$ is known.  If $\eta$ is unknown, one may still solve the equality-constrained problem \R{inf_min_noiseless}  in practice (or the inequality-constrained problem \R{inf_min} with some estimate of the noise).  However, there are no known recovery guarantees for this problem.  See \cite[Chpt.\ 11]{FoucartRauhutCSbook} for some work in this direction in the context of finite-dimensional compressed sensing with random Gaussian matrices.

\subsection{Truncation}
Unfortunately, neither problem \R{inf_min_noiseless} or \R{inf_min} is numerically solvable, since they require optimizing over an infinite-dimensional space.  Let $K \in \bbN$ be a truncation parameter.  To form a computable problem, we replace the space $\ell^1_{w}(\bbN)$ with $\bbC^K$ and truncate the $N \times \infty$ matrix $U$ to the $N \times K$ matrix $U P_K$ spanned by its first $K$ columns.  Hence, we now consider the problem
\be{
\label{fin_min_noiseless}
\min_{ z \in \bbC^K} \| z \|_{1,w}\  \mbox{subject to $U P_{K} z = y $},
}
in the noiseless case, as well as its noisy analogue
\be{
\label{fin_min}
\min_{ z \in \bbC^K} \| z \|_{1,w}\  \mbox{subject to $\| U P_{K} z - y \| \leq  \eta $}.
}
Both of these problems are finite dimensional, and can be solved using standard algorithms.  If $\hat{x} \in \bbC^K$ is a minimizer of either, then the approximation to $f$ is given by
\be{
\label{f_N_K}
\tilde{f} = \sum^{K}_{i=1} \hat{x}_i \phi_i.
}
Note that neither \R{fin_min_noiseless} nor \R{fin_min} modify the constraints of the infinite-dimensional problems \R{inf_min_noiseless} and \R{inf_min}.  In particular, \R{fin_min_noiseless} remains interpolatory and \R{inf_min} is interpolatory up to the noise.  

With this in hand, the general idea is to choose $K$ in such a way to ensure closeness of the solutions of the finite-dimensional problems \R{fin_min_noiseless} and \R{fin_min} to those of the infinite-dimensional problems \R{inf_min_noiseless} and \R{inf_min}.  In \S \ref{s:truncation} we shall show that it is possible to choose $K$ in a function-independent manner, thus overcoming issue (i) of \S \ref{ss:currentapproaches}.

\rem{
It is important that $K$ be chosen suitably large.  To see why, consider Example \ref{ex:Jacobi}.  If $K=N$ then \R{fin_min_noiseless} has a unique solution and $\tilde{f}$ is just the polynomial interpolant of $f$ of degree $N-1$.  However, for equispaced (or more generally, scattered) data this is well known to be a poor approximation to $f$, since it suffers from Runge's phenomenon.  The approximations $\tilde{f}$ will generally diverge and the matrix $U P_N$ will have an exponentially-large condition number.  On the other hand, if one replaces $U P_N$ by $U P_K$ with $K > N$, then provided $K$ is sufficiently large the singular values of $UP_K$ are provably bounded away from zero (see \S \ref{s:truncation}).  
}

\subsection{Comparison to least-squares fitting}\label{ss:LS_comp}
Classical least-squares fitting corresponds to the approximation $\tilde{f} = \sum^{M}_{i=1} \check{x}_i \phi_i$,
where $\check{x}$ is the solution of the problem
\be{
\label{LS_fit}
\min_{z \in \bbC^M} \| U P_M z - y \|.
}
An important difference between least-squares fitting and weighted $\ell^1$ minimization is the choice of the truncation parameters.  In the former, the parameter $M$ affects both the approximation error $\| f - \tilde{f} \|$ and the robustness of the approximation.  In practice, $M$ must be chosen suitably small in relation to $1/h$ to ensure a stability and robustness, while also being sufficiently large to give a good approximation.  The issue of how to best choose $M$, which we discuss further in \S \ref{ss:LS_fit_comp} for the specific case of polynomials, is nontrivial.  While there are many known theoretical estimates for how $M$ should scale for different function systems and datasets (see, for example,  \cite{AdcockPlatteMapped,AdcockNecSamp,DavenportEtAlLeastSquares,DemanetTownsendExtrap,MiglioratiNobileLowDisc,MiglioratiEtAlFoCM} and references therein), optimal selection of $M$ is difficult in practice.  In particular, standard theoretical guarantees usually only determine the asymptotic order of $M$ with $1/h$.  Constants, if known, tend to be overly pessimistic.  This problem becomes more acute in multiple dimensions, since the ordering of the basis functions plays an increasingly important role.

Conversely, the truncation parameter $K$ in the weighted $\ell^1$ minimization formulations \R{fin_min_noiseless} and \R{fin_min} plays a completely different role: namely, it allows one to approximately compute solutions of the infinite-dimensional problems \R{inf_min_noiseless} and \R{inf_min}.  Once $K$ is large enough so that the truncation error when passing to the finite-dimensional problems \R{fin_min_noiseless} and \R{fin_min} is negligible, changing $K$ has little effect on the accuracy of the solution $\tilde{f}$.    
Note also that \R{fin_min_noiseless} leads to interpolatory approximations, which is not the case for the least-squares fit \R{LS_fit}.

\section{The need for weights}\label{s:need}
Before analyzing \R{inf_min_noiseless} and \R{fin_min} in detail, we first examine the role weights play in the minimization.  In particular, we shall show that it is in general necessary for the ratios
\be{
\label{weights_growth2}
w_{i} / \| \phi_i \|_{\infty} \rightarrow \infty,\quad i \rightarrow \infty,
}
in order for the weighted minimization problems to give convergent approximations to $f$ in the case of fixed, deterministic data.  If this is not the case, then the minimization problem can have multiple solutions which aliase the data, leading in general to poor approximations.  

In \S \ref{s:lin_approx_err} we shall prove that \R{weights_growth2} is sufficient to guarantee a good approximation.  To demonstrate its general necessity, we consider Example \ref{ex:Trigonometric}.  Recall that $\| \phi_j \|_{\infty}=1$ in this case.

\prop{
\label{p:aliasing}
Let $D$, $\nu$ and $\{ \phi_j \}_{j \in \bbZ}$ be as in \R{fourier_basis}.  Let $T = \{ t_n \}^{N}_{n=1}$ be a set of $N$ data points such that $t_n P \in \bbZ$ for some $P \in \bbN$ and all $n=1,\ldots,N$.  Suppose that $\hat{x} \in \ell^1(\bbZ)$ is a solution of
\be{
\label{bad_prob}
\inf_{z \in \ell^1(\bbZ)} \| z \|_{1}\ \mbox{subject to $U z = y$},
}
where $U = \{ \phi_j(t_n) \}^{N,\infty}_{n=1,j=-\infty}$ and $y \in \bbC^N$.  Then every shift of the entries of $\hat{x}$ by a multiple of $2P$ is also a solution of \R{bad_prob}.  That is, for every $k \in \bbZ$, the element $z \in \ell^1(\bbZ)$ given by
\be{
\label{aliased_solutions}
z_{i} = \hat{x}_{i-2kP},\quad i \in \bbZ,
}
is a solution of \R{bad_prob}.
}
\prf{
Shifting the entries of $x$ does not affect its $\ell^1$ norm, therefore $\| z \|_1 = \| x \|_1$.  Moreover,
\bes{
(U z)_n = \sum_{j \in \bbZ} z_j \E^{\I j \pi t_n} = \sum_{j \in \bbZ} x_j \E^{\I (j+2k P) \pi t_n} = \sum_{j \in \bbZ} x_j \E^{\I j \pi t_n} = (U x)_n = y_n,\quad n=1,\ldots,N.
}
Hence $z$ is feasible for \R{bad_prob}, and therefore a minimizer.
}
The absence of quadrature weights $\tau_n$ does not change the conclusion here, since we consider the equality-constrained minimization.  We could also consider the inequality-constrained problem with much the same result, but we present the equality-constrained problem to show that the phenomenon is not due in any way to the increased size of the feasible set when $\eta > 0$.

Taken on its own, the fact that \R{bad_prob} has multiple solutions may not be alarming.  After all, convex optimization problems often do.  However, in this case the effect is catastrophic.  Consider the problem where $f =  \phi_0 = 1$ so that its coefficients are $x = e_0$.  Since $y_n = f(t_n) = 1$ in this case, if $z \in \ell^1(\bbZ)$ is feasible for \R{bad_prob} then
\bes{
1 = |(U z)_n| = \left | \sum_{j \in \bbZ} z_j \E^{\I j \pi t_n } \right | \leq \| z \|_1.
}
Hence $x$ itself is a solution of \R{bad_prob}, and by Proposition \ref{p:aliasing} so is every shift $z = e_{2kP}$ of $x$ by a multiple of $2P$.  However, for all these solutions one has $\| x - z \| = 2$.  Thus, although there is one solution of \R{bad_prob} which recovers $x$ (and therefore $f$) exactly there are also infinitely many solutions of \R{bad_prob} that give meaningless approximations to $x$.

\begin{figure}[t]
\begin{center}
$\begin{array}{ccc}
\includegraphics[width=5.00cm]{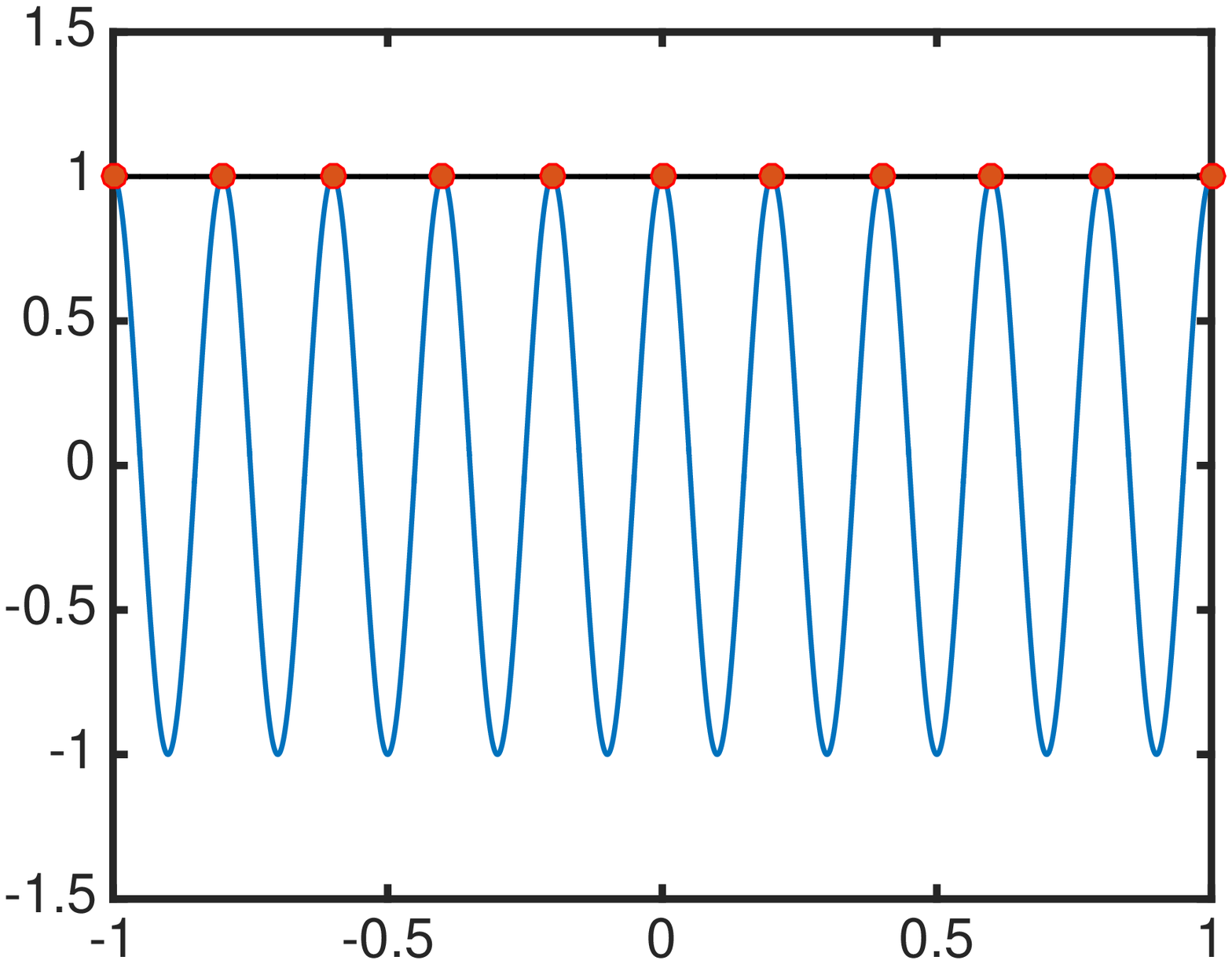}  &  \hspace{2pc} &\includegraphics[width=5.00cm]{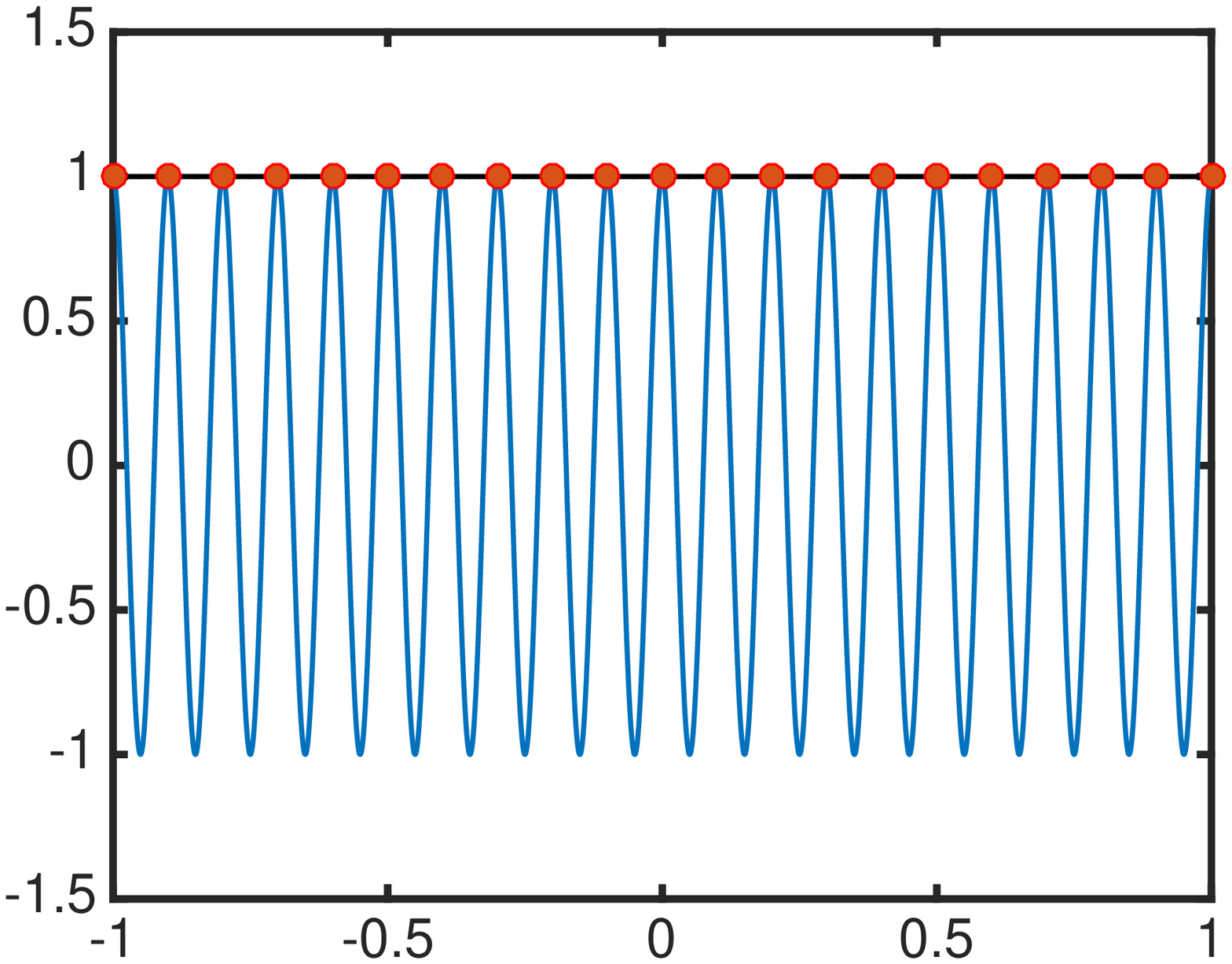} 
\end{array}$
\caption{Aliasing in $\ell^1$ minimization.  The black line is the function $f(t) = \phi_0 = 1$, the red dots are the data points with $N=11$ (left) and $N=21$ (right) and the blue curve is the aliased solution $\phi_{10}$ (left) and $\phi_{20}$ (right).  Although these solutions interpolate $f$, they do not approximate $f$ in between the data points.} \label{f:FourAlias}
\end{center}
\end{figure}

This effect is due to aliasing the data by higher-frequency Fourier modes.  The shifted solutions of Proposition \ref{p:aliasing} correspond to the functions $\phi_{k P}$, $k \in \bbZ$, which interpolate $f$ at the data points but oscillate with frequency proportional to $k P$ in between the data points (see Fig.\ \ref{f:FourAlias}).  Of course, in the simplified scenario described here the aliasing problem could have been avoided by solving a truncated problem with truncation $K = P$.  However, as discussed, this will not work in the general case when truncation with $K \gg N$ is required in order to control the tail and ensure (in the noiseless case) an interpolatory solution.

\begin{figure}[t]
\begin{center}
$\begin{array}{ccc}
\includegraphics[width=4.5cm]{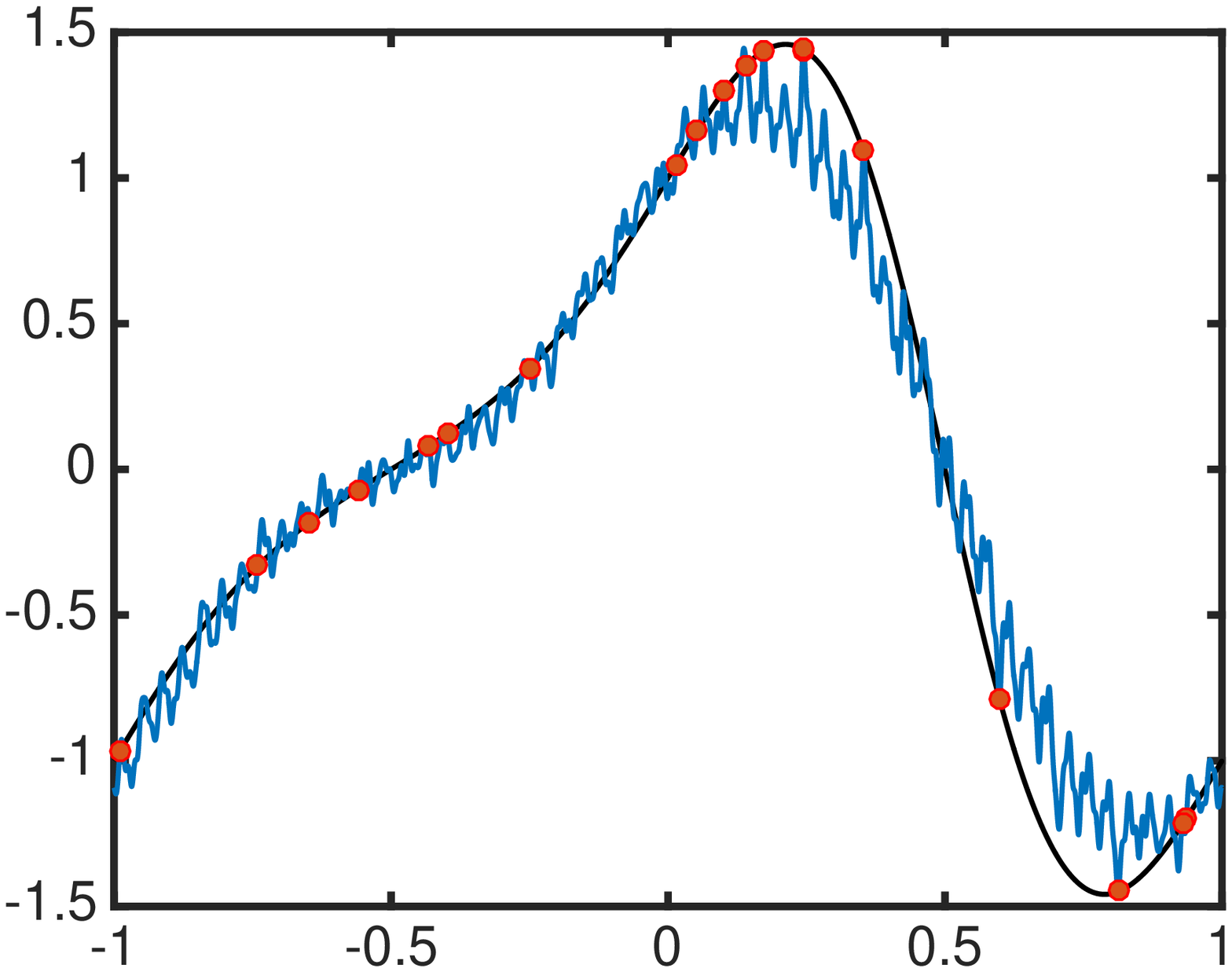}  &  \includegraphics[width=4.5cm]{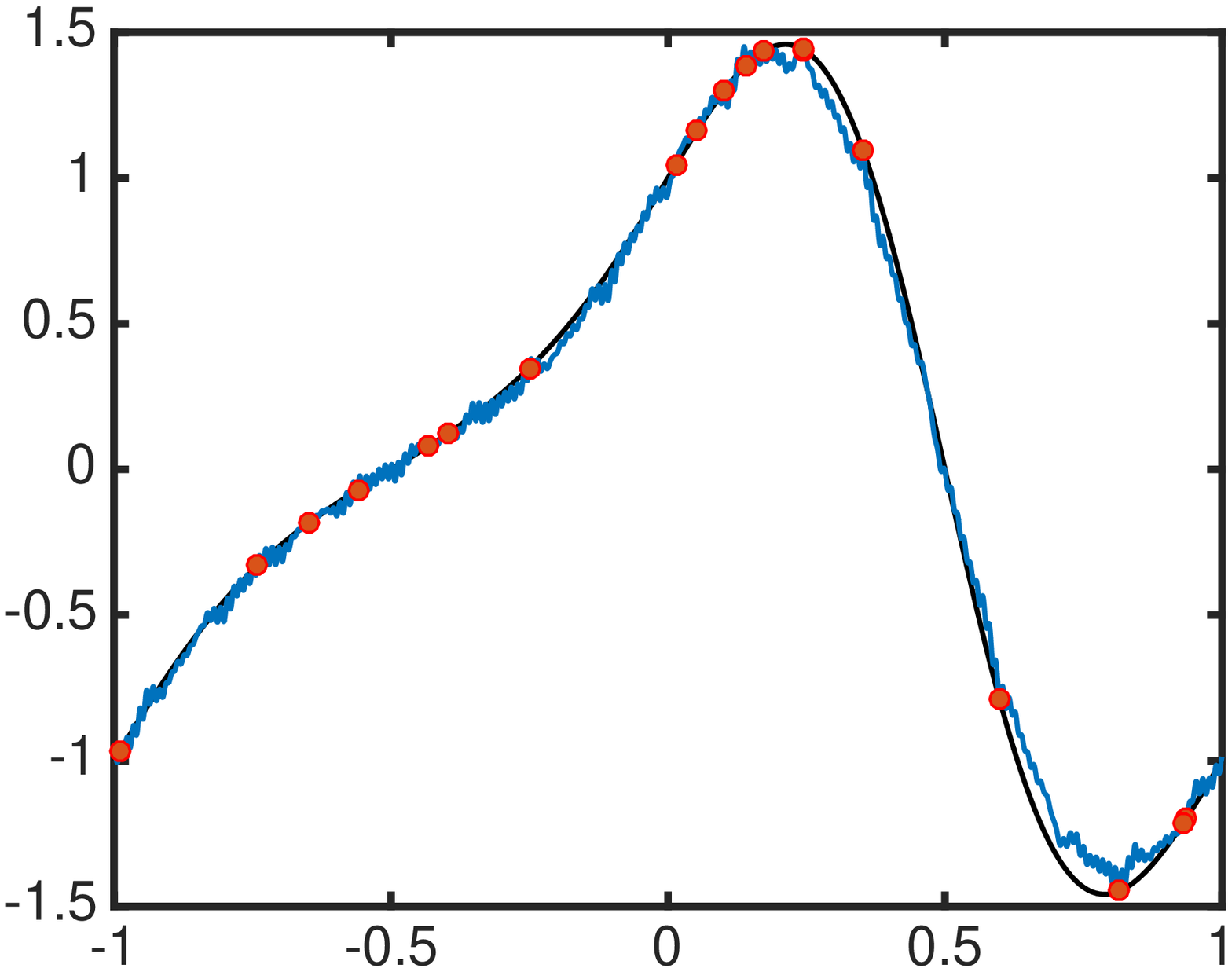}  &\includegraphics[width=4.5cm]{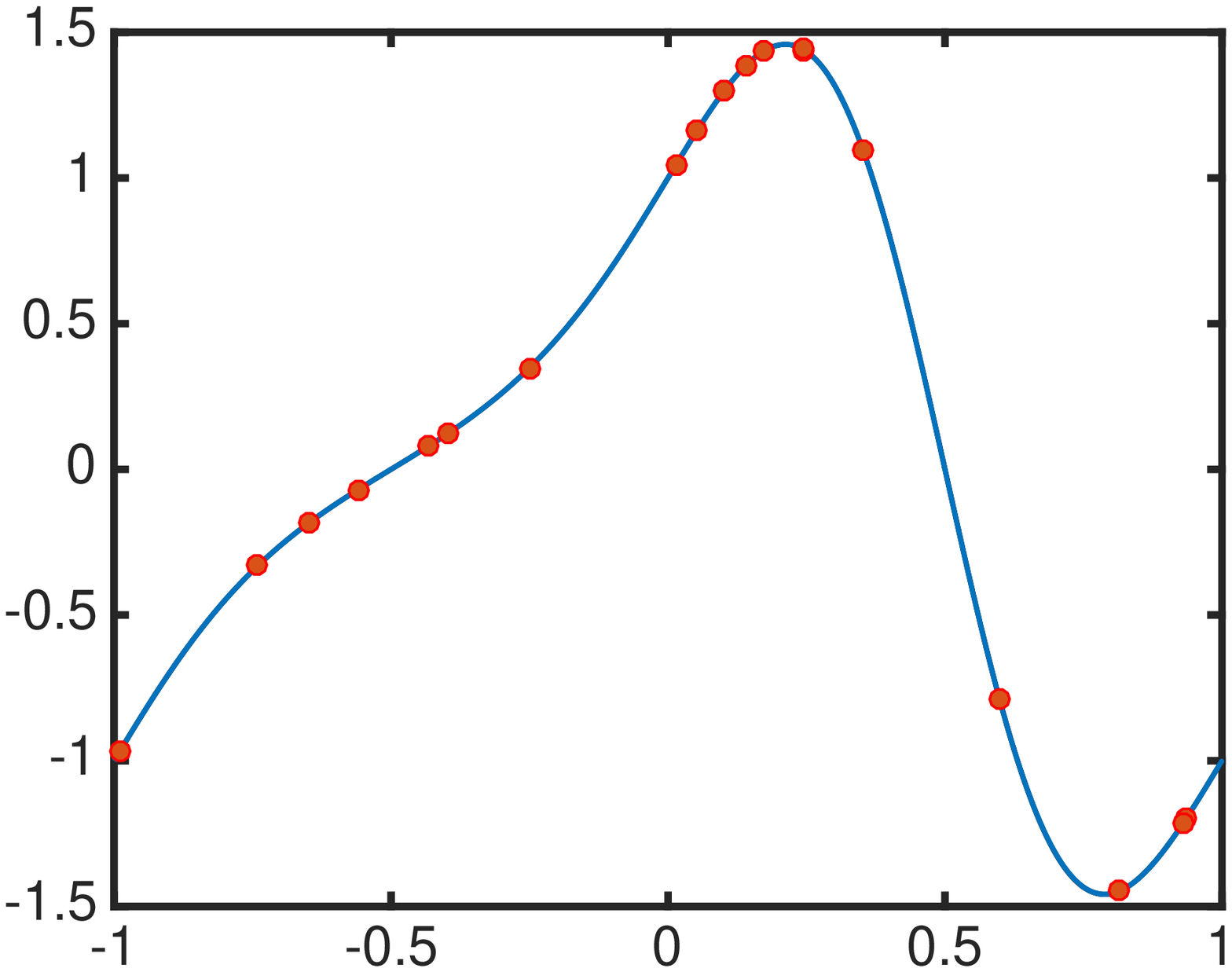} 
\\
\mbox{\small $\| f - \tilde{f} \|_{L^\infty} =\mbox{5.32e-1}$}  &\mbox{\small $ \| f - \tilde{f} \|_{L^\infty} =\mbox{1.63e-01}$} & \mbox{\small $\| f - \tilde{f} \|_{L^\infty} = \mbox{9.66e-04}$}
\\\\
\includegraphics[width=4.5cm]{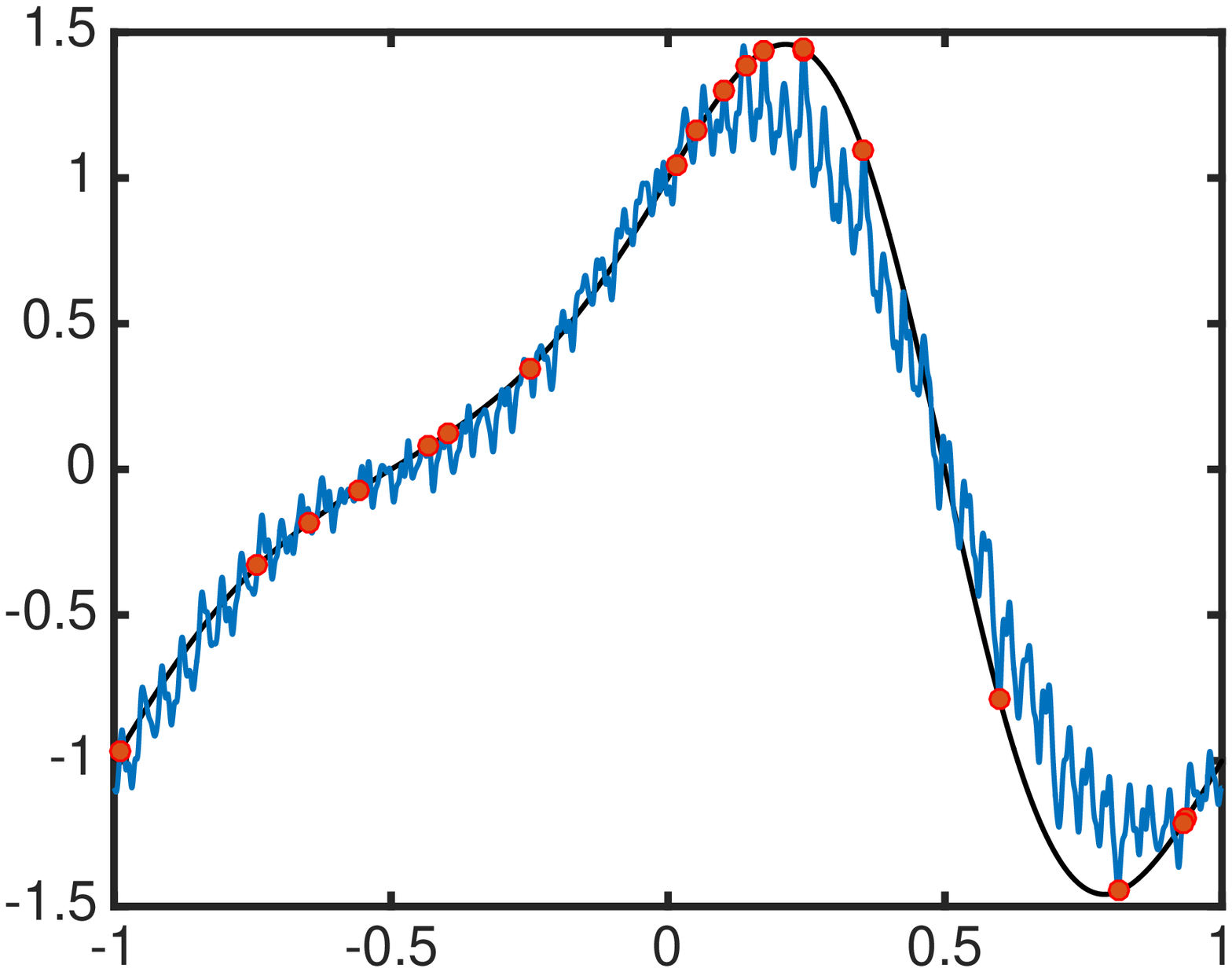}  &  \includegraphics[width=4.5cm]{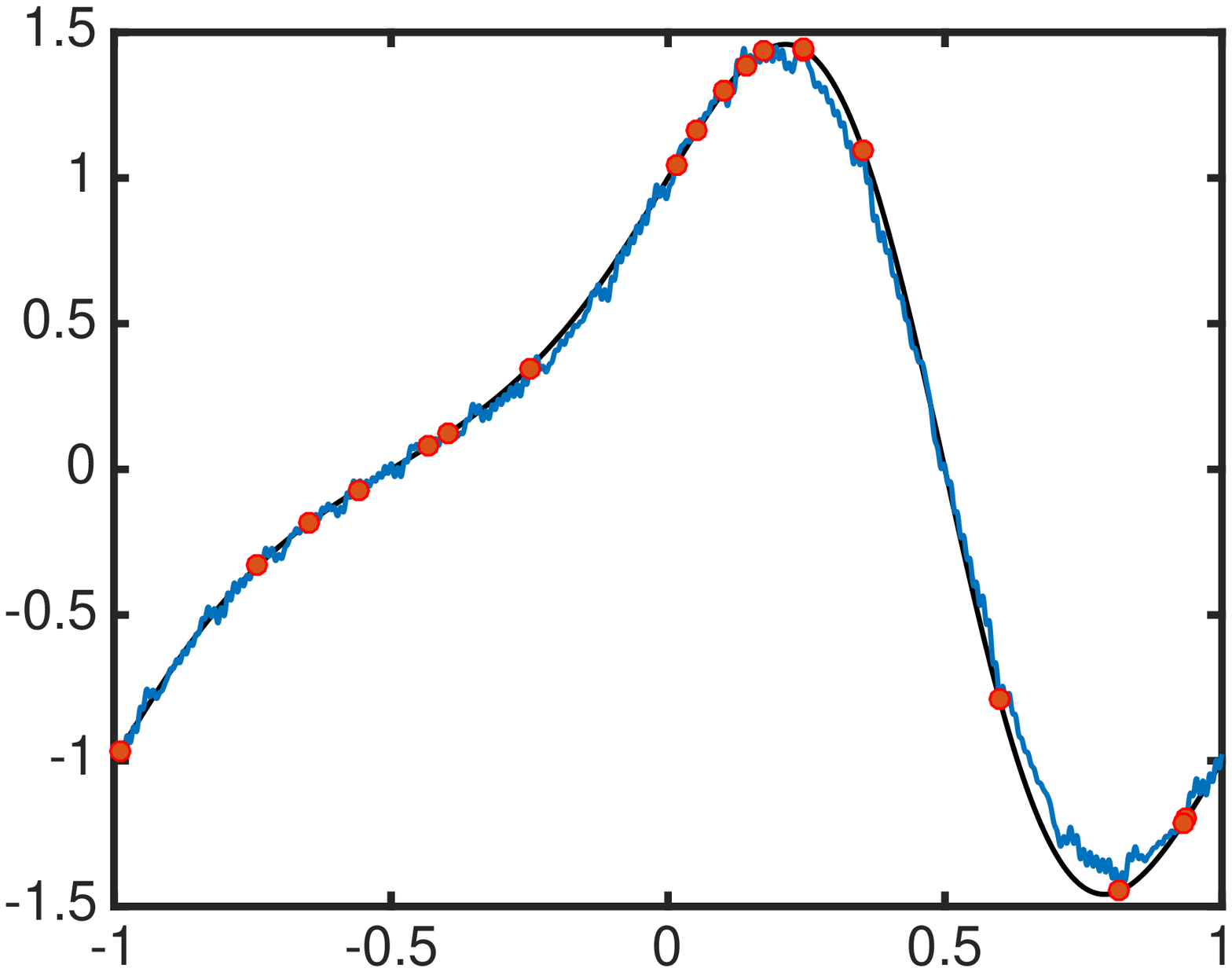}  &\includegraphics[width=4.5cm]{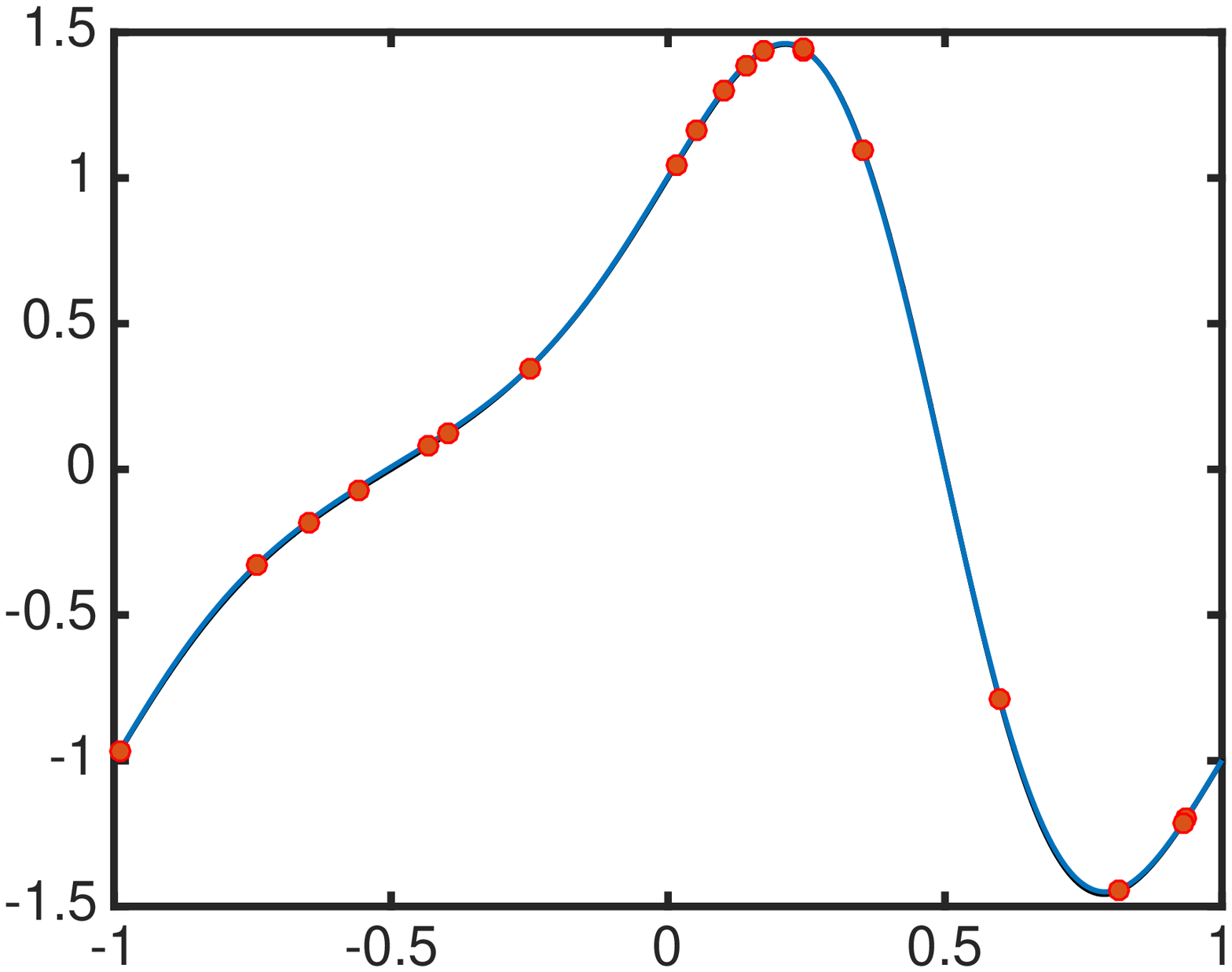} 
\\
\mbox{\small $\| f - \tilde{f} \|_{L^\infty} =\mbox{5.38e-1}$}  & \mbox{\small $\| f - \tilde{f} \|_{L^\infty} =\mbox{1.64e-01}$} & \mbox{\small $\| f - \tilde{f} \|_{L^\infty} = \mbox{1.02e-2}$}
\end{array}$
\caption{Recovery of the function $f(t) = \cos(\pi x) \exp(\sin(\pi x))$ (shown in black) from $N=20$ data points (shown in red).  The blue curve is the function $\tilde{f}$ obtained from (weighted) $\ell^1$ minimization using the Fourier basis with weights $w_{i} = 1$ (left), $w_i =1+|i|^{1/10}$ (middle) and $w_i = 1+|i|^{1/2}$ (right).  Top row: equality-constrained minimization \R{inf_min_noiseless}.  Bottom row: inequality-constrained minimization \R{inf_min} with $\eta = 10^{-2}$. } \label{f:Alias2}
\end{center}
\end{figure}

Now suppose that weights $w_i$ are added, and \R{bad_prob} is replaced by
\be{
\label{good_problem}
\inf_{z \in \ell^1_w(\bbZ)} \| z \|_{1,w}\quad \mbox{subject to $U z = y$}.
}
Assume the weights $w$ satisfy $w_{-i} = w_{i}$, $i \in \bbN$ and $1 \leq w_0 < w_1 < w_2< \ldots$, and consider the case of $f = \phi_0$ once more.  Then none of the aliased solutions of \R{bad_prob} are solutions of \R{good_problem}, since they all have larger weighted $\ell^1$-norm: $\| e_{2kP} \| = w_{2kP} > w_0 = \| e_0 \|$.  Hence, adding growing weights regularizes the problem \R{good_problem} and removes the bad, aliased solutions of \R{bad_prob}.  This improvement is illustrated numerically in Fig.\ \ref{f:Alias2}.  This figure also shows that this phenomenon is not limited to the equality-constrained minimization problem.

\begin{figure}
\begin{center}
$\begin{array}{ccc}
\includegraphics[width=4.5cm]{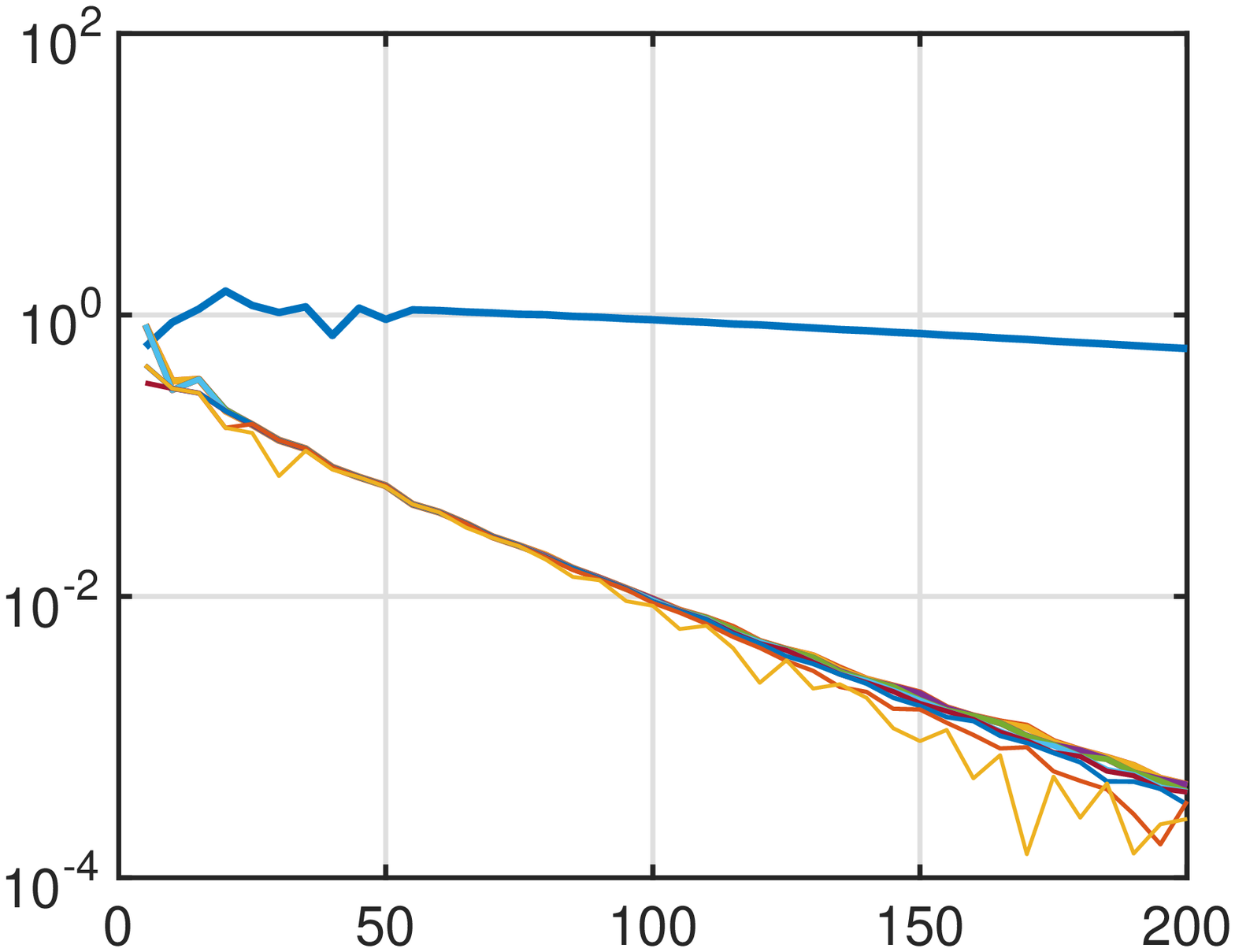}  & 
\includegraphics[width=4.5cm]{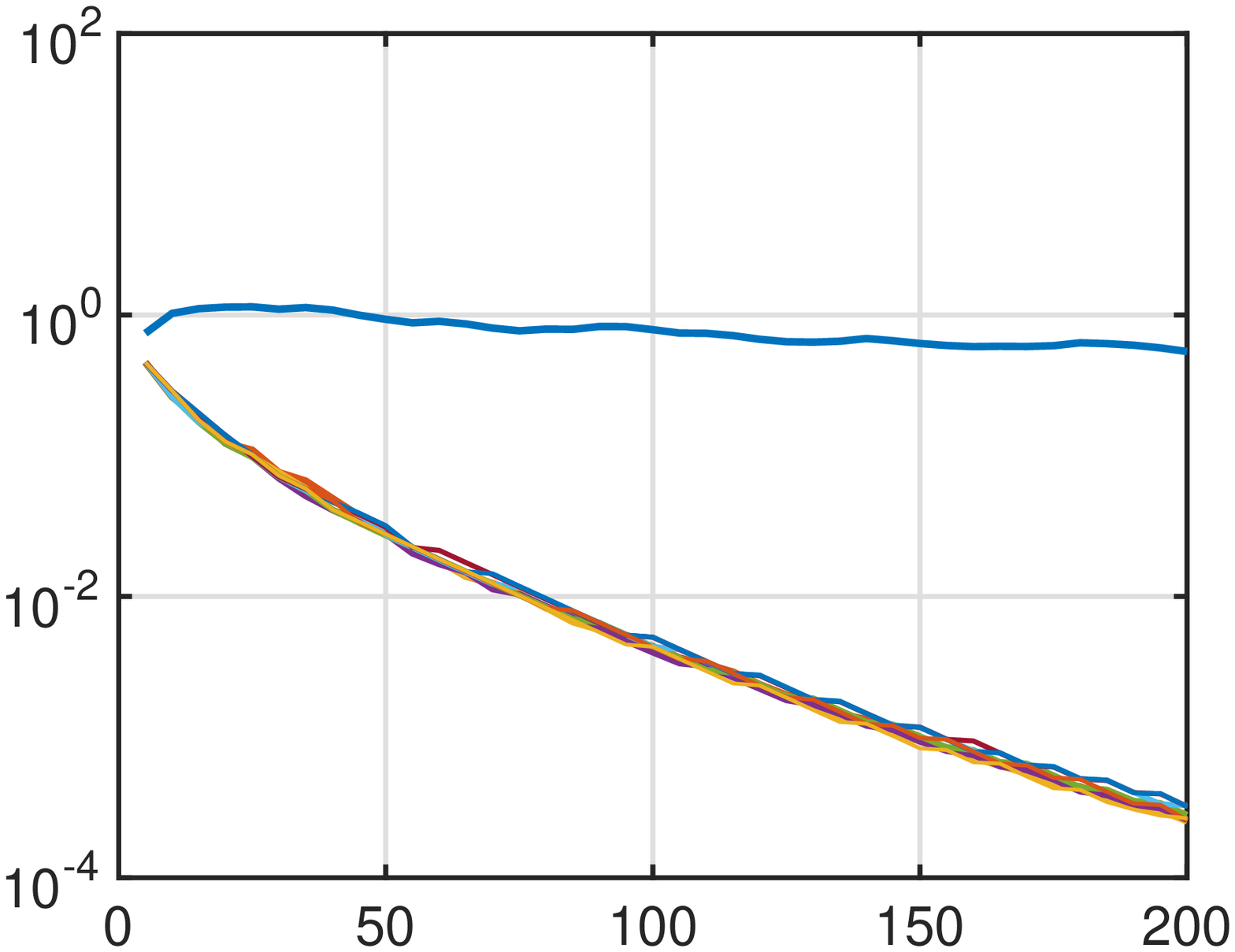}  &  
\includegraphics[width=4.5cm]{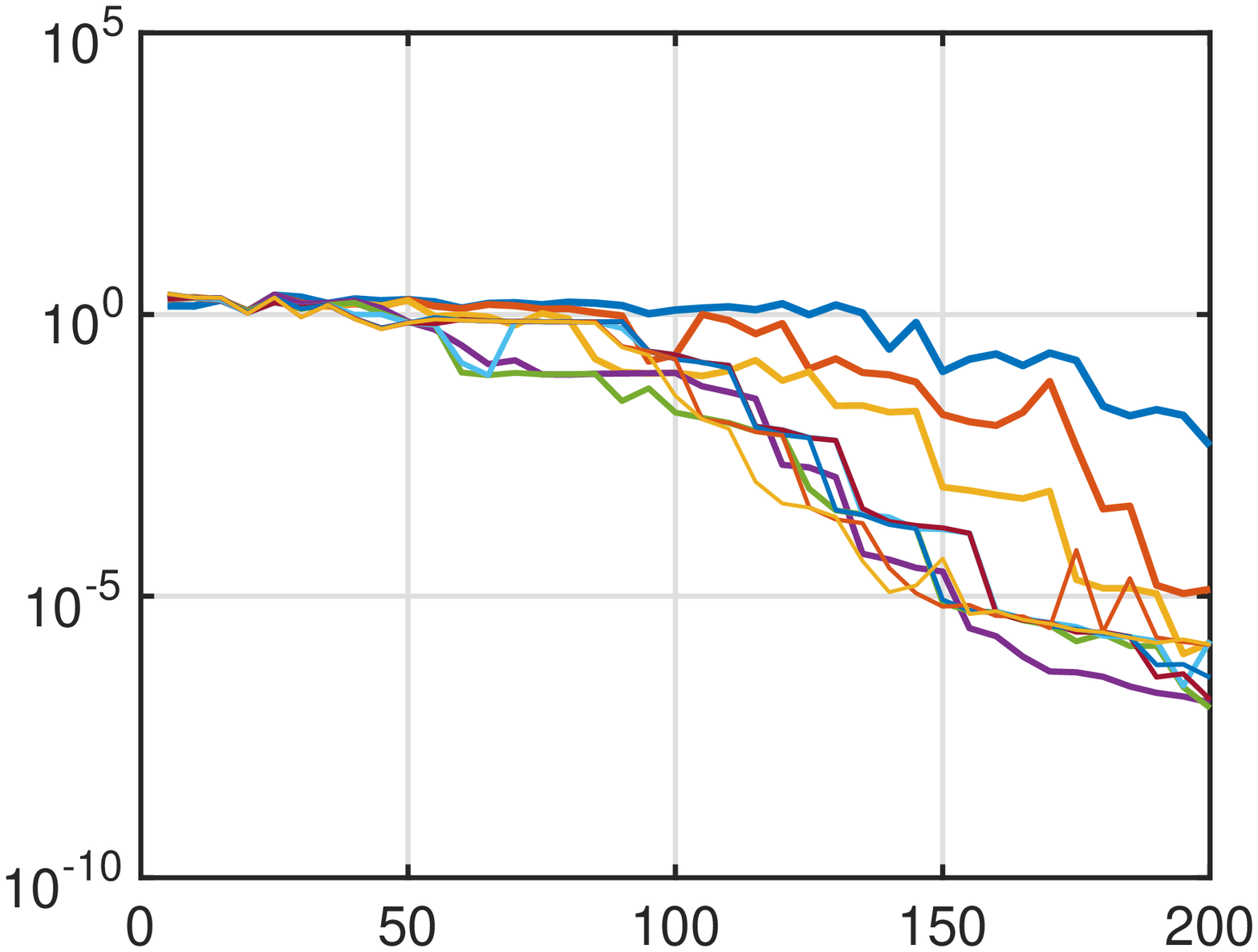}  
\\
\mbox{\small $f(t) = \frac{1}{1+25 t^2}$} & \mbox{\small $f(t) = \frac{1}{35-34t}$} & \mbox{\small $f(t) = \cos(30 t)$}
\\
\includegraphics[width=4.5cm]{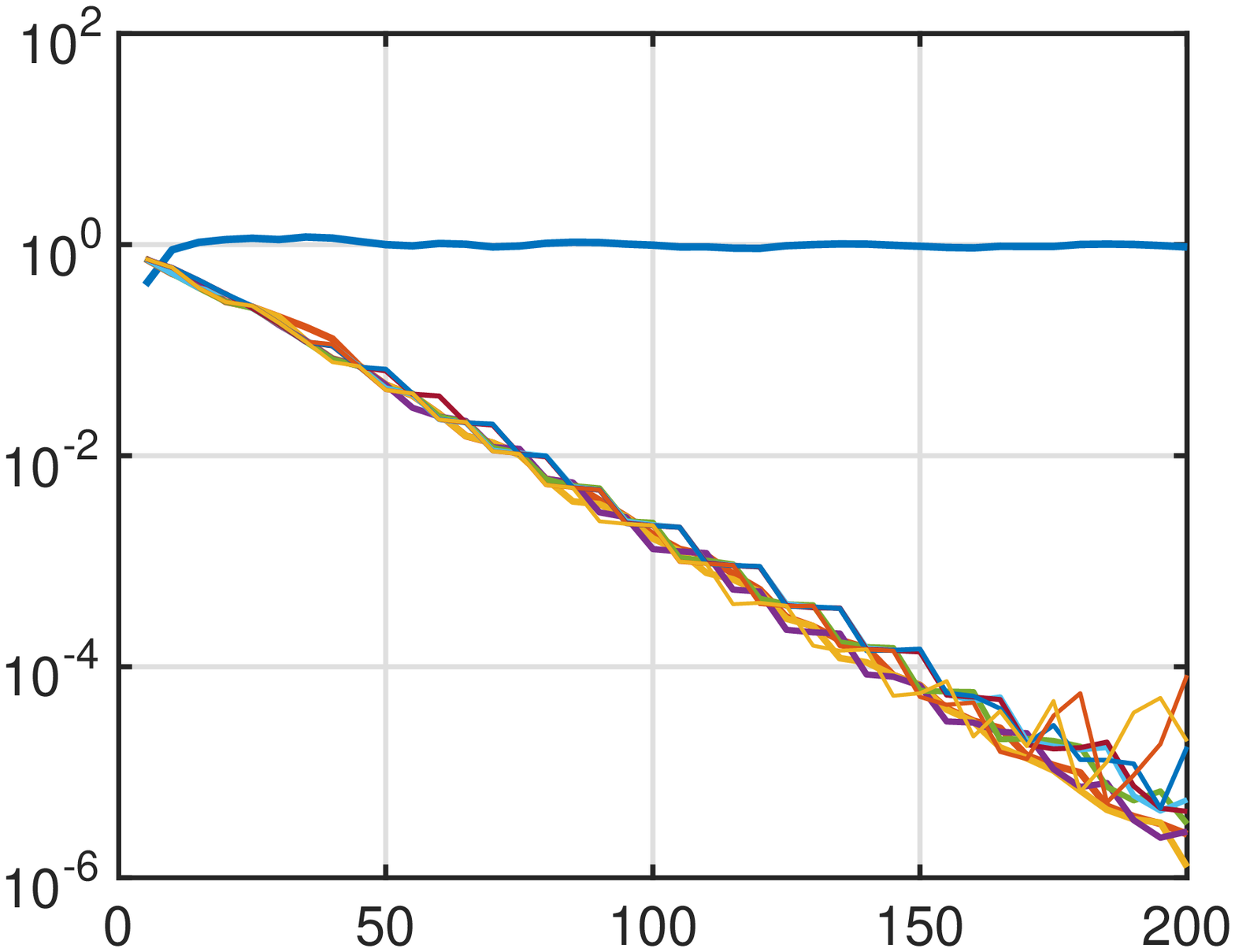} &
\includegraphics[width=4.5cm]{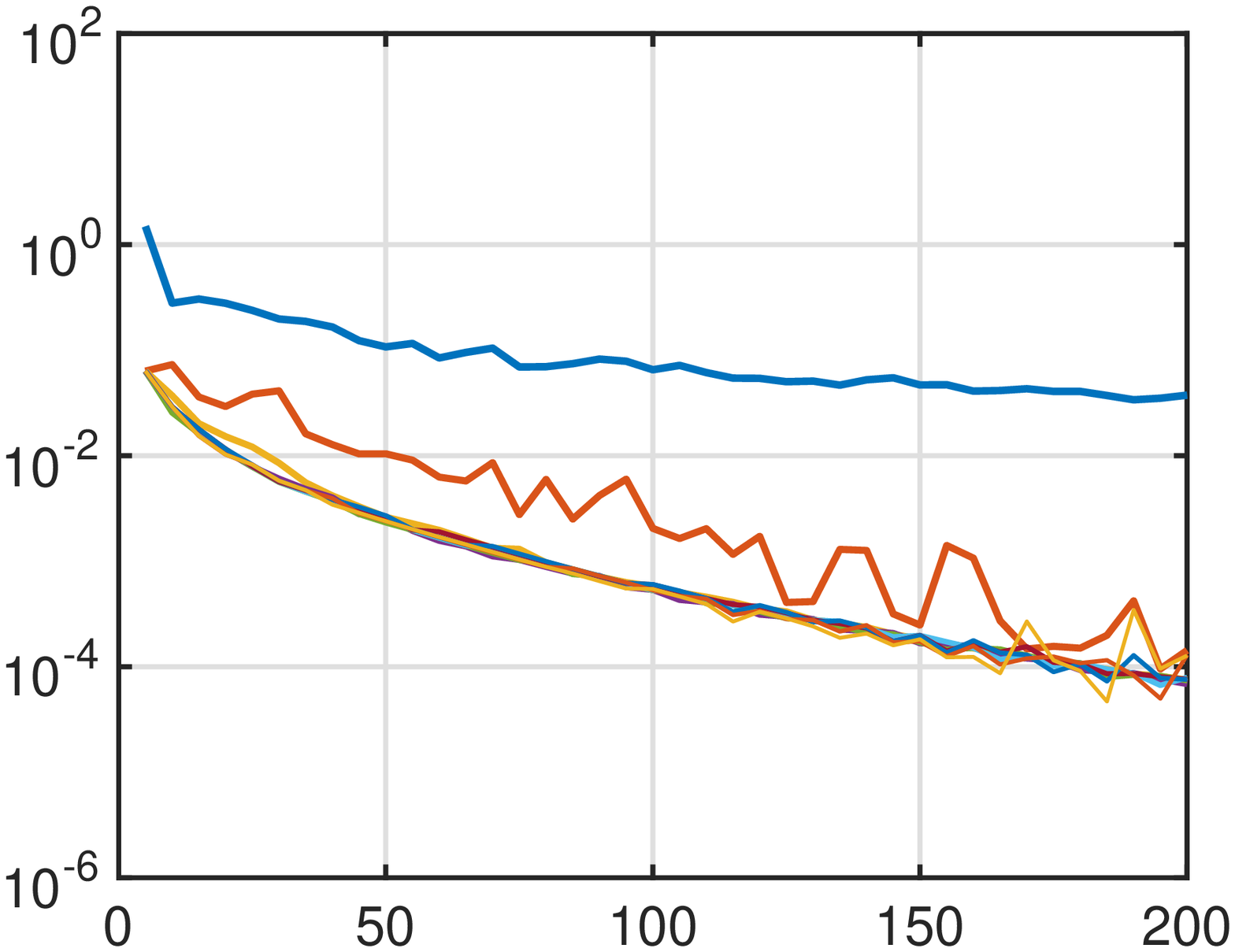} &
\includegraphics[width=4.5cm]{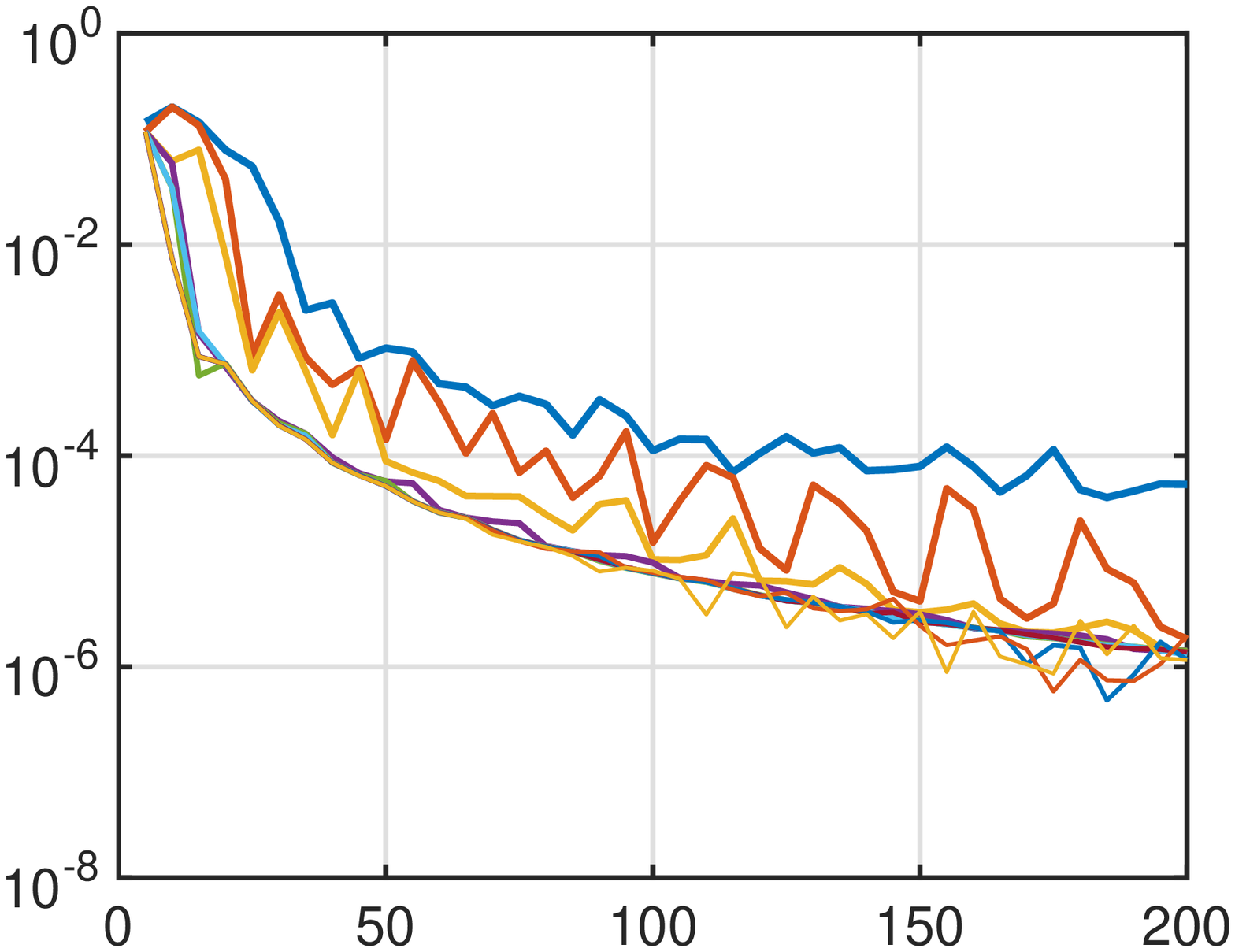}
\\
\mbox{\small $f(t) = \frac{\cosh(30 t^2)}{\cosh(30)}$} & \mbox{\small $f(t) = \sqrt{1.01+t}$} & \mbox{\small $f(t) = t^5 \log(t^2)$}
\end{array}$
\caption{Weighted $\ell^1$ minimization with equispaced data and Chebyshev polynomials.  The error $\| f - \tilde{f} \|_{L^\infty}$ against $N$ is shown for the choice $w_{i} = i^{\gamma}$, where $\gamma = 0.0,0.05,0.1,0.25,0.5,0.75,1.0,1.5,2.0,2.5$ (thickest to thinnest).  The truncation parameter $K=4N$ was used.  As with all numerical results in this paper, the minimization problem \R{fin_min_noiseless} was solved using the CVX optimization package.} \label{f:Alias3}
\end{center}
\end{figure}

\rem{
The use of weighted minimization strategies has been occasionally motivated by the desire to match the decay of the true coefficients $x$ of the unknown function and thereby obtain better approximations \cite{PengHamptonDoostantweighted,RauhutWardWeighted}.  However, this is not the primary role the weights in play in this setting.  To demonstrate this point, in Fig.\ \ref{f:Alias3} we plot the error for weighted $\ell^1$ minimization using Chebyshev polynomials for a number of different test functions and weighting strategies.  As can be seen, increasing the weights does not lead to a consistent improvement across all functions, even though all functions used (beside the final one) are infinitely smooth and thus have coefficients which decay superalgebraically fast.  While weights might help in some small way by promoting smoothness, these results suggest that the effect on the approximation error is much less than the role they play in regularizing the problem.    Furthermore, higher weights may well cause problems for numerical solvers, due to the increasing ill-conditioning of the $N \times K$ system matrix $U W^{-1} P_K$.

We note in passing that this situation is quite unlike the case of weighted $\ell^2$ minimization (see, for example, \cite{MSNTrigPoly,MSNJCP}), in which case the error for a smooth function decays only algebraically fast at a rate dependent on the algebraic growth rate of the weights.  Thus, for $\ell^2$ minimization, rapidly-growing weights promote smoothness.  Conversely, we will prove later that for weighted $\ell^1$ minimization the error decays superalgebraically fast for all smooth functions whenever the weights meet a minimum growth condition (Theorem \ref{t:Leg_poly_full_l1}).
}

\section{Approximation error of weighted $\ell^1$ minimization}\label{s:lin_approx_err}

The remainder of this paper is devoted to the analysis of the problems \R{fin_min_noiseless} and \R{fin_min}.  In this section, we present a linear approximation error analysis.  Truncation and the choice of the parameter $K$ is addressed in \S \ref{s:truncation}, and in \S \ref{s:AlgPoly_Examp} and \S\ref{s:TrigPoly_Examp} we apply these results to the examples of \S \ref{ss:examples}.

\subsection{A general recovery result}

We first require the following result, which bounds the error of weighted $\ell^1$ minimization subject to the existence of a particular dual vector $u$:

\lem{
\label{l:dual_certificate}
Let $\Delta \subseteq \{1,\ldots,K\}$.  Suppose that 
\bes{
(i): \| P_{\Delta} U^*  U P_{\Delta} - P_{\Delta} \| \leq \alpha,\qquad (ii): \max_{i \notin \Delta} \left \{ \|  U e_i \| / w_i \right \} \leq \beta,
}
and that there exists a vector $u = W^{-1} U^*  u'$ for some $u' \in \bbC^N$, such that
\bes{
(iii): \| W(P_{\Delta} u - \sgn(P_{\Delta} x)) \| \leq \gamma,\qquad (iv): \| P^{\perp}_{\Delta} u \|_{\infty} \leq \theta,\qquad (v): \| u' \| \leq L \sqrt{s},
}
where $s = \sum_{i \in \Delta} (w_i)^2$, for constants $0 \leq \alpha , \theta < 1$ and $\beta , \gamma, L \geq 0$ satisfying  
\bes{
\frac{\sqrt{1+\alpha} \beta \gamma}{(1-\alpha)(1-\theta) } < 1.
}
Let $\hat{x}$ be a minimizer of \R{fin_min}.  Then, if $\bar{x} \in \ell^1_w(\bbN)$ is feasible for \R{fin_min}, i.e. $\| U P_K \bar{x} - y \| \leq  \eta$, the error estimate
\be{
\label{l1_error_est}
\| \hat{x} - x \| \leq 2\left ( C_1 + C_2 L \sqrt{s} \right )  \eta + C_2 \left (2 \| P^{\perp}_{\Delta} x \|_{1,w} + \| x - \bar{x} \|_{1,w} \right ),
}
holds, 
where $C_1 = \left ( 1 + \frac{\gamma}{1-\theta} \right ) C_0$, $C_2 =  \frac{\beta}{1-\theta}\left ( 1 + \frac{\gamma}{1-\theta} \right ) C_0 + \frac{1}{1-\theta}$ and $C_0 = \left ( 1 - \frac{\sqrt{1+\alpha} \beta \gamma }{(1-\alpha)(1-\theta)} \right )^{-1} \frac{\sqrt{1+\alpha}}{1-\alpha}$.
}
Recall that the problem \R{fin_min_noiseless} can be viewed as a special case of \R{fin_min} corresponding to the case $\eta= 0$.  Hence this result considers only \R{fin_min}.  The proof of this lemma is given in \S \ref{ss:dual_certificate_proof}.

In practice, we shall use following result, which is a straightforward consequence of Lemma \ref{l:dual_certificate}:

\lem{
\label{l:Delta_recov}
Let $\Delta \subseteq \{1,\ldots,K \}$.  Suppose that there are constants $0 \leq \alpha,\theta < 1$ such that
\bes{
(a): \|  P_{\Delta} U^*U P_{\Delta} - P_{\Delta} \| \leq \alpha,\qquad (b): \| P^{\perp}_{\Delta} W^{-1} U^* U P_{\Delta} A^{-1} W P_{\Delta} \sgn(x) \|_{\infty} \leq \theta,
}
where $A = P_{\Delta} U^* U P_{\Delta}$.  If $\hat{x}$ is a minimizer of \R{fin_min} and $s = \sum_{i \in \Delta} (w_i)^2$, then
\be{
\label{err_bound}
\| \hat{x} - x \| \leq 2 C_3 \left ( 1 + \frac{C_3+1}{1-\theta} \sqrt{s} \right )  \eta  + \frac{C_3+1}{1-\theta}  \left ( 2 \| P^{\perp}_{\Delta} x \|_{1,w} + \| x - \bar{x} \|_{1,w} \right ) ,
}
where $\bar{x} \in \ell^1_{w}(\bbN)$ is any feasible solution of \R{fin_min} and $C_3 = \frac{\sqrt{1+\alpha}}{1-\alpha} $.
}
\prf{
We apply Lemma \ref{l:dual_certificate} with $u = W^{-1} U^* U P_{\Delta} A^{-1} W P_{\Delta} \sgn(x)$.  Note that $(a)$ and $(b)$ imply $(i)$ and $(iv)$ respectively.  Also, by construction, $P_{\Delta} u = P_{\Delta} \sgn(x)$ and therefore $(iii)$ holds with $\gamma = 0$.  Now consider $(ii)$.  By definition
\bes{
\| U e_i \| = \| \phi_i \|_h \leq \| \phi_i \|_{\infty} \sqrt{\sum^{N}_{n=1} \tau_n } \leq w_i,
}
where the last inequality is due to \R{weights_growth} and the fact that $\sum^{N}_{n=1} \tau_n = 1$ when the weights $\tau_n$ are given by \R{quad_weights_1}.  Hence $(ii)$ holds with $\beta = 1$.  Finally, observe that
\eas{
\| u' \| &= \| U P_{\Delta} A^{-1} W P_{\Delta} \sgn(x) \|
\leq \| U P_{\Delta} \| \| A^{-1} \| \| W P_{\Delta} \sgn(x) \|
 \leq \frac{\sqrt{1+\alpha}}{1-\alpha} \sqrt{s},
}
where the final inequality follows from $(a)$.  Hence $(v)$ holds with $L = \frac{\sqrt{1+\alpha}}{1-\alpha}$.
}

\subsection{Linear approximation error preliminaries}

Lemma \ref{l:Delta_recov} gives conditions under which $x$ is approximated with error depending on the magnitude of its coefficients $x_j$ outside some set $\Delta$.  This depends on the conditioning of the corresponding submatrix (condition (a)) and the off-support magnitude of the vector $u$ (condition (b)).  Under a sparsity condition on the coefficients $x$, and appropriate random choices of the points $T$, one may use this result to prove estimates relating the number of measurements to the sparsity \cite{AdcockCSFunInterp}.  However, as discussed in \S \ref{s:introduction}, in practice the data points may not arise from such distributions.  In this section, we consider arbitrary deterministic scattered data points and present a linear error analysis.  This follows by setting $\Delta = \{1,\ldots,M \}$ for some $M \leq K$ and using Lemma \ref{l:Delta_recov} to determine how large $M$ can be chosen in relation to the density $h$ of the points.  Doing this will allow us to make a direct comparison with other techniques (see Remark \ref{r:cooking_time}).

Since such statements will be asymptotic in $h \rightarrow 0$, we first require an additional assumption.  Let $H$ be a subspace of $L^2_{\nu}(D) \cap L^{\infty}(D)$ which is closed under multiplication and complex conjugation and such that $f \in H$ and $\{ \phi_i \}_{i \in \bbN} \subseteq H$.  We now assume that the points $T$ satisfy
\be{
\label{quad_conv}
\sum^{N}_{n=1} \tau_n g(t_n) \rightarrow \int_{D} g(t) \nu(t) \D t,\quad h \rightarrow 0,\qquad \forall g \in H.
}
In particular, since $H$ is closed under multiplication and complex conjugation, one has that
\bes{
\ip{f}{g}_{h} \rightarrow \ip{f}{g}_{L^2_{\nu}},\quad h \rightarrow 0,\qquad \forall f,g \in H.
}
Hence the discrete inner product is equivalent to $\ip{\cdot}{\cdot}_{L^2_{\nu}}$ on finite-dimensional subspaces of $H$ for sufficiently small $h$.  Note that this assumption is by no means stringent.  For example, if $D = (-1,1)$ we may take $H$ to be the space of functions for which $|f(t)|^2 \nu(t)$ is Riemann integrable.

Before stating our main result (Theorem \ref{t:full_samp_recov}), we first require some additional notation and technical lemmas.  For $h >0$ and $M,R \in \bbN$ we define the quantities
\be{
\label{E_2NM}
E_2(h,M) = \| P_M - P_M U^* U P_M \|,\quad E_{\infty}(h,M) = \| P_M - P_M U^* U P_M \|_{\infty} ,
}
and
\be{
\label{F_N_M_R}
F(h,M,R) = \| P^{\perp}_R W^{-1} U^* U P_M \|_{\infty}.
}
We also set $E(h,M) = \max \{ E_2(h,M),E_{\infty}(h,M) \}$.

\lem{
\label{E_decay}
For fixed $M$, we have $E(h,M) \rightarrow 0$ as $h \rightarrow 0$.
}
\prf{
Since all norms on $\bbC^M$ are equivalent, it suffices to show that $(U^* U)_{i,j}\rightarrow \delta_{ij}$ for each $i,j=1,\ldots,M$ as $h \rightarrow 0$.  However, by \R{quad_conv} and orthogonality of the $\phi_j$, we have
$
(U^* U)_{i,j} = \ip{\phi_i}{\phi_j}_{h} \rightarrow \ip{\phi_i}{\phi_j}_{L^2_{\nu}} = \delta_{ij},
$
as required.
}

\lem{
\label{F_decay}
Suppose that the weights $w_i = z_i \| \phi_i \|_{L^\infty}$ for some $z_i \geq 1$.  Then
\bes{
F(h,M,R) \leq \frac{\sqrt{M} \sqrt{1+E(h,M)}}{\inf_{i > R} \{ z_i \} }.
}
}

\prf{
Let $x \in P_M(\ell^2(\bbN))$, $\| x \|_{\infty} = 1$ be arbitrary.  Then $(W^{-1} U^* U P_M x)_{i} = \ip{g}{\phi_i}_{h}/w_i$, where $g = \sum^{M}_{j=1} x_j \phi_j$.  By the Cauchy--Schwarz inequality,
\bes{
\| P^{\perp}_R W^{-1} U^* U P_M x \|_{\infty} \leq \sup_{i > R} \left \{ \frac{1}{w_i} \| g \|_h \| \phi_i \|_h \right \}.
}
Now $\| g \|^2_h = \ip{x}{P_M U^* U P_M x} \leq \left ( 1 + E(h,M) \right ) \| x \|^2 \leq M  \left ( 1 + E(h,M) \right )$.  Also $\| \phi_i \|_h \leq \| \phi_i \|_{L^\infty} = w_i / z_i$ and therefore
\bes{
\| P^{\perp}_R W^{-1} U^* U P_M \|_{\infty} \leq \frac{\sqrt{M} \sqrt{1+E(h,M)}}{\inf_{i > R} \{ z_i \} },
}
as required.
}

\lem{
Suppose that the weights $w_i = z_i \| \phi_i \|_{L^\infty}$, where $z_i \geq 1$ and $z_i \rightarrow \infty$ as $i \rightarrow \infty$.  Then for any $0 < \epsilon < 1/2$ and any $M \in \bbN$ there exists a $R \in \bbN$ and $h > 0$ depending only on $M$ and $\epsilon$ such that
\be{
\label{NMR_condition}
E(h,M) < \epsilon,\quad 
E(h,R) < \epsilon \frac{\min_{M < i \leq R} \{ w_i \} }{\max_{i=1,\ldots,M} \{ w_i \} },
\quad
F(h,M,R) \leq \frac{\epsilon}{\max_{i=1,\ldots,M} \{ w_i \} }.
}
}
\prf{
By Lemma \ref{E_decay} we can find an $h_1$ such that $E(h,M) < \epsilon$ for all $h \leq h_1$, thus satisfying the first condition in \R{NMR_condition}.  Using Lemma \ref{F_decay}, we note that
\bes{
F(h,M,R) \leq \frac{\sqrt{2M}}{\min_{i > R} \{ z_i \}},\quad \forall h \leq h_1.
}
To satisfy the third condition in \R{NMR_condition}, we pick $R$ sufficiently large so that
\bes{
\inf_{i > R} \{ z_i \} > \sqrt{2M} \max_{i=1,\ldots,M} \{ w_i \} / \epsilon.
}
We now pick $h_2$ sufficiently small so that
\bes{
E(h,R) < \epsilon \frac{\min_{M < i \leq R} \{ w_i \} }{\max_{i=1,\ldots,M} \{ w_i \} },\quad \forall h \leq h_2,
}
and then set $h = \min \{ h_1,h_2 \}$.
}

\subsection{Main result}

In order to present our main result, we first introduce the following notation:
\be{
\label{trunc}
T_{h,K,\eta}(x) = \inf \left \{ \| x - \bar{x} \|_{1,w} : \bar{x} \in \bbC^K, \| U P_K \bar{x} - y \| \leq \eta \right \}.
}

\thm{
\label{t:full_samp_recov}
Suppose that the weights $w_i = z_i \| \phi_i \|_{L^\infty}$, where $z_i \geq 1$ and $z_i \rightarrow \infty$ as $i \rightarrow \infty$.  For $0 < \epsilon < 1/2$, let $h>0$ and $M,R \in \bbN$ be such that \R{NMR_condition} holds.  Then there exists a constant $C(\epsilon)$ such that, whenever $\hat{x}$ is a minimizer of \R{fin_min},
\be{
\label{l1w_err_bd}
\| x - \hat{x} \| \leq C(\epsilon) \left [ \left ( 1 + \| P_M w \| \right )  \eta +  \| P^{\perp}_M x \|_{1,w} + T_{h,K, \eta }(x) \right ],
}
where $T_{h,K, \eta }(x)$ is as in \R{trunc}.  Moreover, $\lim_{\epsilon \rightarrow 0^+} C(\epsilon) = 4$.
}
Note that the weights condition is equivalent to \R{weights_growth2} which was shown to be necessary in \S \ref{s:need}.  Theorem \ref{t:full_samp_recov} shows that the same condition is also sufficient.  We note also that the truncation parameter $K$ only influences the term $T_{h,K,\eta}(x)$.  We defer the detailed analysis of this term to \S \ref{s:truncation}.

\prf{
We use Lemma \ref{l:Delta_recov} with $\Delta = \{1,\ldots,M\}$.  Note first that
\be{
\label{alpha_bound}
 E(h,M) < \epsilon <1,
}
and therefore $(a)$ holds with $\alpha \leq \epsilon$.  Now consider $(b)$.  Write 
\be{
\label{u_full_def}
u = W^{-1} U^* U P_M A^{-1} W P_M \sgn{(x)},
}
so that $(b)$ is equivalent to $\| P^\perp_M u \|_{\infty} \leq \theta$.  We have 
\be{
\label{R_splitting}
\| P^{\perp}_M u \|_{\infty} = \max \left \{ \| P_R P^{\perp}_M u \|_{\infty} , \| P^{\perp}_R u \|_{\infty} \right \}.
}
We consider both terms separately.  For $\| P_R P^{\perp}_M u \|_{\infty}$, \R{u_full_def} gives
\eas{
\| P_R P^{\perp}_M u \|_{\infty} &\leq \| P_R P^{\perp}_M W^{-1} \|_{\infty} \| P_R P^{\perp}_M U^* U P_M \|_{\infty} \| A^{-1} \|_{\infty} \| P_M W \|_{\infty} \| \sgn{(x)} \|_{\infty}
\\
& \leq \frac{\max_{1 \leq i \leq M} \{ w_i \} }{\min_{M<i\leq R} \{ w_i \}} \| P_R P^{\perp}_M U^* U P_M \|_{\infty} \| A^{-1} \|_{\infty}.
}
Note that $\| I - A \|_{\infty}  \leq E(h,M)$ and therefore $\| A^{-1} \|_{\infty} \leq \frac{1}{1-E(h,M)}$.  Moreover
\eas{
\| P_R P^{\perp}_M U^* U P_M \|_{\infty} &= \| P_R P^{\perp}_M (I - U^* U ) P_M \|_{\infty}
 \leq \| P_R (I - U^* U ) P_R \|_{\infty}
 \leq E(h,R),
}
where the first equality is due to the fact that $P^{\perp}_M P_M = 0$.  Therefore, we obtain
\be{
\label{M_R_split_term}
\| P_R P^{\perp}_M u \|_{\infty} \leq \frac{\max_{1 \leq i \leq M} \{ w_i \} }{\min_{M<i\leq R} \{ w_i \}} \left ( \frac{E(h,R)}{1-E(h,M)} \right ) \leq \frac{\epsilon}{1-\epsilon}.
}
Now consider the other term in \R{R_splitting}.  By \R{u_full_def} and the definition of $F(h,M,R)$,
\eas{
\| P^{\perp}_R u \|_{\infty} &\leq \| P^{\perp}_R W^{-1} U^* U P_M \|_{\infty} \| A^{-1} \|_{\infty} \| P_M W \|_{\infty} 
 \leq  \frac{\max_{i=1,\ldots,M} \{ w_i \} }{1-E(h,M)} F(h,M,R)
 \leq \frac{\epsilon}{1-\epsilon}.
}
Combining this with \R{M_R_split_term} and substituting into \R{R_splitting} yields $\theta \leq \frac{\epsilon}{1-\epsilon} < 1$.
The result now follows immediately from Lemma \ref{l:Delta_recov}.
}

\subsection{Comparison with least-squares fitting}\label{ss:LS_fit_comp}
The following result is standard (a proof is given for completeness):

\thm{
\label{t:LS_err}
For $0 < \epsilon < 1$ let $M \in \bbN$ and $h>0$ be such that $E_{2}(h,M) \leq \epsilon < 1$, where $E_2(h,M)$ is as in \R{E_2NM}.  Then \R{LS_fit} has a unique solution $\check{x}$, and this satisfies
\be{
\label{LS_err_bd}
\| x - \check{x} \| \leq \left ( 1 + \frac{1}{\sqrt{1-\epsilon}} \right ) \| x - P_M x \|_{1,w} + \frac{1}{\sqrt{1-\epsilon}}  \eta,
}
for any $w = \{ w_i \}_{i \in \bbN}$ with $w_i \geq \| \phi_i \|_{L^\infty}$.
}
\prf{
Observe that
\eas{
\| U P_M z \|^2 &= z^* P_M U^* U P_M z 
= \| z \|^2 - z^* \left ( P_M -  P_M U^* U P_M \right ) z 
\geq \left ( 1 - E_2(h,M) \right ) \| z \|^2.
}
Hence $U P_M$ has full rank since $E_{2}(h,M) \leq \epsilon < 1$, and its minimum singular value satisfies $\sigma_{\min} \geq \sqrt{1-\epsilon}$.  This implies that $\check{x}$ is unique and is given by $\check{x} = (U P_M)^{\dag} y$.  Hence
\eas{
\check{x} - P_M x &= (U P_M)^{\dag} \left ( U x + \{ \sqrt{\tau_n} e_n \}^{N}_{n=1} \right )  - P_M x
 = (U P_M)^{\dag} \left ( U (x - P_M x) + \{ \sqrt{\tau_n} e_n \}^{N}_{n=1} \right ).
}
Therefore
\bes{
\| \check{x} - P_M x \| \leq \frac{1}{\sqrt{1-\epsilon}} \left ( \| U(x-P_M x) \| +  \eta \right ) \leq \frac{1}{\sqrt{1-\epsilon}} \left ( \| x - P_M x \|_{1,w} +  \eta \right ) .
}
Since $\| x - \check{x} \| \leq \| x - P_M x \| + \| \check{x} - P_M x \| \leq \| x - P_M x \|_{1,w} + \| \check{x} - P_M x \|$, the result follows.
}

The error bound \R{LS_err_bd} is similar to the bound \R{l1w_err_bd} for weighted $\ell^1$ minimization.  In the absence of noise and truncation error, both depend on the term $x-P_M x$, i.e.\ the tail of $x$ beyond its first $M$ coefficients.  The primary difference is in the size of $M$, which is determined through the conditions of Theorems \ref{t:full_samp_recov} and \ref{t:LS_err} respectively.  We shall discuss this point further in \S \ref{s:AlgPoly_Examp} and \S \ref{s:TrigPoly_Examp}.  But we first reiterate (see also \S \ref{ss:LS_comp}) that $M$ is a fixed parameter for least squares, required for implementation, whereas for weighted $\ell^1$ it is introduced solely to provide an estimate for approximation error.  Lemma \ref{l:dual_certificate} asserts that weighted $\ell^1$ minimization can recover coefficients of $x$ corresponding to other subsets $\Delta$, provided the various conditions hold.  Note also that for weighted $\ell^1$ minimization the parameter $M$ appearing in Theorem \ref{t:full_samp_recov} is in no way related to the truncation parameter $K$, besides the relation $M \leq K$.

\rem{
\label{r:LSl1weightsdiff}
The error bound \R{LS_err_bd} also differs from \R{l1w_err_bd} in that it holds for the weights $w_i = \| \phi_i \|_{L^{\infty}}$ and its noise term does not involve the factor $\| P_M w \|$.  See \S \ref{s:AlgPoly_Examp} for further discussion.  We also remark that the term $\| x - P_M x \|_{1,w}$ in \R{l1w_err_bd} can be improved slightly to $\| f - \tilde{f} \|_{L^\infty}$, where $\tilde{f} = \sum^{M}_{i=1} x_i \phi_i$ is the projection of $f$ onto $\spn \{ \phi_1,\ldots,\phi_M \}$.  We opt for \R{LS_err_bd} so that a direct comparison can be made with \R{l1w_err_bd}.  Note that one also has $\| x - P_M x \| \leq \| f - \tilde{f} \|_{L^\infty} \leq \| x - P_M x \|_{1,w}$.
}

\rem{
\label{r:cooking_time}
Theorems \ref{t:full_samp_recov} and \ref{t:LS_err} provide a recipe for determining the worst-case behaviour of weighted $\ell^1$ minimization for scattered data.  This is as follows.  Given a orthonormal system $\{ \phi_i \}_{i \in \bbN}$, an $h > 0$ and an $0 < \epsilon < 1/2$, determine:
\bull{
\item[1.] the largest $M = M_1(h)$ such that $E_{2}(h,M) < \epsilon$,
\item[2.] the largest $M = M_2(h)$ such that \R{NMR_condition} holds for some appropriate $R$ and $\{ w_i \}_{i \in \bbN}$.
}
In this case, the errors for both weighted $\ell^1$ minimization and least-squares fitting are determined by $\| x - P_{M} x \|_{1,w}$, where $M = M_1(h)$ for the former and $M = M_2(h)$ for the latter.  Hence, if $M_1(h) \asymp M_2(h)$ as $h \rightarrow 0$ it follows that both least-squares fitting and weighted $\ell^1$ minimization are guaranteed to converge at roughly the same asymptotic rate as $h \rightarrow 0$.  Since $M_1(h)$ and $M_2(h)$ are dependent on the system $\{ \phi_i \}_{i \in \bbN}$ a separate analysis must be carried out in each case.  The last two sections of this paper will be devoted to doing this for the examples of \S \ref{ss:examples}.  
}

\subsection{Proof of Lemma \ref{l:dual_certificate}}
\label{ss:dual_certificate_proof}

To complete this section we give the proof of Lemma \ref{l:dual_certificate}.

\prf{[Proof of Lemma \ref{l:dual_certificate}]
Let $v = \hat{x} - x$.  Then $A v =  P_{\Delta} U^*  U v -  P_{\Delta} U^*  U P^{\perp}_{\Delta} v$, where $A$ is the restriction of $P_{\Delta} U^*  U P_{\Delta} $ to $P_{\Delta}(\ell^2(\bbN))$.  By $(i)$, we have $\| A^{-1} \| \leq \frac{1}{1-\alpha}$ and
\bes{
\|  P_{\Delta} U^*   \|^2 = \|   U P_{\Delta} \|^2 = \|  P_{\Delta} U^*  U P_{\Delta} \| \leq 1 + \alpha.
}
Thus
\eas{
\| P_{\Delta} v \| &\leq \frac{1}{1-\alpha} \|  P_{\Delta} U^*  \| \|  U v \| + \frac{1}{1-\alpha} \|  P_{\Delta} U^*  U  P^{\perp}_{\Delta} v \|
 \leq \frac{\sqrt{1+\alpha}}{1-\alpha} \left ( \|  U v \| + \|   U P^{\perp}_{\Delta} v \| \right ).
}
Observe that
\be{
\label{UN_v}
\|  U v \| = \| U \hat{x} - U x \| \leq 2  \eta .
}
Hence
\bes{
\| P_{\Delta} v \|  \leq \frac{\sqrt{1+\alpha}}{1-\alpha}\left ( 2 \eta + \|   U P^{\perp}_{\Delta} v \| \right ). 
}
The second term can be estimated as follows:
\eas{
\|   U P^{\perp}_{\Delta} v \| \leq \sum_{i \notin \Delta} | v_i | \|   U e_i \| \leq \beta \| P^{\perp}_{\Delta} v \|_{1,w},
}
where the latter inequality is due to $(ii)$.  Hence we get
\be{
\label{v_Delta_Delta_perp}
\| P_{\Delta} v \|  \leq \frac{\sqrt{1+\alpha}}{1-\alpha} \left ( 2 \eta +\beta  \|  P^{\perp}_{\Delta} v \|_{1,w} \right ).
}
We shall return to this inequality later, but let us now consider $\hat{x}$.
\ea{
\nm{\hat{x}}_{1,w} &= \nm{P_\Delta\hat{x}}_{1,w} + \nmu{P^{\perp}_\Delta \hat{x} }_{1,w} \nn
\\
& \geq \Re \ip{P_\Delta W \hat{x}}{\mathrm{sign}(P_\Delta x)} + \nmu{P^{\perp}_\Delta  v }_{1,w} - \nmu{P^{\perp}_\Delta x}_{1,w} \nn
\\
&= \Re \ip{P_\Delta W v}{\mathrm{sign}(P_\Delta x)} + \nm{P_\Delta x}_{1,w} + \nmu{P^{\perp}_\Delta v }_{1,w} - \nmu{P^{\perp}_\Delta x}_{1,w} \nn
\\
& = \Re \ip{P_\Delta W v}{\mathrm{sign}(P_\Delta x)} + \nm{x}_{1,w} + \nmu{P^{\perp}_\Delta v }_{1,w}  - 2 \nmu{P^{\perp}_\Delta x}_{1,w}. \label{sign_ineq}
}
Now let $\bar{x} \in \ell^1_w(\bbN)$ be any feasible solution for \R{fin_min}.  Then $\| \hat{x} \|_{1,w} \leq \| \bar{x} \|_{1,w}$ and we get
\bes{
\|  \bar{x} \|_{1,w} \geq  \Re \ip{P_\Delta W v}{\mathrm{sign}(P_\Delta x)} + \nm{x}_{1,w} + \nmu{P^{\perp}_\Delta v }_{1,w}  - 2 \nmu{P^{\perp}_\Delta x}_{1,w}.
}
After rearranging this gives
\be{
\label{2nd_bound}
\| P^{\perp}_\Delta v \|_{1,w} \leq | \ip{P_\Delta W v}{\mathrm{sign}(P_\Delta x)} | + 2 \| P^{\perp}_\Delta x \|_{1,w} + \| x - \bar{x} \|_{1,w}.
}
We next estimate $| \ip{P_\Delta W v}{\mathrm{sign}(P_\Delta x)} |$.  We have
\bes{
| \ip{P_\Delta W v}{\mathrm{sign}(P_\Delta x)} | \leq | \ip{P_\Delta W v}{\mathrm{sign}(P_\Delta x) - P_{\Delta} u} | +  | \ip{W v}{u} | + | \ip{P^{\perp}_{\Delta} W v}{P^{\perp}_{\Delta} u} |.
}
Note that $| \ip{P_\Delta W v}{\mathrm{sign}(P_\Delta x) - P_{\Delta} u} | \leq \gamma \| P_{\Delta} v \|$ by $(iii)$
and also that $\ip{W v}{u} = \ip{v}{W u} = \ip{v}{U^*  u'} = \ip{ U v}{u'}$.
Hence, \R{UN_v} and $(v)$ give
\bes{
|\ip{W v}{u} | \leq \|  U v \| L \sqrt{s} \leq 2 \eta L \sqrt{s}.
}
Finally, by $(iv)$, we have $| \ip{P^{\perp}_{\Delta} W v}{P^{\perp}_{\Delta} u} | \leq \| P^{\perp}_{\Delta} u \|_{\infty} \| P^{\perp}_{\Delta} v \|_{1,w} \leq \theta \| P^{\perp}_{\Delta} v \|_{1,w}$ and therefore
\bes{
| \ip{P_\Delta W v}{\mathrm{sign}(P_\Delta x)} | \leq \gamma \| P_{\Delta} v \| +2 \eta L \sqrt{s} + \theta \| P^{\perp}_{\Delta} v \|_{1,w}.
}
Substituting into \R{2nd_bound} and rearranging now yields
\bes{
(1-\theta) \| P^{\perp}_{\Delta} v \|_{1,w} \leq  \gamma \| P_{\Delta} v \| + 2 \eta L \sqrt{s} + 2 \| P^{\perp}_{\Delta} x \|_{1,w}+ \| x - \bar{x} \|_{1,w},
}
and applying \R{v_Delta_Delta_perp} gives
\bes{
\| P_{\Delta} v \| \leq \frac{\sqrt{1+\alpha}}{1-\alpha} \left [  2 \eta+ \frac{\beta}{1-\theta} \left ( \gamma \| P_{\Delta} v \| + 2 \eta L \sqrt{s}  + 2 \| P^{\perp}_{\Delta} x \|_{1,w} + \| x - \bar{x} \|_{1,w}\right ) \right ].
}
Hence
\eas{
\| P_{\Delta} v \| &\leq \left ( 1 - \frac{\sqrt{1+\alpha} \beta \gamma}{(1-\alpha)(1-\theta)} \right )^{-1}  \frac{2\sqrt{1+\alpha}}{1-\alpha} \left ( 1 + \frac{\beta}{1-\theta} L \sqrt{s} \right )\eta
\\
& +\left ( 1 - \frac{\sqrt{1+\alpha} \beta \gamma q}{(1-\alpha)(1-\theta)} \right )^{-1} \frac{\sqrt{1+\alpha} \beta}{1-\alpha} \left ( 2 \| P^{\perp}_{\Delta} x \|_{1,w} + \| x - \bar{x} \|_{1,w} \right )
\\
& =2 C_0 \left ( 1 + \frac{\beta}{1-\theta} L \sqrt{s} \right )\eta  + C_0 \beta \left ( 2 \| P^{\perp}_{\Delta} x \|_{1,w} + \| x - \bar{x} \|_{1,w} \right ).
}
Since $\| P^{\perp}_{\Delta} v \| \leq \| P^{\perp}_{\Delta} v \|_1 \leq \| P^{\perp}_{\Delta} v \|_{1,w}$ (recall \R{weights_growth}), we now get
\eas{
\| v \| \leq& \| P_{\Delta} v \| + \| P^{\perp}_{\Delta} v \|_{1,w}
\\
 \leq& \left ( 1 + \frac{\gamma}{1-\theta} \right )  \| P_{\Delta} v \| + \frac{2\eta}{1-\theta} L \sqrt{s} + \frac{1}{1-\theta} \left ( 2 \| P^{\perp}_{\Delta} x \|_{1,w} + \| x - \bar{x} \|_{1,w} \right )
\\
\leq& \left [ \left ( 1 + \frac{\gamma}{1-\theta} \right ) C_0 \left (1 + \frac{\beta}{1-\theta} L \sqrt{s} \right )  + \frac{L\sqrt{s}}{1-\theta}  \right ] (2 \eta + \delta ) 
\\
&+ \left [ C_0 \beta  \left ( 1 + \frac{\gamma}{1-\theta} \right ) + \frac{1}{1-\theta} \right ]  \left ( 2 \| P^{\perp}_{\Delta} x \|_{1,w} + \| x - \bar{x} \|_{1,w} \right ),
}
as required.
}

\section{Handling truncation: the choice of $K$}\label{s:truncation}
We now consider the truncation parameter $K$.  Due to Theorem \ref{t:full_samp_recov}, it suffices to estimate the quantity $T_{h,K,\eta}(x)$ defined in \R{trunc}.

\thm{
\label{p:trunc_err}
For all sufficiently large $K$ we have $\mathrm{Ran}(U) = \mathrm{Ran}(U P_K)$.  In particular, the problems \R{fin_min_noiseless} and \R{fin_min} have solutions for all large $K$.  Moreover, suppose that $r = \mathrm{rank}(U) \leq N$ and $K$ is sufficiently large so that $\mathrm{Ran}(U P_K) = \mathrm{Ran}(U)$.  If $x \in \ell^1_w(\bbN)$ then
\bes{
T_{h,K,\eta}(x) \leq \| x - P_K x \|_{1,w} + \| P_K w \| /\sigma_{\min} \| x - P_K x \|_{1,w},
}
where $\sigma_{\min}$ is the minimum singular value of $U P_K$.  Moreover, if $\{ w_i \}_{i \in \bbN}$ is nondecreasing, and $x \in \ell^1_{\tilde{w}}(\bbN)$, where $\tilde{w} = \{ \tilde{w}_i \}_{i \in \bbN}$ with $\tilde{w}_i \geq \sqrt{i} w^2_{i}$, $\forall i \in \bbN$, then
\bes{
T_{h,K,\eta}(x) \leq \| x - P_K x \|_{1,w} + 1/\sigma_{\min} \| x - P_K x \|_{1,\tilde{w}}.
}
}
\prf{
The first observation is immediate since $U$ has finite rank.  Suppose now that $K$ is such that $\mathrm{Ran}(U P_K) = \mathrm{Ran}(U)$ and write $\bar{x} = P_K x +  (U P_K)^{\dag} U (x-P_K x)$.  Then
\bes{
\| U P_K \bar{x} - y \| = \| U P_K x +  U(x-P_K x) - U x - \{ \sqrt{\tau_n} e_n \}^{N}_{n=1} \| \leq  \eta.
}
Hence $\bar{x}$ is feasible.  Moreover, $\| x - \bar{x} \|_{1,w} \leq \| x - P_K x \|_{1,w} + \| (U P_K)^{\dag} U (x-P_K x) \|_{1,w}$ and
\eas{
\| (U P_K)^{\dag} U (x-P_K x) \|^2_{1,w} & \leq \| P_K w \|^2  \| (U P_K)^{\dag} U (x-P_K x) \|^2
\\
& \leq \| P_K w \|^2 \| U(x-P_K x) \|^2 / \sigma^2_{\min}
\\
& \leq \| P_K w \|^2 \| x - P_K x \|^2_{1,w} / \sigma^2_{\min}.
}
To obtain the final result, we note that
\bes{
\| P_K w \| \| x - P_K x \|_{1,w} \leq w_{K} \sqrt{K} \sum_{i > K} w_i | x_i | \leq \sum_{i > K} \tilde{w}_i | x_i| = \| x - P_K x \|_{1,\tilde{w}},
}
whenever the $w_i$'s are nondecreasing, as required.
}

Note that it need not be the case that \R{fin_min_noiseless} or \R{fin_min} have solutions for arbitrary $K \geq N$, since $\mathrm{Ran}(U P_K) \neq \mathrm{Ran}(U)$ in general.  However, this theorem shows that this holds for all large $K$.  Furthermore, this result shows that once $K$ is chosen so that $1/\sigma_{\min}$ is moderate in size, the effect of truncation is bounded by the decay of the coefficients $x_i$, $i > K$.

\rem{
\label{r:trunc_weak}
Theorem \ref{p:trunc_err} asserts that the truncation parameter $K$ can be chosen independently of the coefficients whenever $x \in \ell^1_{\tilde{w}}(\bbN)$. While it is possible to show that $T_{h,K,\eta}(x) \rightarrow 0$ as $K \rightarrow \infty$ for $x \in \ell^1_{w}(\bbN)$ \cite[Prop.\ 6.6]{BAACHGSCS}, it is currently unknown whether a bound for $T_{h,K,\eta}(x)$ involving only $\| x - P_{K} x \|_{1,w}$ holds.  In other words, if $x \in \ell^1_{w}(\bbN)$ but $x \notin \ell^1_{\tilde{w}}(\bbN)$ the truncation strategy may no longer be independent of $x$.  Improving this result is an open problem.
}

It is important to quantify precisely how large $K$ needs to be in relation to $N$ to ensure a small truncation error.  For this, we shall use the following lemma:

\lem{
\label{l:gamma}
The minimum singular value $\sigma_{\min}$ of $U P_K$ satisfies
\be{
\label{gamma_def}
\sigma_{\min} \geq  \inf_{\substack{y \in \bbC^N \\ \| y \| =1}} \sup_{g \in G_y} \sup_{\substack{\phi \in \Phi_K \\ \phi \neq 0}}  \left \{ \frac{1-\| g - \phi \|_{\infty}}{\| g \|_{\nu} + \| g - \phi \|_{\infty}} \right \},
}
where $\Phi_K = \spn \{ \phi_1,\ldots,\phi_K \}$ and $G_y = \left \{ g \in L^2_{\nu}(D) \cap L^\infty(D) : g(t_n) = y_{n} / \sqrt{\tau_n},\ n=1,\ldots,N \right \}$.
}
\prf{
The minimum singular value is given by $\sigma_{\min} = \inf_{\substack{y \in \bbC^N \\ \| y \| =1}} \| (U P_K)^* y \|$.  Fix $y \in \bbC^N$, $\| y \| =1$ and observe that
\bes{
\| (U P_K)^* y \| = \sup_{\substack{z \in \bbC^K \\ \| z \|=1}} \left | \ip{y}{U P_K z} \right |= \sup_{\substack{\phi \in \Phi_K \\ \| \phi \|_{\nu} = 1}} \left | \sum^{N}_{n=1} \sqrt{\tau_n} y_n \overline{\phi(t_n)} \right |.
}
Let $g \in G_y$ and $\phi \in \Phi_K$.  Then
\eas{
\| (U P_K)^* y \| &\geq \left | \sum^{N}_{n=1} \sqrt{\tau_n} y_n \overline{\phi(t_n) } \right |  \geq \| y \|^2 - \left | \sum^{N}_{n=1} \sqrt{\tau_n} y_n \overline{g(t_n) - \phi (t_n)} \right |
 \geq \| y \|^2 - \| y \| \| g - \phi \|_{L^\infty}.
}
Since $\| \phi \|_{\nu} \leq \| g \|_{\nu} +\| g - \phi \|_{\nu} \leq \| g \|_{\nu} + \| g - \phi \|_{\infty}$ the result now follows immediately.
}

\rem{
\label{r:just_compute}
In practice, rather than performing an analysis of $\sigma_{\min}$ via Lemma \ref{l:gamma}, one may choose $K$ simply by calculating the minimal singular value of $UP_K$.  If this is sufficiently large, then Theorem \ref{p:trunc_err} guarantees that the truncation error is moderate.
}

\subsection{Examples}\label{ss:trunc_examples}
Lemma \ref{l:gamma} allows one to determine the required condition on $K$.  Note that this depends completely on the points $T$ and the basis $\{ \phi_i \}_{i \in \bbN}$.  We now do this for the two examples of \S \ref{ss:examples}:

\thm{
\label{t:Leg_trunc}
For $\alpha,\beta > -1$ let $\{ \phi_i \}_{i \in \bbN}$ be the Jacobi polynomial basis \R{Jacobi_ON} and let $T = \{ t_n \}^{N}_{n=1}$ be an ordered set of points in $[-1,1]$.  Then for every $r \in \bbN$ there exists a $C_r >0$ such that
\bes{
\sigma_{\min} \geq 1 -  \frac{C_r K^{-r} \xi^{-r-1/2}}{\sqrt{\xi / (2 h_N)}-C_r K^{-r} \xi^{-r-1/2}},
}
where $\sigma_{\min}$ is the minimal singular value of $U P_K$,
\bes{
\xi = \frac12\min_{n=0,\ldots,N} \{ t_{n+1}-t_n \},\qquad h = \sup_{-1 \leq t \leq 1} \min_{n=1,\ldots,N} | t - t_n |,
}
and $t_0 = -t_1-2$, $t_{N+1} = 2-t_N$.  In particular, if 
\bes{
K \geq \left ( \sqrt{2}C_r(1+\epsilon^{-1})\right )^{\frac{1}{r}} h^{\frac{1}{2r}} \xi^{-1-\frac1r},
}
for some $0 < \epsilon < 1$ then $\sigma_{\min} > 1- \epsilon$.
}

We defer the proof of this result until \S \ref{ss:Leg_proofs}. Note that in the case of equispaced data, we have $h = \xi = 1/N$ and so it suffices to take $K \gtrsim N^{1+\frac{1}{2r}}$ for any $r \in \bbN$.  In practice, we have found that $K=4N$ is sufficient in all examples (recall also Remark \ref{r:just_compute}).  On the other hand, if the data clusters severely, then this theorem suggests that a larger value of $K$ may be necessary.

\thm{
\label{t:Trig_trunc}
Let $\{ \phi_i\}_{i \in \bbZ}$ be the Fourier basis \R{fourier_basis} and let $T = \{ t_n \}^{N}_{n=1}$ be a set of $N$ ordered points in $[-1,1]$.  Then for every $r \in \bbN$ there exists a $C_r >0$ such that
\bes{
\sigma_{\min} \geq 1 - \frac{C_r K^{-r} \xi^{-r-1/2}}{\sqrt{\xi / (2 h_N)}-C_r K^{-r} \xi^{-r-1/2}},
}
where $\sigma_{\min}$ is the minimal singular value of $U P_K$,
\bes{
\xi = \frac12\min_{n=0,\ldots,N} \{ t_{n+1}-t_n \},\qquad h = \sup_{-1 \leq t \leq 1} \min_{n=1,\ldots,N} | t - t_n |,
}
and $t_0 = -1$, $t_{N+1} = 1$.  In particular, if 
\bes{
K \geq \left ( \sqrt{2}C_r(1+\epsilon^{-1})\right )^{\frac{1}{r}} h^{\frac{1}{2r}} \xi^{-1-\frac1r},
}
for some $0 < \epsilon < 1$ then $0 \leq \gamma < \epsilon$.
}

This result is exactly the same as that for Jacobi polynomials (Theorem \ref{t:Leg_trunc}), except up to a minor change in the definition of the values $t_0$ and $t_{N+1}$, and therefore $\xi$.  Its proof is near identical, and hence is omitted.

\section{Jacobi polynomials on the unit interval}\label{s:AlgPoly_Examp}

We now consider Example \ref{ex:Jacobi}.  For convenience, we recall the growth condition \R{Jacobi_ON_growth}:
\be{
\label{Jacobi_ON_growth2}
\| \phi_j \|_{L^\infty} = \ord{j^{q+1/2}},\quad j \rightarrow \infty,\qquad \mbox{where $q = \max \{ \alpha,\beta , -1/2 \}$}.
}
Our main result is the following:

\thm{
\label{t:Leg_poly_full_l1}
For $\alpha,\beta > -1$ let $\{ \phi_i \}_{i \in \bbN}$ be the orthonormal Jacobi polynomial basis \R{Jacobi_ON}, $T = \{ t_n \}^{N}_{n=1} \subseteq [-1,1]$ be a set of scattered points and suppose that $h$ is as in \R{h_def}.  Suppose that the weights $w_i = i^{\gamma} \| \phi_i \|_{L^{\infty}}$ for some $\gamma > 0$.  Then for $0 < \epsilon < 1/2$ there exists a $c(\epsilon) > 0$ such that if
\be{
\label{hcondPoly}
h \leq c(\epsilon) \left \{ \begin{array}{cl} \frac{1}{M^{2+2(q+1)/\gamma} \log M} & 0 < \gamma \leq 1/2-q \\ \frac{1}{M^2 \log M} & \gamma> 1/2-q \end{array} \right .,
}
where $q$ is given by \R{Jacobi_ON_growth2}, any minimizer $\hat{x}$ of \R{fin_min} satisfies
\bes{
\| x - \hat{x} \| \leq C(\epsilon) \left ( M^{\gamma+q+1}   \eta  + \| x -P_M x \|_{1,w} + T_{N,K, \eta}(x) \right ),
}
for some constant $C$ depending on $\epsilon$ only, where $T_{N,K, \eta }(x)$ is as in \R{trunc}.
}

One also has a similar, albeit simpler, result for least squares:

\thm{
\label{t:Jacobi_poly_LS}
Let $\{ \phi_i \}_{i \in \bbN}$ be the orthonormal Jacobi polynomial basis \R{Jacobi_ON}, $T = \{ t_n \}^{N}_{n=1} \subseteq [-1,1]$ be a set of scattered points and suppose that $h$ is as in \R{h_def}.  Then for each $0 < \epsilon < 1$ there exists a $c(\epsilon) > 0$ such that if
\be{
\label{hcondPolyLS}
h \leq c(\epsilon) M^{-2},
}
then the solution $\check{x}$ of \R{LS_fit} exists uniquely and satisfies
\bes{
\| x - \check{x} \| \leq \left ( 1 + \frac{1}{\sqrt{1-\epsilon}} \right ) \| x - P_M x \|_{1,w} + \frac{1}{\sqrt{1-\epsilon}}  \eta,
}
for any $w = \{ w_i \}_{i \in \bbN}$ with $w_i \geq \| \phi_i \|_{L^\infty}$.
}

Theorems \ref{t:Leg_poly_full_l1} and \ref{t:Jacobi_poly_LS} assert that both techniques guarantee an approximation error that depends on $x - P_M x$ measured in some norm, for the same asymptotic scaling of $h$ with $M$ up to log factors.  Hence, up to the choice of norm, weighted $\ell^1$ minimization with scattered data, Jacobi polynomials and sufficiently large weights $w_i$ is guaranteed a similar approximation rate as least-squares fitting.

It is informative to now consider the following two cases:

\pbk
\textit{Smooth functions.} Let $f \in C^{\infty}([-1,1])$.  Then the coefficients $x_{i} = \ord{i^{-k}}$ as $i \rightarrow \infty$ for any $k >0$.  Hence $\| x - P_M x \|_{1,w} = \ord{M^{-k}}$ as $M \rightarrow \infty$ for any $k >0$ whenever the weights $w_i$ grow at most algebraically fast in $i$.  By Theorems \ref{t:Leg_poly_full_l1} and \ref{t:Jacobi_poly_LS} the approximation errors for weighted $\ell^1$ minimization and least squares both decay superalgebraically fast in $h$ as $h \rightarrow 0$; that is, faster than $h^k$ for any $k > 0$.

\pbk
\textit{Analytic functions.}  Suppose $f$ is analytic so that $x_i = \ord{\rho^{-i}}$ for some $\rho > 1$.  For algebraic weights $w_i$ one has $\| x - P_M x \|_{1,w} = \ord{(\rho')^{-M}}$ as $M \rightarrow \infty$ for any $\rho' < \rho$.  Therefore the approximation error for least squares behaves like $\| x - \check{x} \| = \ordu{(\rho')^{-1/\sqrt{h}}}$ as $h \rightarrow 0$,
and for weighted $\ell^1$ minimization one has the marginally slower convergence $\| x - \hat{x} \| = \ordu{(\rho')^{-1/\sqrt{h \log(h)}}}$, provided the weights $w_{i} =  i^{\gamma} \| \phi_i \|_{L^\infty}$ with $\gamma > 1/2 - q$.

\rem{
On the other hand, for functions of finite smoothness, the need for more rapidly-growing weights in Theorem \ref{t:Leg_poly_full_l1} translates into slower algebraic convergence of the approximation than that of least squares.  We expect that this may be an artifact of the proof, however.
}

Suppose now that the data $T$ is equispaced.  We first note the following general result:

\thm{
\label{t:PTK}
Let $T = \{ t_n \}^{N}_{n=1}$ be an equispaced grid of $N$ points in $[-1,1]$, $E \subseteq \bbC$ be a compact set containing $[-1,1]$ in its interior and let $B(E)$ be the Banach space of functions continuous on $B$ and analytic in its interior, with norm $\| f \|_{B} = \sup_{z \in B} | f(z) |$.  Let 
$F : B(E) \rightarrow L^\infty(-1,1)$ be a mapping such that for each $f \in B(E)$, $F(f)$ depends only on the data $\{ f(t_n) \}^{N}_{n=1}$, and suppose that, for constants $C>0$, $\sigma > 1$ and $1/2 < \tau \leq 1$,
\bes{
\| f - F(f) \|_{L^\infty} \leq C \sigma^{-N^{\tau}} \| f \|_{E},\quad \forall f \in B(E).
}
If $\| f \|_{T,\infty} = \max_{n=1,\ldots,N} | f(t_n) |$ then there exists a constant $\nu > 1$ such that
\be{
\label{condition_replacement}
\Theta(F) = \sup_{\substack{f \in B(E) \\ \| f \|_{T,\infty} \neq 0} } \left \{ \frac{\| F(f) \|_{L^\infty}}{\| f \|_{T,\infty}} \right \} \geq \nu^{N^{2 \tau - 1}}.
}
}
This theorem is due to Platte, Trefethen \& Kuijlaars \cite{TrefPlatteIllCond} (a minor modification is made in \R{condition_replacement} which is more suitable for our purposes).  It states the following: for any method that achieves an error for all functions in $B(E)$ that is exponentially-decaying as $N \rightarrow \infty$ with rate $\tau$ it is possible to find a function $f \in B(E)$ which is bounded on the set $T$, but for which $\|F(f)\|_{L^\infty}$ is exponentially large with index $2 \tau - 1$.  In particular, the best possible convergence rate for a robust method, i.e.\ one for which $\Theta(F) \leq C$ for all $N \in \bbN$, is root-exponential in $N$.  Note that this theorem is very general: the method $F$ can be linear or nonlinear, and $F(f)$ only needs to be defined for extremely smooth (specifically, analytic) functions.

Now consider the cases of weighted $\ell^1$ minimization and least-squares fitting with Jacobi polynomials.  By earlier arguments, the corresponding approximation errors behave like $\ordu{(\rho')^{-\sqrt{N}}}$ and $\ordu{(\rho')^{-\sqrt{N/\log(N)}}}$ respectively as $N \rightarrow \infty$.  Moreover, by setting $x = 0$ in Theorems \ref{t:Leg_poly_full_l1} and \ref{t:Jacobi_poly_LS} respectively, one deduces that $\Theta(F) \lesssim 1$ for the former and $\Theta(F) \lesssim (N/\log(N))^{(\gamma+q+1)/2}$ for the latter.  According to Theorem \ref{t:PTK}, least-squares fitting attains the best possible convergence rate for a robust method, while weighted $\ell^1$ minimization (with sufficiently large weights) attains nearly the best possible convergence rate, with only slow, algebraic growth of $\Theta(F)$.

\rem{
Although Theorem \ref{t:PTK} applies only to equispaced data, it can also be formulated for much more general data.  Loosely speaking, similar conclusions apply unless the data clusters quadratically at the endpoints $x = \pm 1$ \cite{AdcockNecSamp}.
}

\rem{
\label{noise_growth}
For general data, the constant $M^{\gamma+q+1}$ in Theorem \ref{t:Leg_poly_full_l1} implies mild ill-conditioning of weighted $\ell^1$ minimization as $h \rightarrow 0^{+}$.  This is not seen in computations, and we expect it is also an artifact of the proof.  Removing this factor is an open problem.
}

\subsection{Numerical examples}\label{ss:Polynumexamp}

\begin{figure}
\begin{center}
$\begin{array}{ccc}
\includegraphics[width=4.5cm]{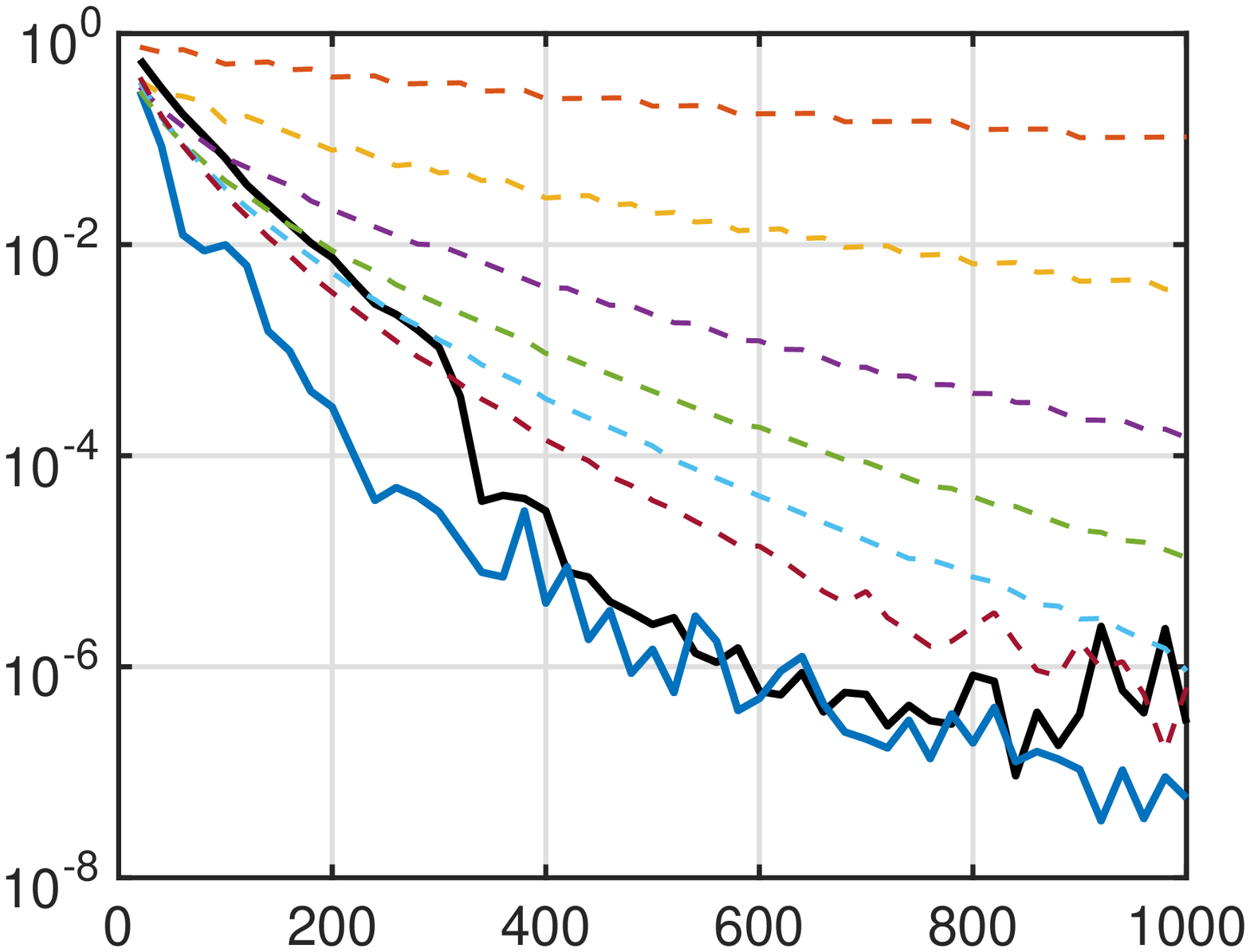}   
&
\includegraphics[width=4.5cm]{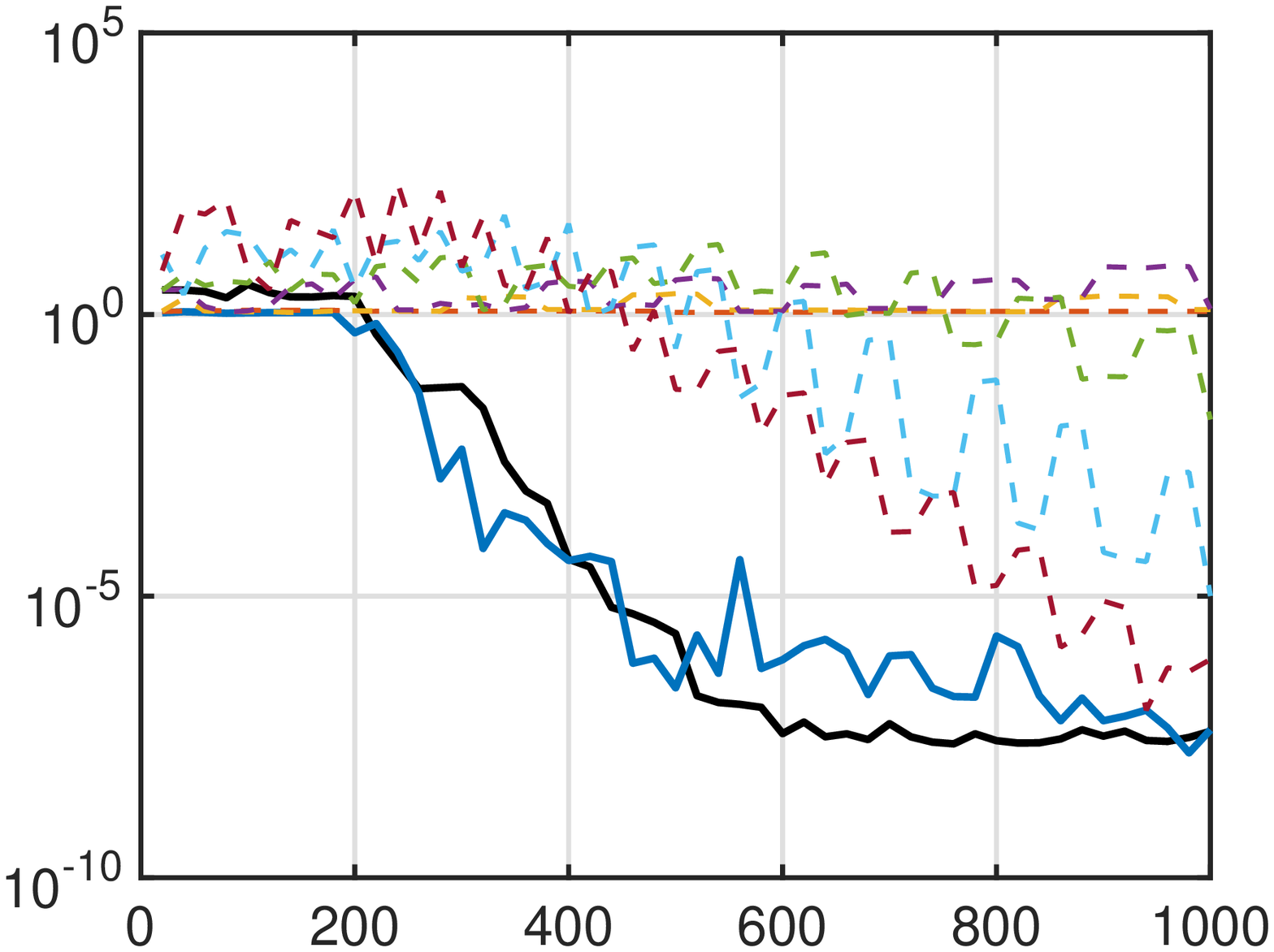}   
&
\includegraphics[width=4.5cm]{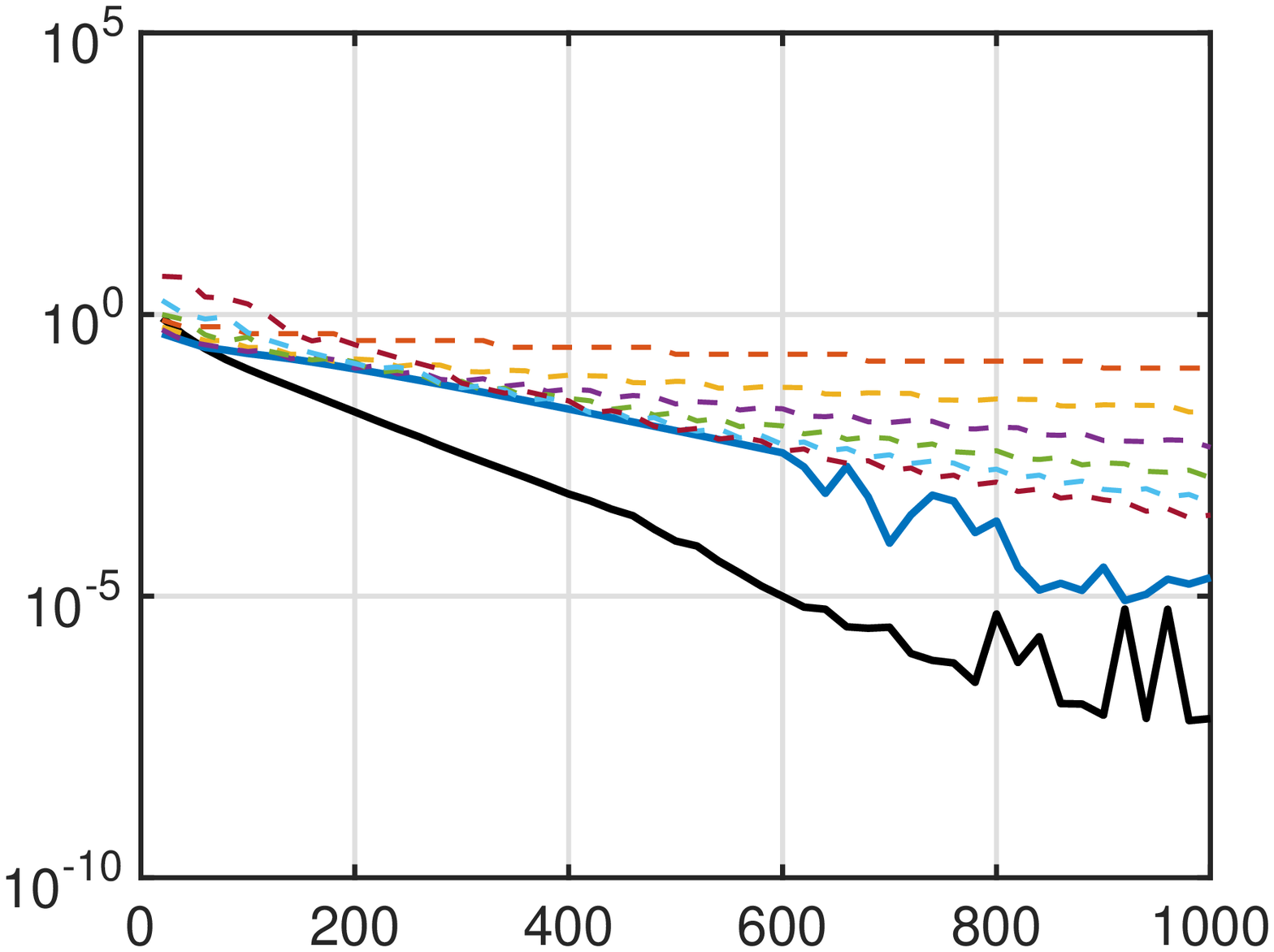}    
\\ 
f(x) = \frac{1}{50/49-\sin(\pi x)} & f(x) = \sin(50 x^2)  &  f(x) = \frac{1}{1+50 x^2}
\\ \\
\includegraphics[width=4.5cm]{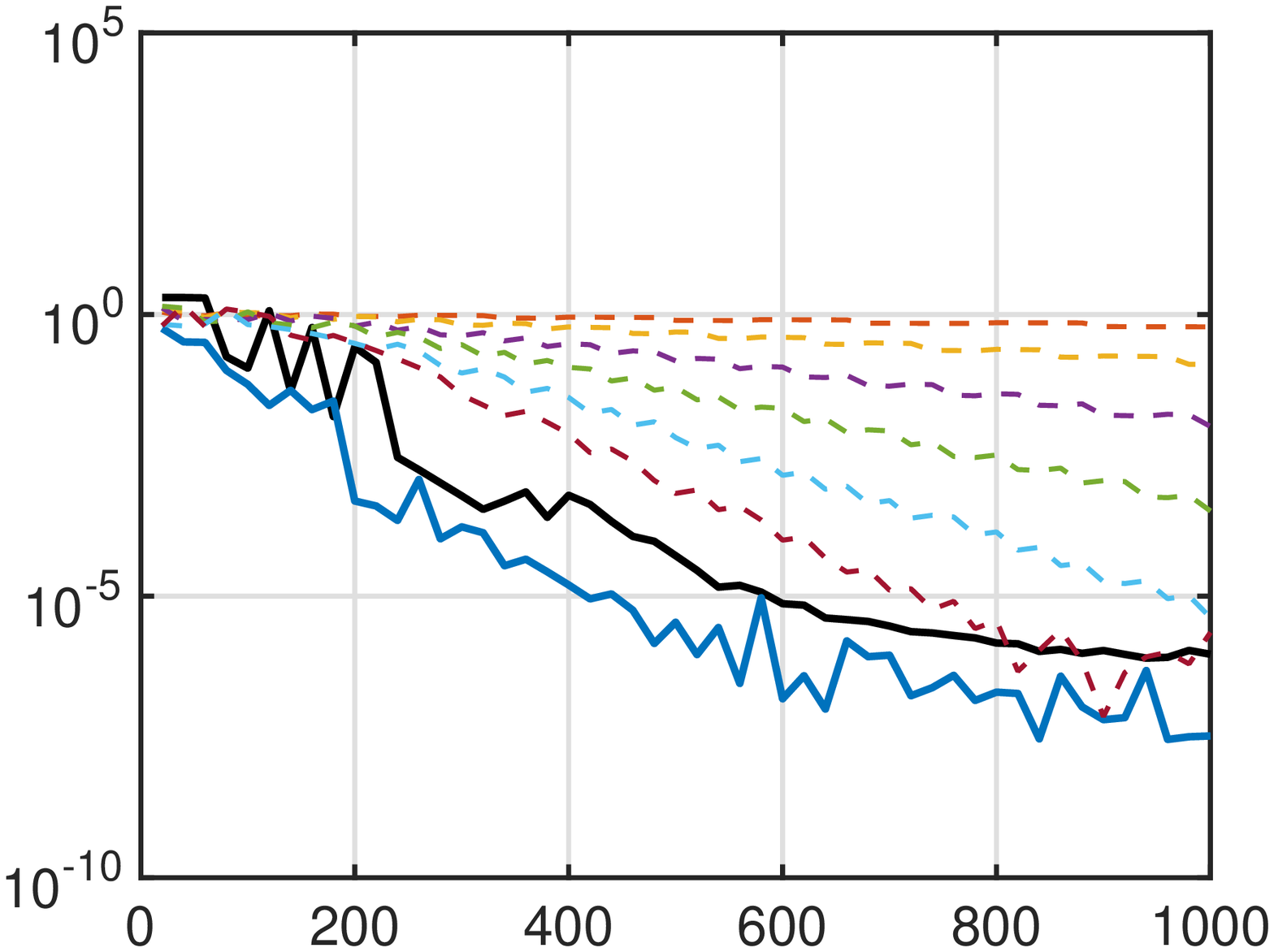}   
&
\includegraphics[width=4.5cm]{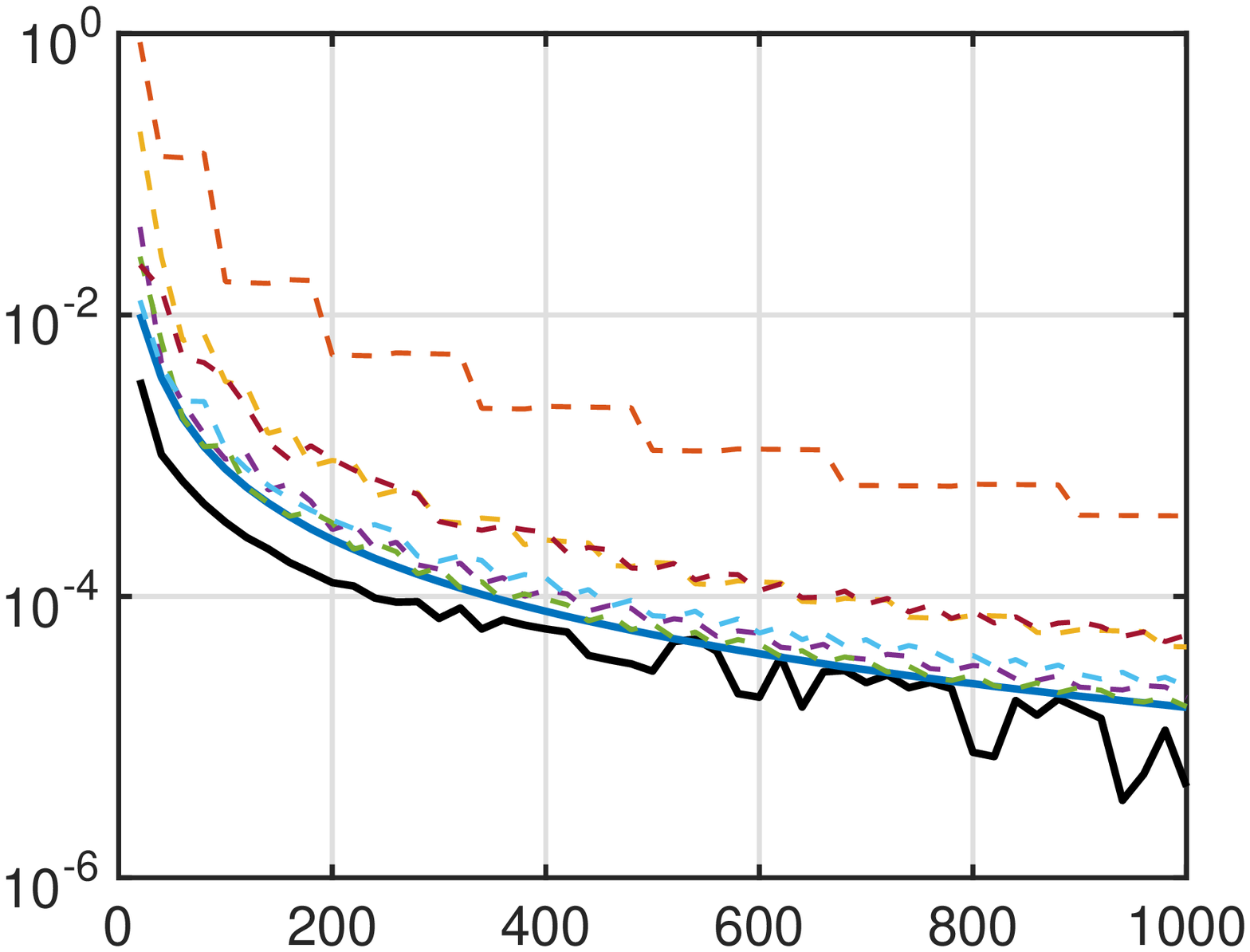}   
&
\includegraphics[width=4.5cm]{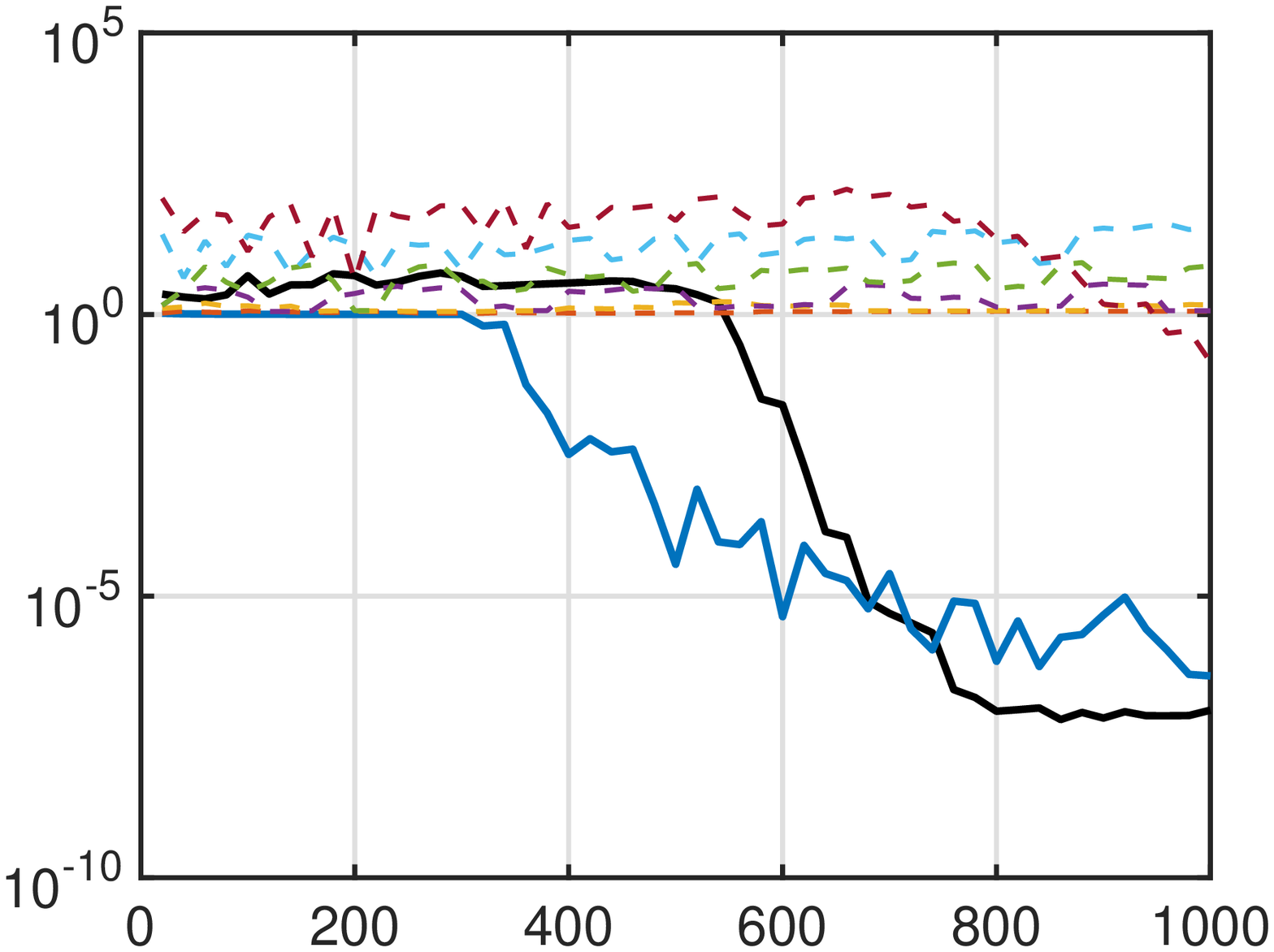}    
\\ 
f(x) = \frac{\cosh(100 x^2)}{\cosh(100)} &  f(x) = |x|^3 &  f(x) = \sin(80 x)
\end{array}$
\caption{
Numerical comparison of weighted $\ell^1$ minimization and least-squares fitting for approximation from equispaced data using Legendre polynomials.  The error against $N$ is plotted for each method.  The solid black line is weighted $\ell^1$ minimization with $K=4N$ and weights $w_{i} = i$.  The dashed lines are least squares with $M = c \sqrt{N}$ and $c=0.5,1.0,1.5,2,2.5,3.0$.  The solid blue line is oracle least squares based on choosing $M$ between $1$ and $N$ which minimizes the error $\| f - \tilde{f} \|_{L^\infty}$ for a given $N$ and $f$.  Random noise of magnitude $10^{-8}$ was added to the data.
}
\label{f:LegEqui}
\end{center}
\end{figure}

\begin{figure}
\begin{center}
$\begin{array}{ccc}
\includegraphics[width=4.5cm]{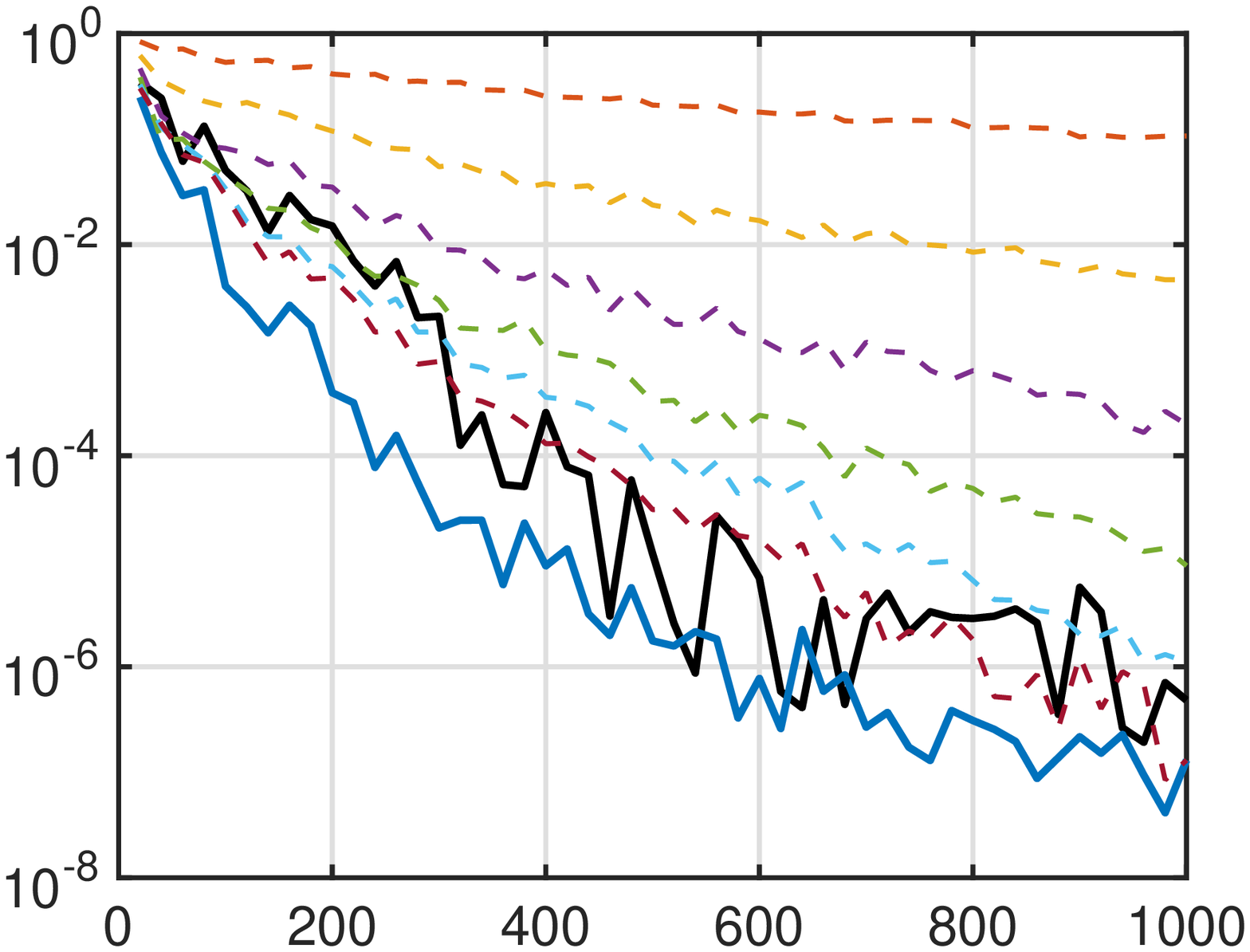}   
&
\includegraphics[width=4.5cm]{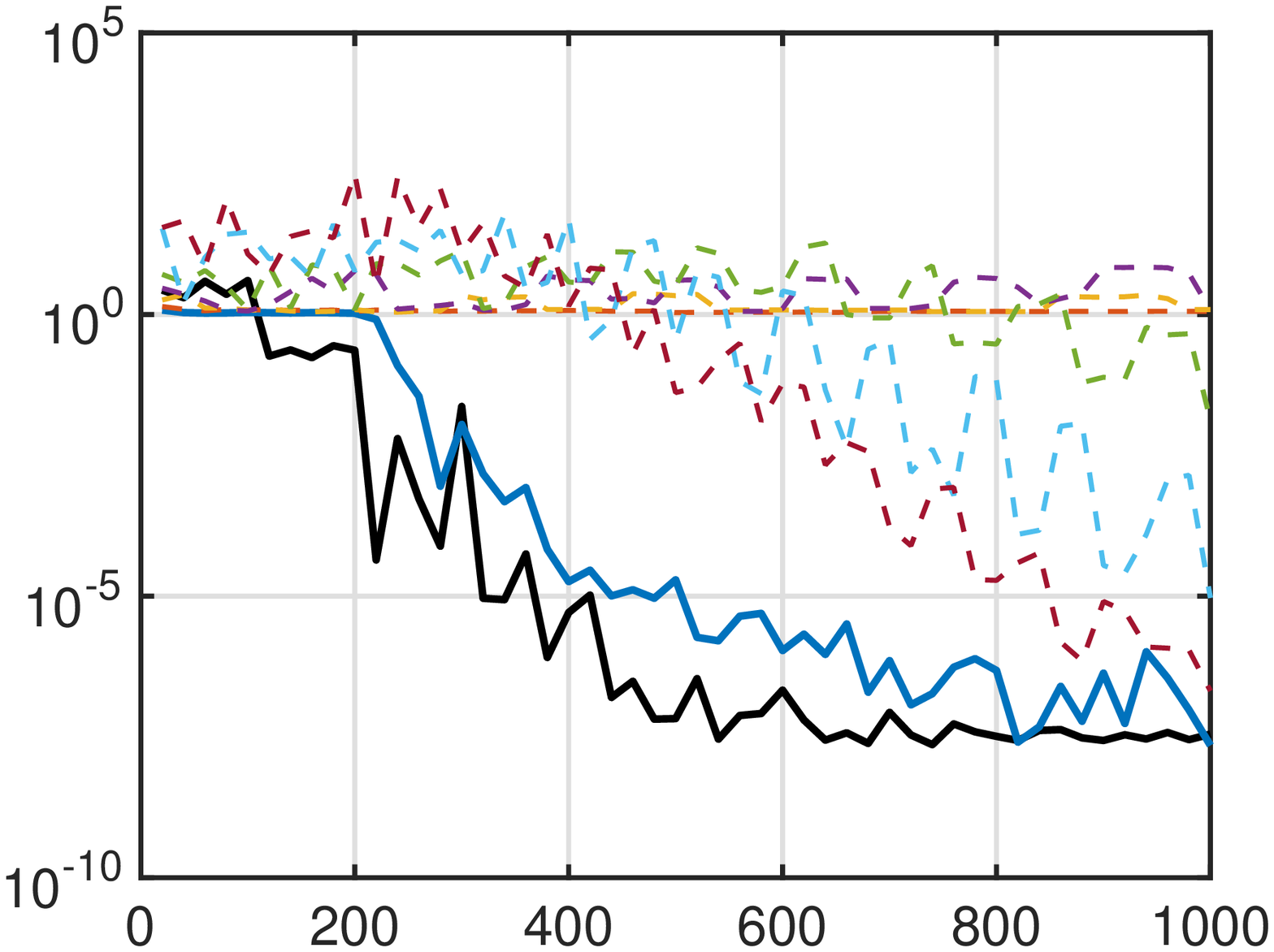}   
&
\includegraphics[width=4.5cm]{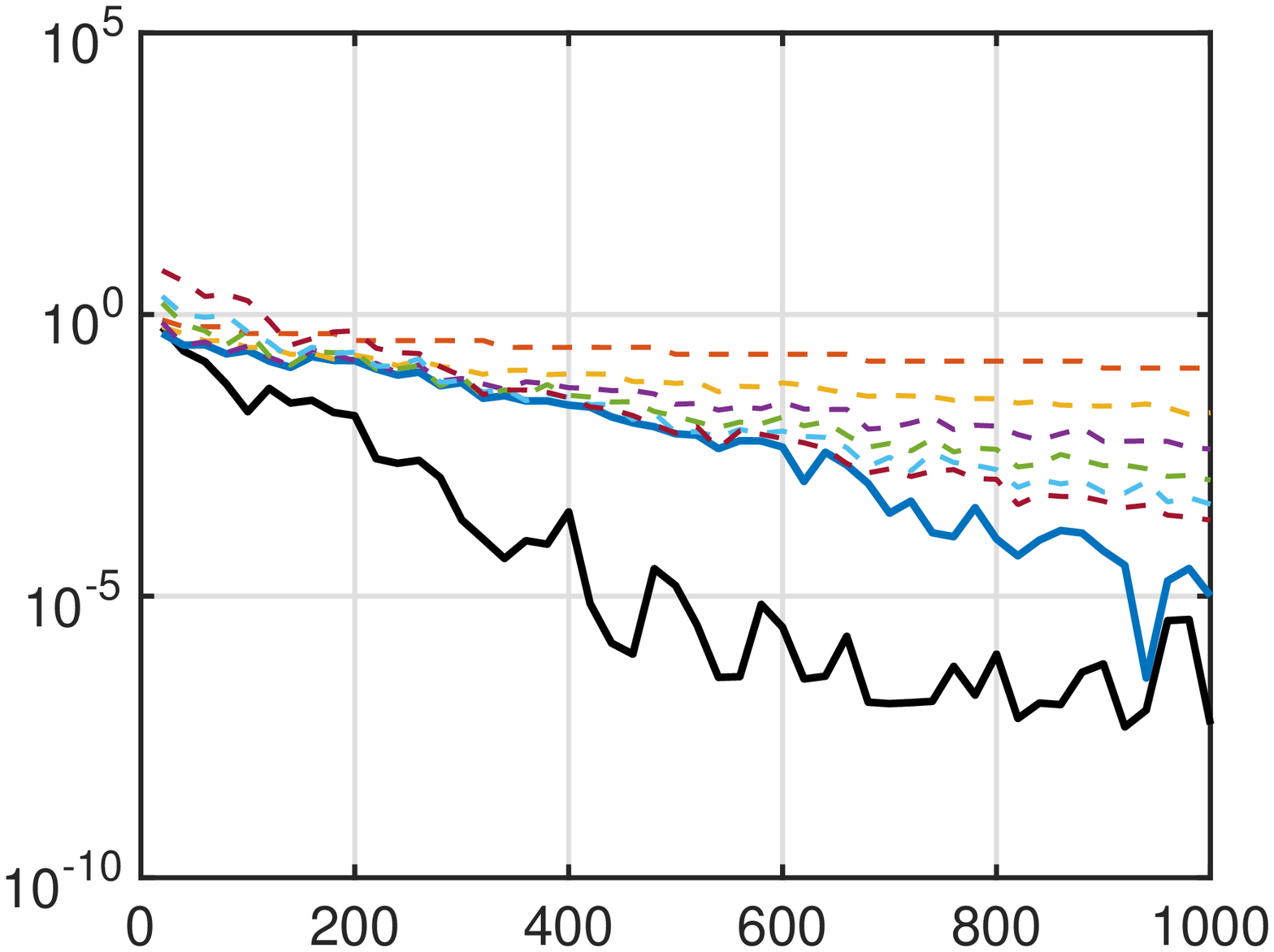}    
\\ 
\mbox{\scriptsize $f(x) = \frac{1}{50/49-\sin(\pi x)}$} & \mbox{\scriptsize $f(x) = \sin(50 x^2) $} & \mbox{\scriptsize $f(x) = \frac{1}{1+50 x^2}$}
\\ \\
\includegraphics[width=4.5cm]{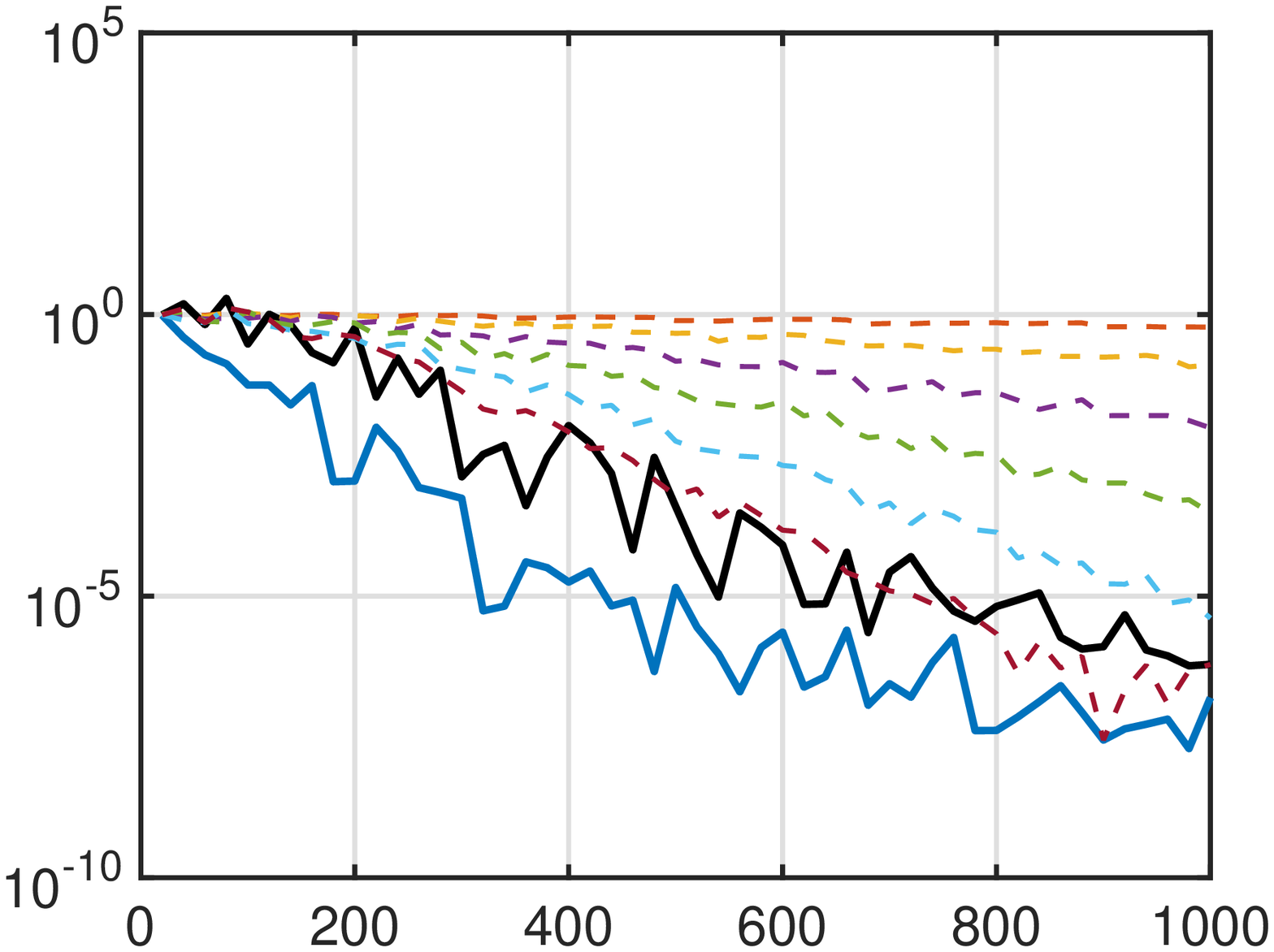}   
&
\includegraphics[width=4.5cm]{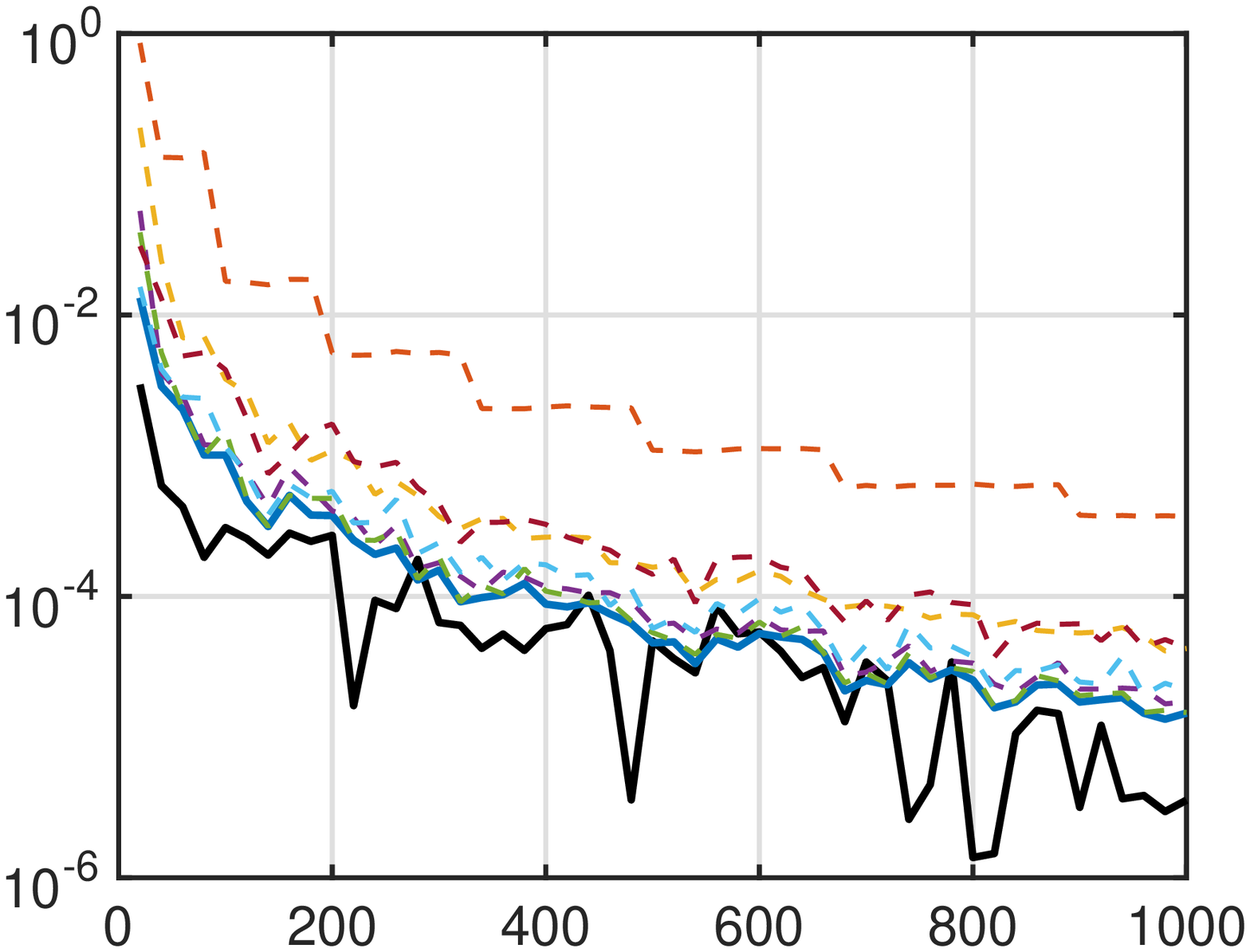}   
&
\includegraphics[width=4.5cm]{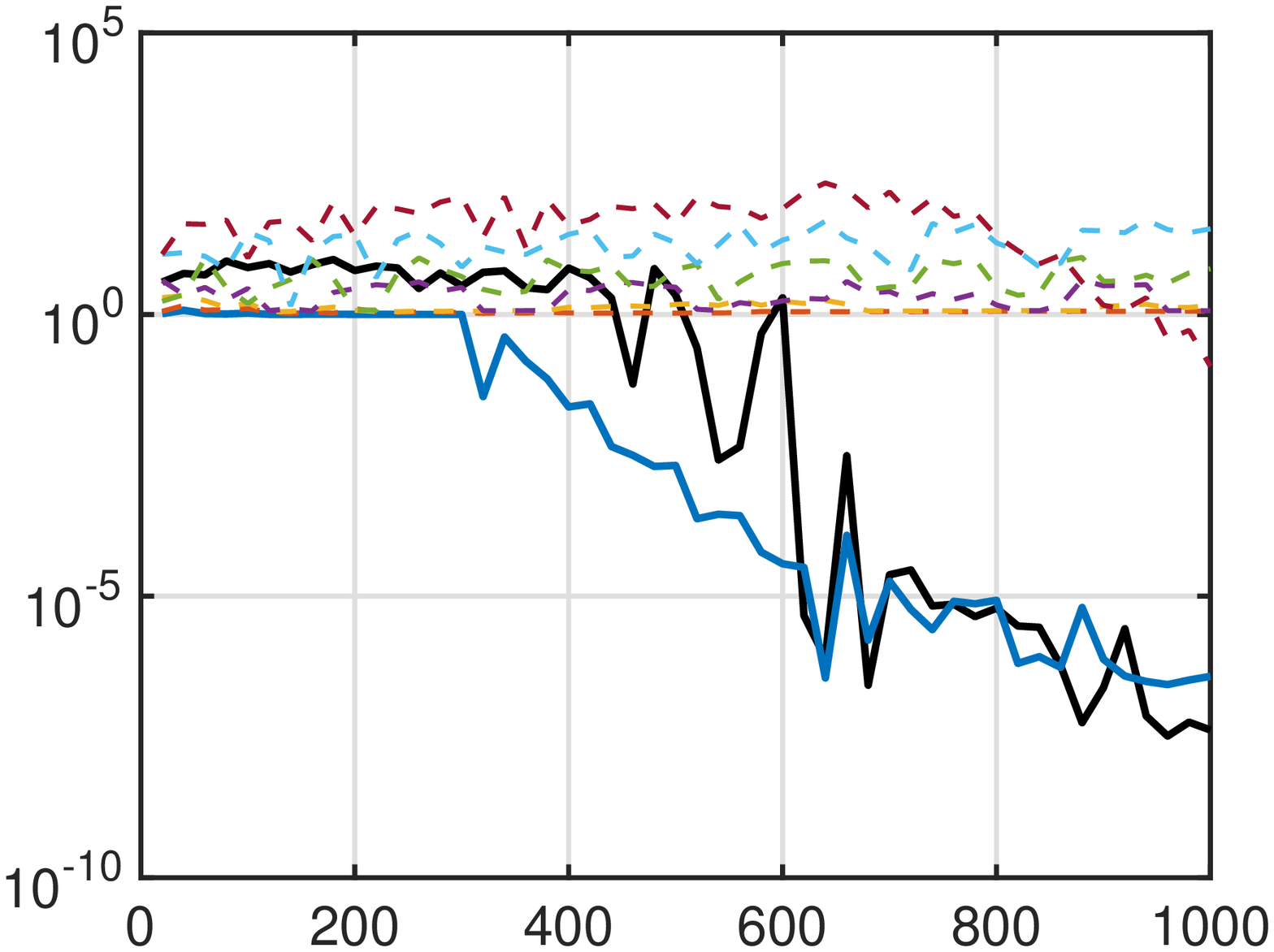}    
\\ 
\mbox{\scriptsize $f(x) = \frac{\cosh(100 x^2)}{\cosh(100)}$} & \mbox{\scriptsize $f(x) = |x|^3$} & \mbox{\scriptsize $f(x) = \sin(80 x)$}
\end{array}$
\caption{
The same as Fig.\ \ref{f:LegEqui} except for randomly jittered data. 
}
\label{f:LegJitt}
\end{center}
\end{figure}

In Figs.\ \ref{f:LegEqui} and \ref{f:LegJitt} we give a results for polynomial approximations from equispaced and jittered data.  Although it has only been proved that weighted $\ell^1$ minimization performs as well (up to a log factor) as least-squares fitting, these results show that it in fact exhibits rather better performance, similar to that of the best possible least-squares fit.  Note that this oracle least squares cannot be implemented in practice since $M$ is calculated by minimizing the approximation error.

\subsection{Proofs}\label{ss:Leg_proofs}
The proof of Theorem \ref{t:Leg_poly_full_l1} relies on the following three lemmas, which provide estimates for the quantities $E_2(h,M)$, $E_{\infty}(h,M)$ and $F(h,M,R)$ respectively.

\lem{
\label{l:E2_Leg}
For $\alpha,\beta > -1$ let $\{ \phi_i \}_{i \in \bbN}$ be the orthonormal Jacobi polynomial basis \R{Jacobi_ON}, $T = \{ t_n \}^{N}_{n=1} \subseteq D$ be a set of scattered data points and suppose that $h$ is as in \R{h_def}.  If $h M^2 \leq 1$ then $E_{2}(h,M) \lesssim \sqrt{h} M$, where $E_{2}(h,M)$ is as in \R{E_2NM}.
}

\prf{
By self-adjointness,
\be{
\label{norm_identity}
\| P_M - P_M U^* U P_M \| = \sup_{\substack{x \in P_M(\ell^2(\bbN)) \\ \| x \|=1}} | \ip{(P_M - P_M U^* U P_M) x}{x} |.
}
Let $x \in P_M(\ell^2(\bbN))$, $\| x \| =1$ be arbitrary and set $g = \sum^{M}_{j=1} x_j \phi_j \in \bbP_{M-1}$ so that $\| g \|_{L^2_{\nu^{(\alpha,\beta)}}} = 1$.  Let $\{ V_n \}^{N}_{n=1}$ be the Voronoi cells of the points $\{ t_n \}^{N}_{n=1}$ and set $\chi(t) =\sum^{N}_{n=1} g(t_n) \bbI_{V_n}(t)$.
Then, by the definition \R{quad_weights_1} of the weights $\tau_n$, we have  $\| \chi \|^2_{L^2_{\nu^{(\alpha,\beta)}}} = \sum^{N}_{n=1} \tau_n | g(t_n) |^2$, and therefore
\bes{
| \ip{(P_M - P_M U^* U P_M) x}{x} | = \left | 1 - \sum^{N}_{n=1}\tau_n | g(t_n) |^2 \right | = \left | 
\| g \|^2_{L^2_{\nu^{(\alpha,\beta)}}} - \| \chi \|^2_{L^2_{\nu^{(\alpha,\beta)}}} \right |.
}
Hence 
\be{
\label{diff_ineq}
| \ip{(P_M - P_M U^* U P_M) x}{x} | \leq \| g - \chi \|_{L^2_{\nu^{(\alpha,\beta)}}} \left ( 2 \| g \|_{L^2_{\nu^{(\alpha,\beta)}}} + \| g - \chi \|_{L^2_{\nu^{(\alpha,\beta)}}} \right ),
}
and so it suffices to show that
\be{
\label{desired_result}
 \| g - \chi \|_{L^2_{\nu^{(\alpha,\beta)}}} \lesssim \sqrt{h} M \| g \|_{L^2_{\nu^{(\alpha,\beta)}}}.
 }
We have
\bes{
\| g - \chi \|^2_{L^2_{\nu^{(\alpha,\beta)}}} = \sum^{N}_{n=1} \int_{V_n} \left | g(t) - g(t_n) \right |^2 \nu^{(\alpha,\beta)}(t) \D t = \sum^{N}_{n=1} I_n.
}
Let $n_0$ be the unique number such that $0 \in V_{n_0}$.  Then we write this as 
\be{
\label{int_sum_split}
\| g - \chi \|^2_{L^2_{\nu^{(\alpha,\beta)}}} = \sum^{N}_{n=n_0+1} I_n + \sum^{n_0-1}_{n=1}  I_n + I_{n_0} = S_{1} + S_{-1} + S_{0}.
}
We shall address each term separately.  Consider $S_1$.  Since $\nu^{(\alpha,\beta)}(t) \lesssim (1-t)^{\alpha}$ on $[0,1]$, we have 
\bes{
I_n \lesssim \int_{V_n} | g(t) - g(t_n) |^2 (1-t)^{\alpha} \D t,\quad n=n_0,\ldots,N.
}
We now consider three cases: (i) $-1 < \alpha < 0$, (ii) $\alpha = 0$ and (iii) $\alpha > 0$.  Consider case (i).  Then
\eas{
I_n \lesssim \int_{V_n} \left| \int_{V_n} |g'(s)| \D s \right |^2 (1-t)^{\alpha} \D t
 \lesssim \int_{V_n} (1-t)^{\alpha} \D t \int_{V_n} (1-t)^{-\alpha-1} \D t \int_{V_n} |g'(t)|^2 (1-t)^{\alpha+1} \D t
}
By construction, $V_n$ is of width at most $2h$.  Hence, after a short calculation we get that
\bes{
I_n  \lesssim h \int_{V_n} |g'(t)|^2 (1-t)^{\alpha+1} \D t,\qquad -1 < \alpha < 0.
}
Now consider case (ii).  We have
\bes{
I_n \lesssim \int_{V_n} \D t \left| \int_{V_n} |g'(s)| \D s \right |^2 \leq \left ( \int_{V_n} \D t \right )^2 \int_{V_n} |g'(t)|^2 \D t \lesssim h^2 \int_{V_n} | g'(t) |^2 \D t,\qquad \alpha = 0.
}
Final, consider case (iii).  By similar arguments, we get that
\bes{
I_n \lesssim \left ( \int_{V_n} (1-t)^{\alpha} \int^{t}_{t_n} (1-s)^{-\alpha-1} \D s \D t \right )  \int_{V_n} |g'(t)|^2 (1-t)^{\alpha+1} \D t.
}
Write $V_n = (a,b)$ where $0 \leq a \leq t_n$ and $t_n \leq b \leq 1$.  Then
\eas{
\int_{V_n} (1-t)^{\alpha} \int^{t}_{t_n} (1-s)^{-\alpha-1} \D s \D t  &= \frac{1}{\alpha} \int^{b}_{a} (1-t)^{\alpha} \left ( (1-t)^{-\alpha} - (1-t_n)^{-\alpha} \right ) \D t
\\
& = \frac{1}{\alpha} \left [ (b-a) + \frac{1}{\alpha+1} \frac{(1-b)^{\alpha+1} - (1-a)^{\alpha+1}}{(1-t_n)^{\alpha}} \right ]
}
Note that $1-b \leq 1-t_n$ and $1-a \geq 1-t_n$.  Therefore 
\bes{
\int_{V_n} (1-t)^{\alpha} \int^{t}_{t_n} (1-s)^{-\alpha-1} \D s \D t \leq \frac{1}{\alpha} \left [ b-a - \frac{b-a}{\alpha+1} \right ] \lesssim h.
}
Hence
\bes{
I_n \lesssim h  \int_{V_n} |g'(t)|^2 (1-t)^{\alpha+1} \D t,\qquad \alpha > 0.
}
With these estimates to hand, we now deduce the following bound for the term $S_1$ in \R{int_sum_split}:
\be{
\label{S1_bd}
S_1 \lesssim \left \{ \begin{array}{lc} h^2 \| g' \|^2_{L^2(0,1)} & \alpha = 0 \\ h \| g' \|^2_{L^2_{\nu^{(\alpha+1,\beta+1)}}} & \alpha \neq 0  \end{array}  \right .,
}
This follows from the fact that the $V_n$ form a partition of $(-1,1)$, the definition of $n_0$ and the fact that $(1+t)^{\beta+1} \geq 1$ for $t \in [0,1]$.  Identical arguments give a similar result for $S_{-1}$:
\be{
\label{Sm1_bd}
S_{-1} \lesssim \left \{ \begin{array}{lc} h^2 \| g' \|^2_{L^2(-1,0)} & \beta = 0 \\ h \| g' \|^2_{L^2_{\nu^{(\alpha+1,\beta+1)}}} & \beta \neq 0  \end{array}  \right . .
}
Finally, consider $S_0$.  Let $V_{n_0} = J_{-1} \cup J_{1}$, where $J_{1} \subseteq [0,1]$ and $J_{-1} \subseteq [-1,0]$.  If $h M^2 \lesssim 1$ we may assume that $h \leq 1/2$ so that $|1-t| \geq 1/2$ and $|1+t| \geq 1/2$ for $t \in V_{n_0}$.   Then we have
\bes{
J_{\pm 1} \lesssim \left ( \int_{J_{\pm 1}} \D t \right )^2 \int_{J_{\pm 1}} |g'(t)|^2 \D t \lesssim h^2 \int_{J_{\pm 1}} |g'(t)|^2 \D t \lesssim h^2 \int_{J_{\pm 1}} |g'(t)|^2 \nu^{(\gamma,\delta)}(t)\D t,
}
for any $-1 < \gamma,\delta \leq 0$.  It now follows that $J_{\pm1}$ satisfy exactly the same bounds as \R{S1_bd} and \R{Sm1_bd} for $S_{\pm 1}$.  Therefore, in order to estimate the left-hand side of \R{int_sum_split} it suffices from now on to consider only $S_{\pm 1}$.  For this, we shall use Markov's inequality:
\be{
\label{Markov}
\| g' \|_{L^2(I)} \lesssim M^2/|I| \| g \|_{L^2(I)},\quad \forall g \in \bbP_{M},
}
where $I$ is an arbitrary bounded interval, as well as the following Markov-type inequality:
\be{
\label{General_Markov}
\| g' \|_{L^2_{\nu^{(\alpha+1,\beta+1)}}} \lesssim M \| g \|_{L^2_{\nu^{(\alpha,\beta)}}},\quad \forall g \in \bbP_{M}.
}
Markov's inequality \R{Markov} is well-known.  A proof of \R{General_Markov} is given in Appendix \ref{a:Jacobi}.

There are now four cases: (a) $\alpha = \beta = 0$, (b) $\alpha = 0$, $\beta \neq 0 $, (c) $\alpha \neq 0 $, $\beta = 0$ and (d) $\alpha \neq 0$, $\beta \neq 0$.  Consider case (a).  Then by \R{int_sum_split}, \R{S1_bd}, \R{Sm1_bd} and \R{Markov} we find that $\| g - \chi \|^2_{L^2} \lesssim h^2 M^4 \| g \|^2_{L^2}$.  Since $h M^2 \lesssim 1$, this now gives \R{desired_result} for case (a).  Now consider case (b).  By \R{int_sum_split}, \R{S1_bd}, \R{Sm1_bd}, \R{Markov} and \R{General_Markov} we get that
\bes{
\| g - \chi \|^2_{L^2_{\nu^{(0,\beta)}}} \lesssim h^2 M^4 \| g \|^2_{L^2(0,1)} + h M^2 \| g \|^2_{L^2_{\nu^{(0,\beta)}}} \lesssim h M^2 \| g \|^2_{L^2_{\nu^{(0,\beta)}}}.
}
Here in the final inequality we use the facts that $h M^2 \lesssim 1$ and $(1+t)^{\beta} \geq 1$ for $t \in [0,1]$.  Thus we get \R{desired_result} in this case as well.  Case (c) is near-identical to case (b).  To complete the proof, we consider case (d).  Using \R{int_sum_split}, \R{S1_bd}, \R{Sm1_bd} and \R{General_Markov}, we get 
\bes{
\| g - \chi \|^2_{L^2_{\nu^{(\alpha,\beta)}}} \lesssim h M^2 \| g \|^2_{L^2_{\nu^{(0,\beta)}}},
}
which yields \R{desired_result}.
}

\lem{
\label{l:Einf_Leg}
For $\alpha,\beta > -1$ let $\{ \phi_i \}_{i \in \bbN}$ be the orthonormal Jacobi polynomial basis \R{Jacobi_ON}, $T = \{ t_n \}^{N}_{n=1} \subseteq D$ be a set of scattered data points and suppose that $h$ is as in \R{h_def}.  If $h M^2 \leq 1$ then
\bes{
E_{\infty}(h,M) \lesssim h M^2 \log M,
}
where $E_{\infty}(h,M)$ is as in \R{E_2NM}.
}

\prf{
Let $x \in P_M(\ell^2(\bbN))$, $\| x \|_{\infty} = 1$ and set $g = \sum^{M}_{j=1} x_j \phi_j \in \bbP_{M-1}$ as in the previous proof.  Then
\be{
\label{unif_norm_identity}
\| (P_M - P_M U^* U P_M )x \|_{\infty} = \max_{i=1,\ldots,M} \left | \ip{g}{\phi_i}_{L^2_{\nu^{(\alpha,\beta)}}} - \ip{g}{\phi_i}_h \right |.
}
Observe that
\eas{
 \left | \ip{g}{\phi_i}_{L^2_{\nu^{(\alpha,\beta)}}} - \ip{g}{\phi_i}_h \right | &= \left | \sum^{N}_{n=1} \int_{V_n} \left ( g(t) \phi_i(t) - g(t_n) \phi_i(t_n) \right ) \nu^{(\alpha,\beta)}(t) \D t \right |
 \\
 & \leq \sum^{N}_{n=1} \int_{V_n} \left | \int_{t_n}^t \left ( g(s) \phi_i(s) \right )' \D s \right | \nu^{(\alpha,\beta)}(t) \D t 
 \\
& \lesssim \sum^{N}_{n=1} \int_{V_n} \nu^{(\alpha,\beta)}(t) \D t \int_{V_n} | (g(s) \phi_i(s))' | \D s .
}
Since $\| x \|_{\infty} = 1$, we have $| (g(s) \phi_i(s))' | \leq \sum^{M}_{j=1} | ( \phi_i(s) \phi_j(s))' |$ and we now substitute into \R{unif_norm_identity} to deduce that
\ea{
\label{unif_norm_ineq}
\| (P_M - P_M U^* U P_M )x \|_{\infty} &\lesssim \max_{i=1,\ldots,M} \sum^{M}_{j=1} \sum^{N}_{n=1}  \int_{V_n} \nu^{(\alpha,\beta)}(t) \D t  \int_{V_n} | (\phi_i(s) \phi_j(s))' | \D s  \nn
\\
& = \max_{i=1,\ldots,M} \sum^{M}_{j=1} \sum^{N}_{n=1} I_n  = \max_{i=1,\ldots,M} \sum^{M}_{j=1} \left ( S_{1}+S_{-1}+S_0 \right ),
}
where, as in the previous lemma, $S_1$ corresponds to $[0,1]$, $S_{-1}$ corresponds to $[-1,0]$ and $S_0$ corresponds to the term $I_{n_0}$ where $ 0 \in V_{n_0}$.

Consider the term $S_1$.  By \R{h_asymp} and \R{Jacobi_Local}, we have
\bes{
| \phi^{(k)}_i(t) | \lesssim \min \left \{ (\sqrt{1-t^2})^{-\alpha-k-1/2} i^{k} , i^{2k+\alpha+1/2} \right \},\quad 0 \leq t \leq 1.
}
Since $1 \leq i,j \leq M$, we have that
\bes{
S_1 \lesssim \sum^{N}_{n=n_0+1}  \int_{V_n} (1-t)^{\alpha} \D t \int_{V_n}\min \left \{ (\sqrt{1-t^2})^{-2\alpha-2} M , M^{2\alpha+3} \right \} \D t.
}
Let $n_0+1 \leq n^* < N$ be arbitrary (its value will be chosen later) and split this sum into two according to $n^*$.  Then
\ea{
S_1 &\lesssim M \sum^{n^*}_{n=n_0+1}  \int_{V_n} (1-t)^{\alpha} \D t  \int_{V_n} (\sqrt{1-t^2})^{-2\alpha-2} \D t + M^{2\alpha+3} \sum^{N}_{n=n^*+1}  \int_{V_n} (1-t)^{\alpha} \D t \int_{V_n}1 \D t \nn
\\
& = M S^{-}_{1} + M^{2\alpha+3} S^{+}_{1}. \label{S1}
}
We consider $S^{\pm}_1$ separately.  For $S^{+}_{1}$, we have
\be{
\label{Splus1}
S^{+}_{1} \lesssim h \int^{1}_{z^*} (1-t)^{\alpha} \D t \lesssim h (1-z^*)^{\alpha+1},
}
where $z^*$ is the right endpoint of $V_{n^*}$.  Now consider $S^{-}_{1}$:
\bes{
S^{-}_{1} \lesssim h \sum^{n^*}_{n=n_0+1}  \int_{V_n} (1-t)^{\alpha} \D t \sup_{t \in V_n}(1-t)^{-\alpha-1} = h \sum^{n^*}_{n=n_0+1}  \int_{V_n} (1-t)^{\alpha} (1-z_n)^{-\alpha-1} \D t.
}
where $z_n$ is the right endpoint of $V_n$.  However, if $t \in V_n$ is arbitrary then $z_n \leq t + h$.  Thus $(1-z_n)^{-\alpha-1} \leq (1-t-h)^{-\alpha-1}$, and this gives
\bes{
S^{-}_{1} \lesssim h \int^{z^*}_{0} (1-t)^{\alpha}(1-t-h)^{-\alpha-1} \D t \leq h \int^{1}_{1-z^*} s^{-1} (1-h/s)^{-\alpha-1} \D s 
}
Suppose now that $n^*$ is chosen so that
\be{
\label{nstar_choice}
1 - 2/M^2 - 2 h \leq z^* \leq 1 - 2/M^2,
}
(recall that the Voronoi cells are of width at most $2h$, hence such a choice is possible).
Then, since $h M^2 \leq 1$, we find that $1-h/s \gtrsim 1$ for $1-z^* \leq s \leq 1$.  Therefore, if \R{nstar_choice} holds we have
\be{
\label{Sminus1}
S^{-}_{1} \lesssim h \int^{1}_{1-z^*} s^{-1} \D s \lesssim h \log(1-z^*) \lesssim h \log (M).
}
Combining this with \R{Splus1} and \R{nstar_choice}, and using the fact that $h M^2 \leq 1$ once more, we now get that $S^{+}_{1} \lesssim h (M^{-2} + h )^{\alpha+1} \lesssim h M^{-2\alpha-2}$.
Substituting this and \R{Sminus1} back into \R{S1} now gives
\be{
\label{S1bound}
S_1 \lesssim S^{-}_{1} \lesssim h M \log(M) + h M \lesssim h M \log(M),
}
which completes the estimate for $S_{1}$.  The estimate for $S_{-1}$ is near-identical, except that we use \R{Jacobi_negative} as well as \R{Jacobi_Local} since $S_1$ sums over integrals contained in the negative portion of the interval.  Hence we get
\be{
\label{Sm1bound}
S_{-1} \lesssim h M \log(M),
}
for this term as well.  Next we need to estimate
\bes{
S_0 = \int_{V_{n_0}} \nu^{(\alpha,\beta)}(t) \D t \int_{V_{n_0}} | (\phi_i(s) \phi_j(s))' | \D s
}
Since $V_{n_0} \subseteq [-2h,2h]$, we have that $\nu^{(\alpha,\beta)}(t)  \lesssim 1$ for $t \in V_{n_0}$.  Also, by \R{Jacobi_Local} and \R{Jacobi_negative}, we have $| (\phi_i(s) \phi_j(s))' | \lesssim M$, $s \in V_{n_0}$.  Hence we get $S_0 \lesssim h M$.
Combining this with \R{S1bound} and \R{Sm1bound} and substituting into \R{unif_norm_ineq} now gives
\bes{
\| (P_M - P_M U^* U P_M )x \|_{\infty} \lesssim \sum^{M}_{j=1} h M \log M = h M^2 \log M,
}
from which the result follows immediately.
}

\lem{
\label{l:F_N_M_R_poly}
For $\alpha,\beta > -1$ let $\{ \phi_i \}_{i \in \bbN}$ be the orthonormal Jacobi polynomial basis \R{Jacobi_ON}, $T= \{ t_n \}^{N}_{n=1} \subseteq D$ be a set of $N$ scattered data points and suppose that $h$ is as in \R{h_def}.  Suppose that $h M^2 \leq 1$.  If the weights $w_i \gtrsim i^{q+1/2}$, where $q$ is as in \R{Jacobi_ON_growth2}, then the quantity $F(h,M,R)$ defined by \R{F_N_M_R} satisfies
\be{
\label{F_poly_gen}
F(h,M,R) \lesssim \sqrt{M} \sup_{i > R} \left \{ i^{q+1/2} / w_i  \right \},
}
Moreover, if the weights $w_i \gtrsim i \log i$ and $R \geq M$ then
\be{
\label{F_poly_log}
F(h,M,R) \lesssim h M \sup_{i > R} \left \{ \frac{i \log i}{w_i} \right \}.
}
}
\prf{
As before, let $x \in P_M(\ell^2(\bbN))$, $\| x \|_{\infty} = 1$ and set $g = \sum^{M}_{j=1} x_j \phi_j \in \bbP_{M-1}$.  Then
\bes{
\| P^{\perp}_R W^{-1} U^* U P_M x \|_{\infty} = \sup_{i > R} \left | \frac{1}{w_i} \ip{g}{\phi_i}_h \right | \leq \| g \|_h \sup_{i > R} \left \{ \frac{\| \phi_i \|_h}{w_i} \right \}.
}
Note that $\| \phi_i \|_h \leq \| \phi_i \|_{L^\infty} \lesssim i^{q+1/2}$ by \R{Jacobi_ON_growth} and also
\eas{
\| g \|^2_h &= \ip{P_M U^* U P_M x}{x}
 \leq \left ( 1+\| P_M - P_M U^* U P_M \| \right ) \| x \|^2
 \lesssim M,
}
where in the final inequality we use Lemma \ref{l:E2_Leg} and the fact that $\| x \| \leq \sqrt{M} \| x \|_{\infty} = \sqrt{M}$.  This now gives \R{F_poly_gen}.

For \R{F_poly_log} we use orthogonality and the fact that $R \geq M$ to get
\be{
\label{FNMR_bison}
\| P^{\perp}_R W^{-1} U^* U P_M x \|_{\infty} = \sup_{i > R} \left | \frac{1}{w_i} \ip{g}{\phi_i}_h \right | = \sup_{i > R} \left \{ \frac{1}{w_i} \left | \ip{g}{\phi_i}_{L^2_{\nu}} -   \ip{g}{\phi_i}_h \right | \right \},
}
We now proceed in a similar manner to the proof of Lemma \ref{l:Einf_Leg}.  First, since $\| x \|_{\infty}=1$ we have
\bes{
\left | \ip{g}{\phi_i}_{L^2_{\nu}} -   \ip{g}{\phi_i}_h \right | \leq \sum^{M}_{j=1} \sum^{N}_{n=1} \int_{V_n} \nu^{(\alpha,\beta)}(t) \D t \int_{V_n} | (\phi_i(s) \phi_j(s))' | \D s.
}
We now argue in a similar way, using the fact that $j \leq M \leq R \leq i$.  This gives
\bes{
\left | \ip{g}{\phi_i}_{L^2_{\nu}} -   \ip{g}{\phi_i}_h \right | \lesssim \sum^{M}_{j=1} h i \log(i) \lesssim M h i \log(i).
}
Substituting back into \R{FNMR_bison} now gives the required result.
}

We are now ready to prove Theorem \ref{t:Leg_poly_full_l1}:

\prf{[Proof of Theorem \ref{t:Leg_poly_full_l1}]
We use Theorem \ref{t:full_samp_recov} and the estimates of Lemmas \ref{l:E2_Leg}--\ref{l:F_N_M_R_poly}.  Note that
\bes{
\min_{M < i \leq R } \{ w_i \},\ \max_{i=1,\ldots,M} \{ w_i \} \asymp M^{\gamma+q+1/2},\quad M \rightarrow \infty.
}
Hence, we require $h$, $M$ and $R \geq M$ such that 
\bes{
E(h,M) < \epsilon,\quad E(h,R) \lesssim \epsilon,\quad F(h,M,R) \lesssim \epsilon M^{-\gamma-q-1/2}.
}  
Since $R \geq M$, Lemmas \ref{l:E2_Leg} and \ref{l:Einf_Leg} give that the first two conditions are satisfied provided
\bes{
h \lesssim \frac{\epsilon}{R^2 \log R}.
}
We now consider $F(h,M,R)$.  Suppose first that $0 < \gamma \leq 1/2-q$.  Then Lemma \ref{l:F_N_M_R_poly} gives that $F(N,M,R) \lesssim \epsilon M^{-\gamma-q-1/2}$ provided $R \gtrsim \epsilon^{-1/\gamma} M^{1+(q+1)/\gamma}$.
Hence this results in the condition
\be{
\label{hN_argh}
h \lesssim \frac{c(\epsilon)}{M^{2+2(q+1)/\gamma} \log M},
}
which gives the result for $0 < \gamma \leq 1/2-q$.  Now suppose $\gamma > 1/2-q$.  Then $w_i \gtrsim i^{\gamma+q+1/2} \gtrsim i \log i$ and therefore Lemma \ref{l:F_N_M_R_poly} gives $F(N,M,R) < \epsilon M^{-\gamma-q - 1/2}$ whenever
\bes{
h \frac{\log R}{R^{\gamma+q-1/2}} < \epsilon M^{-3/2 - \gamma-q}.
}
Suppose that $R = c M$ for some $c$.  Then the above condition holds, provided
\bes{
h\lesssim \frac{\epsilon}{M^2 \log M}.
}
Moreover, substituting $R = c M$ into \R{hN_argh} gives precisely the same condition on $h$ in terms of $M$.  Hence the result follows.
}

The proof of Theorem \ref{t:Jacobi_poly_LS} is straightforward:
\prf{[Proof of Theorem \ref{t:Jacobi_poly_LS}]
We use Theorem \ref{t:LS_err} in combination with Lemma \ref{l:E2_Leg}.
}

Finally, we also give the proof of Theorem \ref{t:Leg_trunc}:
\prf{[Proof of Theorem \ref{t:Leg_trunc}]
We shall use Lemma \ref{l:gamma}.  Let $y \in \bbC^N$, $\| y \| =1$ be given.  Let $\chi \in C^{\infty}_{c}(-1,1)$ be a smooth compactly-supported function in $(-1,1)$ with the properties
\bes{
\| \chi \| = 1,\quad \chi(0) = 1,\quad 0 \leq \chi(x) \leq 1,\ x \in (-1,1).
}
Define
\bes{
g(t) = \sum^{N}_{n=1} \frac{y_{n}}{\sqrt{\tau_n}} g_{n}(t),\quad g_{n}(t) = \chi \left ( \frac{x-t_n}{\xi_{n}} \right ),
}
where $\xi_{n} = \mathrm{dist}(t_n,\partial V_n)$ is the distance of the point $t_n$ from the boundary of its Voronoi cell (in this one-dimensional setting, $V_n$ is an interval and its boundary is the set of the two endpoints).  Observe that
\bes{
\xi_{n} = \frac12 \min \left \{ t_{n+1}-t_n , t_n - t_{n-1} \right \},\quad n=1,\ldots,N,
}
and therefore $\xi_{n} \geq \xi$ for each $n=1,\ldots,N$.  By construction $\mathrm{supp}(g_{n}) \subseteq V_n$, $n=1,\ldots,N$, and therefore $\mathrm{supp}(g_{n}) \cap \mathrm{supp}(g_{m}) = 0$, $n \neq m$.  It follows that $g \in G_y$.  Since $g \in C^{\infty}[-1,1]$ a standard result in polynomial approximation gives that
\bes{
\inf_{\substack{\phi \in \Phi_K \\ \phi \neq 0}} \| g - \phi \|_{\infty} \leq C_r K^{-r} \| g^{(r)} \|_{\infty},\quad K \geq r,
}
for some constant $C_r > 0$ independent of $K$ and $g$ (see \cite[(5.4.16)]{SMSD}).  Observe that
\bes{
\| g^{(r)} \|_{\infty} = \max_{n=1,\ldots,N} \frac{|y_{n}|}{\sqrt{\tau_n}} \| g^{(r)}_{n} \|_{\infty} \leq \| y \| \max_{n=1,\ldots,N} \frac{(\xi_{n})^{-r}}{\sqrt{\tau_n}} \leq \xi^{-r-1/2},
}
where in the last inequality we use that fact that $\xi_{n} \leq \tau_n$.  Hence
\be{
\label{interp_poly_err}
\inf_{\substack{\phi \in \Phi_K \\ \phi \neq 0}} \| g - \phi \|_{\infty} \leq C_r K^{-r} \xi^{-r-1/2},\quad K \geq r.
}
In order to apply \R{gamma_def}, it remains to estimate $\| g \|_{\nu} = \| g \|$.  Since the $g_{n}$'s have disjoint supports, we have
\bes{
\| g \|^2 = \sum^{N}_{n=1} \frac{|y_{n}|^2}{\tau_n} \| g_{n} \|^2 = \sum^{N}_{n=1} \frac{|y_{n}|^2}{\tau_n} \xi_{n} \geq \min_{n=1,\ldots,N} \{ \xi_{n} / \tau_n \} \geq \xi_{N} / \max_{n=1,\ldots,N} \tau_n.
}
Observe that $\tau_n \leq 2 h$.  Therefore, $\| g \| \geq \sqrt{\xi / (2 h)}$.  Substituting this and \R{interp_poly_err} into \R{gamma_def} now gives the result.  
}

Finally, we now note that Theorem \ref{t:intro_thm} follows immediately from Theorems \ref{t:Leg_trunc} and \ref{t:Leg_poly_full_l1}.

\section{Trigonometric polynomials on bounded intervals}\label{s:TrigPoly_Examp}
We now consider Example \ref{ex:Trigonometric}.  Note that in this case we define the projections $P_N : \ell^2(\bbZ) \rightarrow \ell^2(\bbZ)$ by $P_N x = \{\ldots,0,0,x_{-N},x_{-N+1},\ldots,x_{N-1},0,0,\ldots\}$.  Our main result is as follows:

\thm{
\label{t:Trig_Scatter}
Let $\{ \phi_i\}_{i \in \bbZ}$ be the Fourier basis \R{fourier_basis}, $T = \{ t_n \}^{N}_{n=1} \subseteq [-1,1]$ be a set of $N$ scattered data points and suppose that $h$ is as in \R{h_def}.  Suppose that the weights $w_{i} = 1+ |i|^{\gamma}$ for some $\gamma >0$.  Then for each $0 < \epsilon < 1/2$ there exists a $c(\epsilon) > 0$ such that if
\bes{
h \leq c(\epsilon) \left \{ \begin{array}{cc} M^{-3/2-3/(4\gamma)} & 0 < \gamma <1  \\  M^{-3/2} & \gamma > 1 \end{array} \right . ,
}
then any minimizer $\hat{x}$ of \R{fin_min} satisfies
\bes{
\| x - \hat{x} \| \leq C(\epsilon) \left [ M^{\gamma+1/2} \eta + \| P^{\perp}_M x \|_{1,w} + T_{N,K, \eta}(x) \right ).
}
for some constant $C$ depending on $\epsilon$ only, where $T_{N,K, \eta}(x)$ is as in \R{trunc}.
}

For least-squares fitting, we have the following:

\thm{
\label{t:Trig_LS}
Let $\{ \phi_i\}_{i \in \bbZ}$ be the Fourier basis \R{fourier_basis}, $T = \{ t_n \}^{N}_{n=1} \subseteq [-1,1]$ be a set of $N$ scattered data points and let $h$ be as in \R{h_def}.  Then for each $0 < \epsilon < 1$ there exists a $c(\epsilon) > 0$ such that if
\be{
\label{hcondTrigLS}
h \leq c(\epsilon) M^{-1},
}
then the solution $\check{x}$ of \R{LS_fit} exists uniquely and satisfies
\bes{
\| x - \check{x} \| \leq \left ( 1 + \frac{1}{\sqrt{1-\epsilon}} \right ) \| x - P_M x \|_{1,w} + \frac{1}{\sqrt{1-\epsilon}}  \eta,
}
for any $w = \{ w_i \}_{i \in \bbN}$ with $w_i \geq 1$.
}

Unlike the case of Jacobi polynomials, these results give a worse recovery guarantee for weighted $\ell^1$ minimization than that of least-squares fitting.  However, we do not believe the scaling $h \lesssim M^{-3/2}$ is sharp, and instead we conjecture that the true scaling is $h \lesssim (M \log M)^{-1}$.  Proving this conjecture is an open problem.  We remark in passing that this conjecture holds in the special case where the data is equispaced (we omit the proof for brevity's sake).

\subsection{Numerical examples}

\begin{figure}[t]
\begin{center}
$\begin{array}{ccc}
\includegraphics[width=4.5cm]{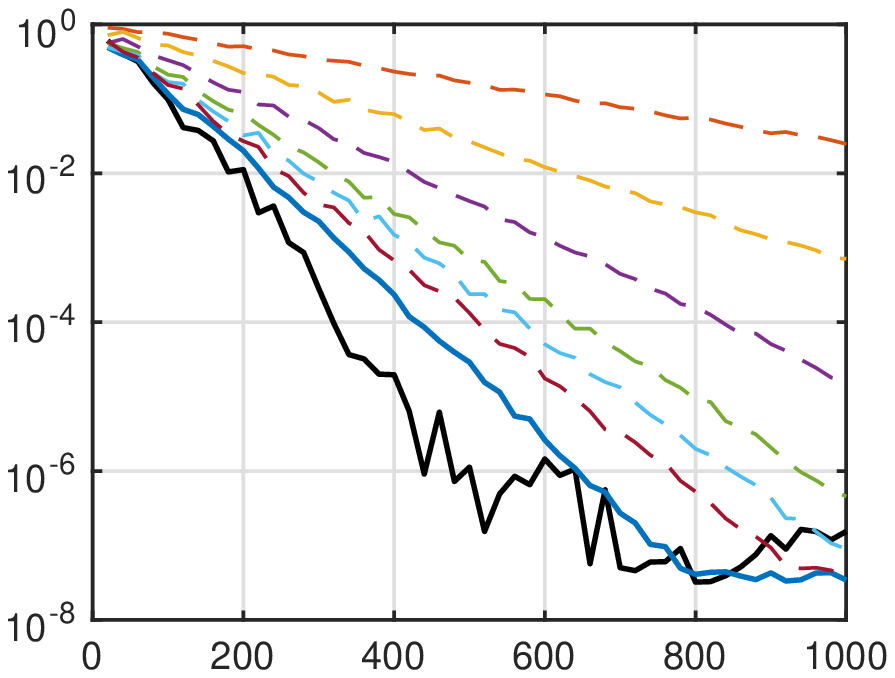}   
&    
\includegraphics[width=4.5cm]{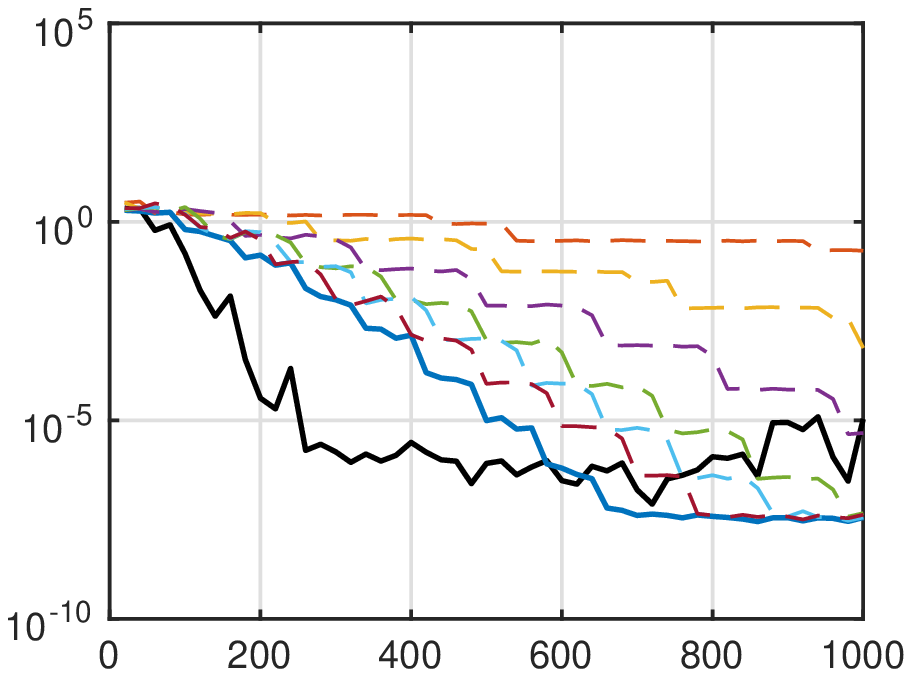}   
&
\includegraphics[width=4.5cm]{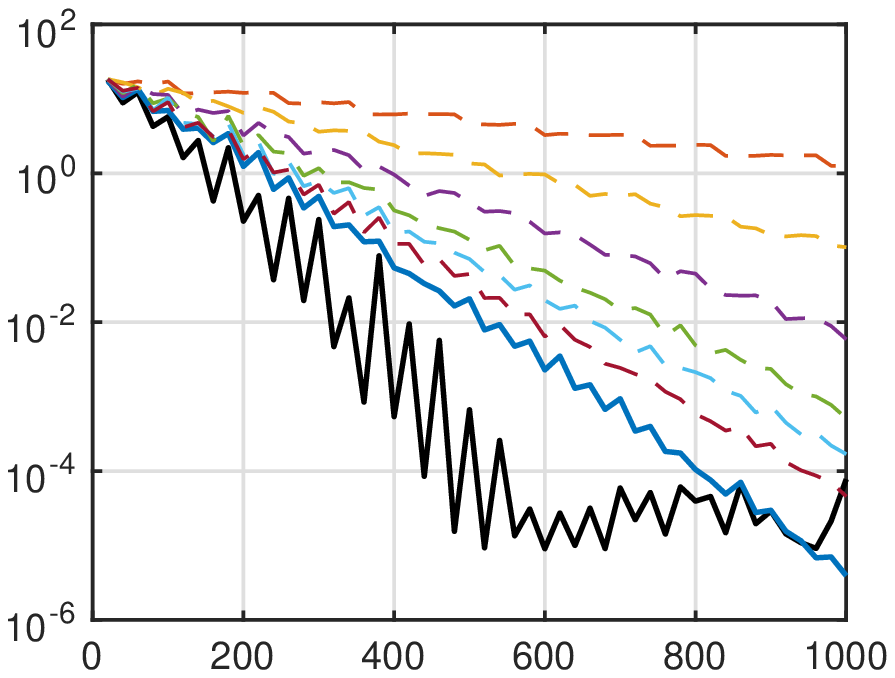} 
\\ 
f(x) = \frac{1}{1+500 \cos^2(\pi x)} & 
f(x) = \cos(4 \pi x) \exp(\sin(40 \pi x)) & f(x) = \frac{1}{20/19-\sin(10 \pi x)}
\end{array}$
\caption{
Numerical comparison of weighted $\ell^1$ minimization and least-squares fitting for approximation from jittered data.  The error against $N$ is plotted for each method.  The solid black line is weighted $\ell^1$ minimization with $K=4N$ and weights $w_{i} = \sqrt{i}$.  The dashed lines are least squares with $M = c N$ and $c=\frac16,\frac14,\frac12,\frac23,\frac34,\frac56$.  The solid blue line is oracle least squares based on choosing $M$ to minimize the error for a given $N$ and $f$.  Random noise of magnitude $10^{-8}$ was added to the data.  
}
\label{f:FourJitt}
\end{center}
\end{figure}

In Fig.\ \ref{f:FourJitt} we give a comparison of the two techniques for jittered data.  In alignment with the discussion above, these results suggest the scaling predicted by Theorem \ref{t:Trig_Scatter} is not optimal.  In fact, weighted $\ell^1$ minimization performs better than least-squares fitting (including the oracle case) in all the examples.  We suspect this strong performance is due in part to the presence of some sparsity in the functions considered, and the fact that jittered points are near-optimal points for the recovery of sparse trigonometric polynomials \cite{FoucartRauhutCSbook}.

\subsection{Proofs}
We first require the following lemma:

\lem{
Let $\{ \phi_i\}_{i \in \bbZ}$ be the orthonormal Fourier basis \R{fourier_basis}, $T = \{ t_n \}^{N}_{n=1} \subseteq D$ be a set of $N$ scattered data points and suppose that $h$ is as in \R{h_def}.  Suppose that $h M \leq 1$.  If $E_2(h,M)$ and $E_{\infty}(h,M)$ are as in \R{E_2NM} and \R{E_2NM} respectively, then
\bes{
E_{2}(h,M) \lesssim h M,\quad E_{\infty}(h,M) \lesssim h M^{3/2}.
}
}
\prf{
Consider $E_{2}(h,M)$ first.  As in the proof of Lemma \ref{l:E2_Leg}, let $x \in P_M(\ell^2(\bbZ))$, $\| x \|=1$ be arbitrary and set $g = \sum^{M-1}_{j=-M} x_j \phi_j$ so that $\| g \|_{L^2} = 1$.  Arguing in an identical manner, we see that it suffices to show that
\be{
\label{fourier_E2_desired}
\sum^{N}_{n=1} \int_{V_n} | g(t) - g(t_n) |^2 \D t \lesssim h^2 M^2.
}
Observe that $| g(t) -g(t_n) |^2 \lesssim h \int_{V_n} |g'(s) |^2 \D s$ and therefore $\sum^{N}_{n=1} \int_{V_n} | g(t) - g(t_n) |^2 \D t \lesssim h^2 \| g' \|^2_{L^2}$.
To get \R{fourier_E2_desired} we recall Bernstein's inequality $\| g' \|_{L^2} \lesssim M \| g \|_{L^2}$ for trigonometric polynomials.

For $E_{\infty}(h,M)$ we let $x \in P_M(\ell^2(\bbZ))$ with $\| x \|_{\infty}=1$ and defined $g$ as before.  As in the proof of Lemma \ref{l:Einf_Leg} it suffices to estimate
\bes{
\max_{i=-M,\ldots,M-1} \left | \ip{g}{\phi_i}_{L^2} - \ip{g}{\phi_i}_h \right |.
}
Arguing in the standard way, we see that
\bes{
\max_{i=-M,\ldots,M-1} \left | \ip{g}{\phi_i}_{L^2} - \ip{g}{\phi_i}_h \right | \lesssim h \| (g \phi_i)' \|_{L^1} \leq  h\| (g \phi_i)' \|_{L^2}
}
Since $g \phi_{i}$ is a trigonometric polynomial of degree at most $2M$ Bernstein's inequality gives
\bes{
\max_{i=-M,\ldots,M-1} \left | \ip{g}{\phi_i}_{L^2} - \ip{g}{\phi_i}_h \right | \lesssim h M \| x \|_{2} \leq h M^{3/2},
}
where in the final inequality we use the Cauchy--Schwarz inequality the fact that $\| x \|_{\infty} = 1$.  This now gives the estimate for $E_{\infty}(h,M)$.
}

\lem{
Let $\{ \phi_i\}_{i \in \bbZ}$ be the orthonormal Fourier basis \R{fourier_basis}, $T = \{ t_n \}^{N}_{n=1} \subseteq D$ be a set of $N$ scattered data points and suppose that $h$ is as in \R{h_def}.  Suppose that $h M \leq 1$.  Then the quantity $F(h,M,R)$ defined by \R{F_N_M_R} satisfies
\bes{
F(h,M,R) \lesssim \sqrt{M} \sup_{i > R} \left \{ 1/w_{i} \right \}.
}
Moreover, if the weights $w_i \gtrsim i$ and $R \geq M$, then
\bes{
F(h,M,R) \lesssim h \sqrt{M} \sup_{i>R} \left \{ i / w_i \right \}.
}
}
\prf{
We argue as in the proof of Lemma \ref{l:F_N_M_R_poly}.  In the first case, since the functions $\phi_i$ are uniformly bounded we have $\| P^{\perp}_R W^{-1} U^* U P_M x \|_{\infty} \leq \| g \|_{h} \sup_{i>R} \{ 1/w_{i} \}$,
where $g = \sum^{M-1}_{j=-M} x_j \phi_j$ and $\| x \|_{\infty} = 1$.  By the same argument, we find that $\| g \|_{h} \lesssim \| g \|_{L^2} \lesssim \sqrt{M}$, which gives the first result.

Now consider the second.  With $x$ and $g$ as before, we have
\bes{
\| P^{\perp}_R W^{-1} U^* U P_M x \|_{\infty} = \sup_{i > R} \left \{ \frac{1}{w_i} \left | \ip{g}{\phi_i}_{L^2} - \ip{g}{\phi_i}_h \right | \right \}.
}
As in the previous lemma, we note that $\left | \ip{g}{\phi_i}_{L^2} - \ip{g}{\phi_i}_h \right |  \leq h \| (g \phi_i)' \|_{L^2}$,
and therefore by Bernstein's inequality $\left | \ip{g}{\phi_i}_{L^2} - \ip{g}{\phi_i}_h \right | \lesssim h i \| x \|_{2} \lesssim h i \sqrt{M}$.
This gives the second result.
}

\prf{[Proof of Theorem \ref{t:Trig_Scatter}]
With this lemma in hand, the proof is identical in manner to that of Theorem \ref{t:Leg_poly_full_l1}.  We omit the details.
}

\section{Conclusions}\label{s:conclusions}
We have presented an infinite-dimensional framework for weighted $\ell^1$ minimization.  Its advantages are that it does not require \textit{a priori} knowledge of the expansion tail in order to be implemented and in the absence of noise it leads to interpolatory approximations.  We have discussed the role weights play in the minimization in resolving the aliasing phenomenon, as opposed to promoting smoothness, and provided an explicit way to choose the truncation parameter.  In the second half this paper we performed a linear error analysis for this framework valid for arbitrary scattered data, and used it to show near-optimal performance for Jacobi polynomial bases.

There are several topics for future research.  Three immediate problems are (i) to obtain a better scaling than \R{hcondTrigLS} in the trigonometric polynomial case, (ii) to estimate the truncation error $T_{h,K,\eta}(x)$ in a way that does not require additional regularity of $x$ (see Theorem \ref{p:trunc_err} and Remark \ref{r:trunc_weak}), and (iii) to improve the noise bound in Theorem \ref{t:Leg_poly_full_l1} (see Remark \ref{noise_growth}).  Besides these, a question of singular importance are the extensions of \S \ref{s:AlgPoly_Examp} to higher dimensions and to unbounded intervals (using Laguerre and Hermite polynomials, for example).  Other higher-dimensional problems can also be investigated, such as approximations in spherical harmonics (see also \cite{RauhutWardSpher}).

Another topic is the optimal selection of weights.  The results of this paper suggest that weights aid approximations from deterministic, scattered data by resolving the aliasing phenomenon and not necessarily by matching the decay of the expansion coefficients.  In particular, slowly growing weights seem sufficient, at least in the one-dimensional setting.  The situation may be different however in the multidimensional case when the samples are random.  It has recently been shown in \cite{AdcockCSFunInterp,ChkifaDownwardsCS} that weighted $\ell^1$ minimization with a specific choice of weights leads to optimal approximation rates for certain classes of multivariate functions when the samples are drawn randomly from the orthogonality measure of the polynomial basis.
We also note the possibility of using reweighted $\ell^1$ minimization, where weights iteratively updated to get a better estimation of the support set of $x$ \cite{PengHamptonDoostantweighted,KarniadakisUQCS}.  We expect this technique can be combined with our framework.

\section*{Acknowledgements}
The work was supported by the Alfred P. Sloan Foundation and the Natural Sciences and Engineering Research Council of Canada through grant 611675.  A preliminary version of this work was presented during the Research Cluster on ``Computational Challenges in Sparse and Redundant Representations'' at ICERM in November 2014.  The author would like to thank the participants for the useful feedback received during the program.  He would also like to thank Alireza Doostan, Anders Hansen, Rodrigo Platte, Aditya Viswanathan, Rachel Ward and Dongbin Xiu.

\appendix

\section{Jacobi polynomials}\label{a:Jacobi}
Given $\alpha,\beta > -1$ let $P^{(\alpha,\beta)}_j$ be the Jacobi polynomial of degree $j$.  These polynomials are orthogonal on $D=(-1,1)$ with respect to $\nu^{(\alpha,\beta)}(t) = (1-t)^\alpha (1+t)^{\beta}$,
with
\bes{
\ip{P^{(\alpha,\beta)}_j}{P^{(\alpha,\beta)}_k}_{L^2_{\nu^{(\alpha,\beta)}}} = \delta_{j,k} \kappa^{(\alpha,\beta)}_j,
}
where
\be{
\label{Jacnorm_def}
\kappa^{(\alpha,\beta)}_j = \frac{2^{\alpha+\beta+1}}{2j+\alpha+\beta+1} \frac{\Gamma(j+\alpha+1) \Gamma(j+\beta+1)}{j! \Gamma(j+\alpha+\beta+1)},
}
and have the normalization
\bes{
P^{(\alpha,\beta)}_{j}(1) = \left ( \begin{array}{c} j + \alpha \\ j \end{array} \right ).
}
The corresponding orthonormal polynomials are defined by $\phi_j(t) = \left ( \kappa^{(\alpha,\beta)}_{j-1} \right )^{-1/2} P^{(\alpha,\beta)}_{j-1}(t)$, $j \in \bbN$.  Note that
\be{
\label{h_asymp}
\kappa^{(\alpha,\beta)}_{j} \sim 2^{\alpha+\beta} j^{-1} ,\quad j \rightarrow \infty.
}
and also that
\bes{
P^{(\alpha,\beta)}_{j}(1) \sim \frac{j^{\alpha}}{\Gamma(\alpha+1)}.
}
The polynomials $P^{(\alpha,\beta)}_j$ satisfy the differential equation
\be{
\label{Jacobi_DE}
-\left ( \nu^{(\alpha+1,\beta+1)} \left ( P^{(\alpha,\beta)}_{j} \right )' \right )' + \lambda^{(\alpha,\beta)}_j \nu^{(\alpha,\beta)} P^{(\alpha,\beta)}_j = 0,
}
where $\lambda^{(\alpha,\beta)}_{j} = j(j+\alpha+\beta+1)$.  In particular, the derivatives $(P^{(\alpha,\beta)}_{j})'$ are orthogonal with respect to $\nu^{(\alpha+1,\beta+1)}$ and satisfy
\be{
\label{Jacobi_Derivative}
\left ( P^{(\alpha,\beta)}_{j} \right )' = \sqrt{\frac{\lambda^{(\alpha,\beta)}_j \kappa^{(\alpha,\beta)}_{j}}{\kappa^{(\alpha+1,\beta+1)}_{j-1}} } P^{(\alpha+1,\beta+1)}_{j-1}.
}

\lem{
Let $\alpha,\beta > -1$.  Then $\| p' \|_{L^2_{\nu^{(\alpha+1,\beta+1)}}} \leq \lambda^{(\alpha,\beta)}_M \| p \|_{L^2_{\nu^{(\alpha,\beta)}}}$, $\forall p \in \bbP_M$.
}
\prf{
Let $p \in \bbP_M$ be arbitrary and observe that
\bes{
p (t)= \sum^{M}_{j=0} \frac{x_j}{\kappa^{(\alpha,\beta)}_j} P^{(\alpha,\beta)}_j(t),\qquad x_j = \ip{p}{P^{(\alpha,\beta)}_j}_{L^2_{\nu^{(\alpha,\beta)}}}.
}
Note that
\be{
\label{Jacobi_Pars_x}
\| p \|^2_{L^2_{\nu^{(\alpha,\beta)}}} = \sum^{M}_{j=0} \frac{|x_j |^2}{\kappa^{(\alpha,\beta)}_j}.
}
Similarly,
\bes{
p'(t) = \sum^{M-1}_{j=0} \frac{y_j}{\kappa^{(\alpha+1,\beta+1)}_{j}} P^{(\alpha+1,\beta+1)}_{j}(t),\qquad y_j = \ip{p'}{P^{(\alpha+1,\beta+1)}_j}_{L^2_{\nu^{(\alpha+1,\beta+1)}}},
}
and
\be{
\label{Jacobi_Pars_y}
\| p' \|^2_{L^2_{\nu^{(\alpha+1,\beta+1)}}} = \sum^{M-1}_{j=0} \frac{|y_j |^2}{\kappa^{(\alpha+1,\beta+1)}_j}.
}
Consider $x_j$.  By the differential equation \R{Jacobi_DE} and the fact that $\nu^{(\alpha+1,\beta+1)}(\pm 1 ) = 0$, we have
\eas{
x_j = \int^{1}_{-1} p(t) P^{(\alpha,\beta)}_j(t) \nu^{(\alpha,\beta)}(t) \D t 
 = \frac{1}{\lambda^{(\alpha,\beta)}_j} \int^{1}_{-1} p'(t) \left ( P^{(\alpha,\beta)}_{j}(t) \right )' \nu^{(\alpha+1,\beta+1)}(t) \D t.
}
Hence, by \R{Jacobi_Derivative},
\eas{
x_j &= \sqrt{\frac{\kappa^{(\alpha,\beta)}_{j}}{\lambda^{(\alpha,\beta)}_j \kappa^{(\alpha+1,\beta+1)}_{j-1}} } \ip{p'}{P^{(\alpha+1,\beta+1)}_{j-1}}_{L^2_{\nu^{(\alpha+1,\beta+1)}}} = \sqrt{\frac{\kappa^{(\alpha,\beta)}_{j}}{\lambda^{(\alpha,\beta)}_j \kappa^{(\alpha+1,\beta+1)}_{j-1}} } y_{j-1}.
} 
Using \R{Jacobi_Pars_x} and \R{Jacobi_Pars_y} we now get that
\bes{
\| p \|^2_{L^2_{\nu^{(\alpha,\beta)}}} \geq \sum^{M}_{j=1} \frac{|y_{j-1}|^2}{ \lambda^{(\alpha,\beta)}_j \kappa^{(\alpha+1,\beta+1)}_{j-1}} \geq \frac{1}{\lambda^{(\alpha,\beta)}_{M} } \| p' \|^2_{L^2_{\nu^{(\alpha+1,\beta+1)}}},
}
as required.
}

We also require several results concerning the asymptotic behaviour of Jacobi polynomials.  The first is as follows (see \cite[Thm.\ 7.32.1]{SzegoOrthPolys}):
\bes{
\| P^{(\alpha,\beta)}_j \|_{L^\infty} = \ord{j^{q}} ,\quad n \rightarrow \infty,\qquad 
q = \max \{ \alpha,\beta , -1/2 \} .
}
Hence, using \R{h_asymp} we find that the normalized functions $\phi_j$ defined by \R{Jacobi_ON} satisfy
\be{
\label{Jacobi_unif}
\| \phi_j \|_{L^\infty} = \ord{j^{q+1/2}},\quad j \rightarrow \infty,
}
which gives \R{Jacobi_ON_growth}.  We also note the following local estimates for Jacobi polynomials.  If $k=0,1,2,\ldots$ and $c >0$ is a fixed constant then
\be{
\label{Jacobi_Local}
\left | \frac{\D^kP^{(\alpha,\beta)}_j(t)}{\D t^k} \Bigg |_{t = \cos \theta} \right |= \left \{ \begin{array}{ll} \theta^{-\alpha-k-1/2} \ord{j^{k-1/2}} & c j^{-1} \leq \theta \leq \pi/2 \\ \ord{j^{2k+\alpha}} & 0 \leq \theta \leq c j^{-1} \end{array} \right . ,
}
as $j \rightarrow \infty$.  See \cite[Thm.\ 7.32.4]{SzegoOrthPolys}.  This estimate bounds the Jacobi polynomial and its derivatives for $0 \leq t \leq 1$.  For negative $t$, we may use the relation
\be{
\label{Jacobi_negative}
P^{(\alpha,\beta)}_j(-t) = (-1)^j P^{(\beta,\alpha)}_j(t).
}
Hence behaviour of $P^{(\alpha,\beta)}_j(t)$ and its derivatives for $t<0$ is given by \R{Jacobi_Local} with $\alpha$ replaced by $\beta$.

\bibliographystyle{abbrv}
\small
\bibliography{l1_pointwise_refs}

\end{document}